\documentclass[11pt,reqno]{amsart}
\usepackage{amsmath, amssymb, amsthm}
\usepackage[linktocpage=true]{hyperref}
\hypersetup{colorlinks,linkcolor=blue,urlcolor=blue,citecolor=blue}
\usepackage{mathtools} %\underbracket
\usepackage[table]{xcolor}
\usepackage{stmaryrd} %\llbracket, \rrbracket
\usepackage{bbm} % \mathbbm
\usepackage{leftidx} %\ltrans, \leftidx
\usepackage{mathabx} %\widecheck
\usepackage{rotating} %\rotatebox
\usepackage{wasysym} %\lhd
\usepackage[mmddyy]{datetime}
\usepackage{qtree}
\usepackage[all]{xy}
\usepackage[bbgreekl]{mathbbol} %\bbm for Greek letters
\DeclareSymbolFontAlphabet{\mathbbm}{bbold}
\DeclareSymbolFontAlphabet{\mathbb}{AMSb}
\usepackage{arydshln}
\makeatletter
\renewcommand*\env@matrix[1][*\c@MaxMatrixCols c]{%
  \hskip -\arraycolsep
  \let\@ifnextchar\new@ifnextchar
  \array{#1}}
\makeatother
\usepackage{youngtab,young}
\makeatletter
\newcommand{\labeltarget}[1]{\Hy@raisedlink{\hypertarget{#1}{}}}
\makeatother
\usepackage{todonotes}
\usepackage{pgf, tikz, colortbl}
\usetikzlibrary{cd, braids}
\numberwithin{equation}{subsection}
\usepackage[margin = 0.9in]{geometry}
\newcommand{\tcH}{\widetilde{\cH}}

\def\<{\langle}
\def\>{\rangle}
\DeclareMathSymbol{\hp}{\mathord}{AMSa}{"39}
\newcommand{\ba}{\mbf a}
\newcommand{\bA}{\mbf A}
\newcommand{\bc}{\mbf c}

\newcommand{\bH}{\mbf H}

\newcommand{\bS}{\mbf S}
\newcommand{\bT}{\mbf T}

\newcommand{\bbC}{\mathbb{C}}

\newcommand{\bbK}{\mathbb{K}}
\newcommand{\bbq}{\mathbbm{q}}

\newcommand{\bbH}{\mathbb{H}}

\newcommand{\D}{\mrm{D}}

\newcommand{\cA}{\mathcal{A}}
\newcommand{\cB}{\mathcal{B}}
\newcommand{\cC}{\mathcal{C}}

\newcommand{\cF}{\mathcal F}
\newcommand{\cG}{\mathcal G}
\newcommand{\cH}{\mathcal{H}}

\newcommand{\cO}{\mathcal{O}}

\newcommand{\CC}{\mathbb{C}}

\newcommand{\crl}{\curvearrowleft}
\newcommand{\crr}{\curvearrowright}

\newcommand{\End}{\operatorname{End}}
\newcommand{\Ext}{\mrm{Ext}}

\newcommand{\KZ}{\mrm{KZ}}
\newcommand{\hd}{\operatorname{hd}}

\newcommand{\fh}{\mathfrak{h}}

\newcommand{\Hh}{\hbar}
\newcommand{\Hom}{\mrm{Hom}}

\newcommand{\Irr}{\mrm{Irr}}

\newcommand{\id}{\mathbbm{1}}
\newcommand{\inv}{^{-1}}

\newcommand{\KK}{\mathbb{K}}

\newcommand{\mbf}{\mathbf}
\newcommand{\mrm}{\mathrm}
\newcommand{\Mod}{\mrm{Mod}}
\newcommand{\mMod}{\mrm{mod}}

\newcommand{\Mor}{\underset{\mrm{Mor}}{\simeq}}
\newcommand{\ind}{{\textup{ind}}}

\newcommand{\NN}{\mathbb{N}}
\newcommand{\QQ}{\mathbb{Q}}

\newcommand{\RR}{\mathbb{R}}

\newcommand{\sgn}{{\textup{sgn}}}

\newcommand{\tif}{\textup{if }}

\newcommand{\triv}{{\textup{triv}}}

\newcommand{\veq}{\rotatebox{90}{$\,=$}}
\newcommand{\usubseteq}{\protect\rotatebox[origin=c]{90}{$\subseteq$}}
\newcommand{\dsubseteq}{\protect\rotatebox[origin=c]{270}{$\subseteq$}}

\newcommand{\ZZ}{\mathbb{Z}}

\usepackage[shortlabels]{enumitem}
\newenvironment{enu}{\begin{enumerate}}{\end{enumerate}}
\newenvironment{enua}{\begin{enumerate}[label=\textup{(\alph*)}]
}{\end{enumerate}}

\newenvironment{eq}{\begin{equation}}{\end{equation}}

\theoremstyle{definition}

\newtheorem{Def}{Definition}[subsection]

\newtheorem{exa}[Def]{Example}
\newtheorem{rmk}[Def]{Remark}

\theoremstyle{plain}
\newtheorem{prop}[Def]{Proposition}
\newtheorem{thm}[Def]{Theorem}
\newtheorem{mainthm}{Theorem}

\newtheorem{theorem}[Def]{Theorem}
\newtheorem{lemma}[Def]{Lemma}
\newtheorem{lem}[Def]{Lemma}
\newtheorem{cor}[Def]{Corollary}

\newtheorem{prob}[Def]{Problem}

\begin{document}

\title{Quantum wreath products and Schur--Weyl duality III}

\begin{abstract}
In the representation theory of finite-dimensional algebras, constructing quasi-hereditary 1-covers is a natural problem. For Hecke algebras of complex reflection groups, one approach required developing the substantial theory of Cherednik algebras. In this paper, the authors develop dual Specht theory, a homological framework designed to detect when Schur functors yield 1-covers for finite-dimensional algebras. Departing from existing Cherednik algebra techniques, our theory combines a natural generalization of dual Specht modules with a spectral sequence argument. In the resulting cover, standard and projective indecomposable modules arise from dual Specht and Young modules, respectively, while simple modules are classified via Schur duality. As our main application, the authors solve the Ginzburg--Guay--Opdam--Rouquier conjecture for Hecke algebras of type D in the separated case, whose representation theory is controlled by a special class of quantum wreath products called the Hu algebras.

Our results build on the first two papers of this series in two essential ways. First, the (dual) Specht modules are constructed by parabolic induction from tensor products of the ``wreath modules'' introduced in the second paper. Second, the explicit formula for a special bar-invariant element of the Hu algebra, obtained in the first paper, controls how the algebra interacts with a natural automorphism. This formula yields the duality between Specht and dual Specht modules, allowing the dual Specht theory to produce 1-covers of the Hu algebras and type D Hecke algebras. Finally, by demonstrating the Schur-concavity of Lusztig's a-function, the authors prove that the c-function order refines the partial order from the wreath product construction.  Under natural conditions, Rouquier's unicity theorem then identifies our Schur functor with the Knizhnik--Zamolodchikov functor from the rational Cherednik algebra.
%last update: \currenttime,\today
\end{abstract}

\author{\sc Chun-Ju Lai}
\address{Institute of Mathematics \\ Academia Sinica \\ Taipei 106319, Taiwan} 
\email{cjlai@gate.sinica.edu.tw}
\thanks{Research of the first author was supported in part by
NSTC grant 114-2628-M-001-001}

\author{\sc Daniel K. Nakano}
\address{Department of Mathematics\\ University of Georgia \\
Athens\\ GA~30602, USA}
\thanks{Research of the second author was supported in part by
NSF grant DMS-2401184}
\email{nakano@math.uga.edu}

\author{\sc Ziqing Xiang}
\address{ Department of Mathematics\\ Southern University of Science and Technology \\Shenzhen 518055, China} 
\thanks{Research of the third author was supported in part by
NSFC (12350710787) and NSFC (12471311)}
\email{xiangzq@sustech.edu.cn}
\keywords{Quantum wreath products, quasi-hereditary covers, Schur duality}
\subjclass{Primary 20C08; Secondary 16G99, 20G43, 20F55}

\maketitle

%\tableofcontents

%=============================================
\section{Introduction}
%=============================================

%=============================================
\subsection{1-Faithful Quasi-Hereditary Covers}
%=============================================
Let $\bH$ be a finite-dimensional algebra over $K$.
A fundamental problem in representation theory is to understand the module category of $\bH$ when it is not semisimple.
One standard approach is to construct another finite-dimensional $K$-algebra $\bS$, which has a quasi-hereditary cover $\cF: \bS$-mod $\to \bH$-mod.
That is, $\bS$-mod is a highest weight category (as defined by Cline, Parshall, and Scott in \cite{CPS88}), and $\cF$ satisfies specific homological properties, including
\begin{equation}\label{eq:HW-cover}
\Hom_\bS(M,N) \cong \Hom_\bH(\cF(M), \cF(N)),
\end{equation}
 for all projective modules $M, N \in \bS$-mod.

However, the condition \eqref{eq:HW-cover} yields limited information regarding higher extensions.
In \cite{Ro08}, Rouquier introduced the notion of an $i$-faithful quasi-hereditary cover (or $i$-cover), which is a quasi-hereditary cover $\cF: \bS$-mod $\to \bH$-mod satisfying additional properties.
Specifically, for any $j \leq i$:
\begin{equation}\label{eq:i-cover}
\Ext^j_\bS(M,N) \cong \Ext^j_\bH(\cF(M), \cF(N)),
\end{equation}
for all modules $M, N \in \bS$-mod which are filtered by the class $\{\Delta(\lambda)\}_\lambda$ of standard modules in the highest weight category $\bS$-mod.
Rouquier also established a unicity theorem asserting that highest weight categories admitting 1-covers are equivalent as categories, under a natural condition on the poset structure.

Constructing 1-covers has become an important problem.
See \cite{RSVV16, Lo16} for technical constructions of 1-covers leading to independent proofs of a conjecture by Varagnolo and Vasserot regarding the connection between cyclotomic rational Cherednik algebras and affine category $\mathcal{O}$.
Recent developments have further explored the properties of such quasi-hereditary covers via Hemmer--Nakano dimension and Ringel duality \cite{C24a,C24b}.

For an arbitrary finite-dimensional algebra $\bH$, the first main contribution of this paper is to simultaneously (1) establish an axiomatic framework which produces a 1-cover of $\bH$, and (2) obtain a classification of simple $\bH$-modules via a Schur duality argument.
Our approach is motivated by the special case \cite{HN04} when $\bH$ is the Hecke algebra $\cH_q(\Sigma_d)$ of type A.

Roughly speaking, we require $\cF$ to be the Schur functor given by $M \mapsto eM$ for some idempotent $e=e^2$ in some algebra $\bS$ such that $e\bS e \cong \bH$.
Consequently, there exists an inverse Schur functor $\cG = \cG_{\textup{Hom}}$ which is left exact,
as well as a spectral sequence (cf. \cite{DEN04})
\begin{equation}
E_2^{i,j} = \Ext^i_\bS( M, R^j\cG(N))
\Rightarrow
\Ext^{i+j}_\bH(\cF(M), N).
\end{equation}
We introduce a new homological notion of {\em dual Specht modules} (see Definition \ref{def:dS}) for $\bH$ which allows one to construct natural indecomposable modules (see Definition \ref{def:YoungMod}) that generalize the Young modules for $\cH_q(\Sigma_d)$.
However, in quantum characteristic two, the classical dual Specht modules of $\mathcal H_q(\Sigma_d)$ do not satisfy Definition~\ref{def:dS}.

An algebra $\bH$ with dual Specht modules is said to admit a {\em dual Specht theory} if it satisfies crucial properties (see Definition \ref{def:dST}) for $\cF$ to be a 1-cover, including that $R^1\cG$ vanishes on all {\em dual Specht modules}.

\begin{mainthm}[Theorems \ref{thm:Young}, \ref{thm:1-cover}, Propositions \ref{prop:indec}, \ref{prop:PDfil}, and Corollary \ref{cor:simple}] \label{thmB}
Let $\bH$ be a finite-dimensional $K$-algebra that admits dual Specht modules $\{S_\lambda\}_{\lambda\in\Lambda}$.
\begin{enua}
\item Each $\lambda$ is associated with an indecomposable module $Y^\lambda$, called the Young module, which has a dual Specht filtration. %both Specht and dual Specht filtrations
\item Additionally, suppose that $\bH$ admits a dual Specht theory. Then,
the sets 
$\{P(\lambda) := \cG(Y^\lambda) ~|~ \lambda \in \Lambda\}$ and 
$\{L(\lambda) := \hd P(\lambda) ~|~ \lambda \in \Lambda\}$ form complete sets of pairwise non-isomorphic projective indecomposable and simple modules in $\bS$-mod, respectively.
Moreover, the functor $\cF$ is a 1-faithful quasi-hereditary cover of $\bH$,
with the set of standard modules given by 
$\{\Delta(\lambda):= \cG(S_\lambda) ~|~ \lambda \in \Lambda\}$.

As a result, 
$\{D^\lambda := \cF(L(\lambda)) ~|~ \lambda \in \Lambda, D^\lambda \neq 0\}$  
forms a complete set of pairwise non-isomorphic  simple $\bH$-modules.
\end{enua}
\end{mainthm}
\vskip .25cm
\begin{rmk}~
\begin{enua}
\item
Even for the special case $\bH = \cH_q(\Sigma_d)$, Theorem \ref{thmB} provides a new proof of (and refines) the results in \cite{HN04}, which relied on the representation theory of algebraic groups and quantum groups (e.g., Ext orthogonality between standard and costandard modules).
\item
In \cite{RSVV16}, a result similar to our Theorem \ref{thmB} is presented, and it describes conditions to lift certain highest weight covers to 1-covers.
Precisely speaking, the covers considered there are the Knizhnik--Zamolodchikov (KZ) functor for cyclotomic Hecke algebras and the functor $\Phi$ arising from Kazhdan--Lusztig fusion products.
Their assumptions, including 0-faithfulness at points of codimension 1, appear closely tied to the nature of these functors and specific knowledge about the standard modules in the corresponding highest weight categories.
See \cite{Lo16} for an independent approach built on checking $(-1)$-faithfulness in codimension 1.
 
In contrast, our assumptions reduce to an explicit calculation of $R^1\cG$ on the dual Specht modules.
We utilize the classical properties of a Schur functor to construct the standard modules from the new ``dual Specht modules''.
We expect our theory will produce new 1-covers for the cyclotomic Hecke algebras when they are non-semisimple, though this is beyond the scope of the present article.
\end{enua}
\end{rmk}

Our general theory suggests that the pursuit of a 1-cover should begin with constructing the dual Specht modules,
which leads to the construction of Young modules.
Utilizing the inverse Schur functor $\cG$, one obtains the modules $\Delta(\lambda) := \cG(S_\lambda)$ and $P(\lambda) := \cG(Y^\lambda)$.  
Furthermore, one must then verify that these form  complete sets of standard and projective indecomposable modules, respectively.

{In the following subsection, we explain how to apply our 1-cover theory to solve the GGOR conjecture for type D in the separated case,
utilizing the wreath module construction for the dual Specht modules over the Hu algebras.}

%=============================================
\subsection{Application: Ginzburg--Guay--Opdam--Rouquier Problem}
%=============================================
Let $\cH_q(W)$ be the Hecke algebra associated with a complex reflection group $W$.
A family of quasi-hereditary covers of $\cH_q(W)$ was constructed by Ginzburg, Guay, Opdam, and Rouquier \cite{GGOR03} using the KZ functor $\KZ_W: \cO(\mathbb{H}(W)) \to \cH_q(W)$-mod from the category $\cO$ of the corresponding rational Cherednik algebra $\mathbb{H}(W)$.
It is natural to compare  these quasi-hereditary covers with those  arising from the quasi-hereditary algebra structure of the $q$-Schur algebra due to Dipper and James \cite{DJ89}.
The following (rephrased) conjecture was proved in type A by Rouquier in \cite{Ro08}, in combination with a 1-cover result established by Hemmer and Nakano in \cite{HN04}.
\begin{prob}[GGOR Conjecture]
For  any complex reflection group $W$,
there exists a finite-dimensional Schur algebra $S_q^W$ containing an idempotent $e\in S_q^W$ such that (1) $eS_q^We \cong \cH_q(W)$, and (2)
the corresponding Schur functor $\mathcal{F}_W: M \mapsto eM$ makes the following diagram commute:
\begin{equation}\label{cd:GGOR}
\begin{tikzcd}
\cO(\mathbb{H}(W)) \arrow{rd}[swap]{\KZ_{W}} \arrow{rr}{\simeq} &{}& S_q^W\textup{-mod} \arrow{ld}{\mathcal{F}_W}
\\
{}&\cH_{q}(W)\textup{-mod}
\end{tikzcd}
\end{equation}
\end{prob}

From the 1990s through the early 2000s, 
the search for generalizations of $q$-Schur algebras beyond type A 
involved various research groups trying distinct constructions, each with a different set of justifications and claimed merits.
For types B and C, these constructions include the work \cite{Gr97} of Green on the hyperoctahedral Schur algebra $S^{\textup{Gr}}$ (see also \cite{BKLW18}),
%for a generalization to arbitrary $d$),
the cyclotomic Schur algebra $S^{\textup{DJM}}$ due to Dipper, James, and Mathas in  \cite{DJM98a,DJM98b} (see also the work \cite{Ar99} of Ariki),
the work \cite{DS00b} of Du and Scott about the $q$-Schur$^2$ algebras $S^2$,
and another Schur algebra $S^{\textup{SS}}$ introduced in \cite{SS99} by Sakamoto and Shoji.

The first solution to the GGOR problem was obtained by Rouquier \cite{Ro08} 
using a Schur algebra $S^{\textup{R}}$ which is Morita equivalent to $S^{\textup{DJM}}$,
if the quantum characteristic $\ell:=\textup{ord}(q)$ (with $\ell=\infty$ if $q$ is not a root of unity) is  greater than the rank $d$ of the complex reflection group.
Rouquier's assumptions on the Cherednik algebras for \eqref{cd:GGOR} to commute were later relaxed by Gordon and Losev \cite{GL14} when $\ell = \infty$.
Both approaches \cite{Ro08, GL14} are valid over the family of complex reflection groups $G(m,1,d) \cong C_m \wr \Sigma_d$.
The authors obtained a solution in \cite{LNX22} for type B (i.e., $m=2$) using a natural realization of $S^{\textup{Gr}}$ in terms of the tensor space for $\ell \geq 4$.

For type D, constructions in the literature include the Schur algebra $S^{\textup{GH}}$ due to Gruber and Hiss in \cite[\S7.2]{GH97},
and three (Morita equivalent) Schur algebras $S^{1.5}, S^{1\D},$ and $S^{2.5}$ introduced by Du and Scott in \cite{DS00a}.
They observed that $S^{1.5}$ is not quasi-hereditary in general ---  counterexamples appear already when $d=4$.
In the work  \cite{FL15} of Fan and Li, two (non-isomorphic) algebras $S^{\textup{FL}}$ and   $S^{\textup{FL},m}$ are studied.
A dimension formula for $S^{\textup{FL}}$ can be obtained\footnote{\cite[(25)]{FL15} was subsequently corrected in private correspondence with the authors. See \cite{DLZ25}.} via a lengthy computation of $\mathrm{SO}(2d)$-orbits.

A type-free construction was developed by Du, Parshall, and Scott in a related program \cite{DPS98, DPS17, DPS24} whose goal was to find the correct enlargement of the tensor space so that its corresponding endomorphism algebra $S^{\textup{DPS}}$ is standardly stratified and admits a highest weight cover after inverting bad primes.

To the best of our knowledge, 1-covers are not available for the aforementioned Schur algebras in type D.
The authors will show that the type B construction for $S^{\textup{Gr}}$ does generalize to type D,
thanks to the following Morita equivalence in \cite{Hu02}:
\begin{equation}\label{eq:MorD}
F_M:\cH_{q}(W^\D_{2m}) \Mor  
\cA(m)
\oplus
\bigoplus\nolimits_{i=0}^{m-1}
\cH_q(\Sigma_i  \times \Sigma_{2m-i}),
\end{equation}
which holds under the invertibility assumption (i.e., the separation condition in \cite{HM09}) that  
\begin{equation}\label{eq:separate}
q^i \neq -1
\quad
\textup{for}
\quad
1\leq i \leq 2m-1,    
\end{equation}
where $\cA(m)$, which we call the {\it Hu algebra}, is a motivating example of the theory of quantum wreath products.
Observe that when $\ell$ is even, \eqref{eq:separate} holds when $\ell \geq 4m$.
In contrast, \eqref{eq:separate} holds automatically for odd $\ell$.
In other words, $\ell$ does not have to be large
for \eqref{eq:MorD} to hold.
The counterpart for $\cH_q(W^\D_{2m+1})$, due to Pallikaros, only involves a direct sum of parabolic subalgebras of type A \cite{Pa94}.

A similar Morita equivalence theorem was proved by Hu and Mathas \cite{HM09} for the Hecke algebra of type $G(d,d,md)$.

%=============================================
\subsection{Quantum Wreath Products}
%=============================================
{In the first part \cite{LNX24} of this series, the authors established a theory of Schur duality for certain finite-dimensional quantum wreath products.
In particular, a Schur functor for the Hu algebra was constructed.

Additionally, \cite{LNX24}  provides an explicit formula \eqref{eq:preciseH1} for the special generator $\Hh\in\cH_q(\Sigma_{2m})$ of the Hu algebra $\cA(m)$. 
This formula turns out to be the key input for the duality between Specht modules for the Hu algebra.
It allows us to compute the effect of the automorphism $\#$ on the central element $Z \in \cA(m)$ explicitly (see Proposition~\ref{prop:Zsharp}).
This identity is what makes the automorphism $\#$ available for the Hu algebras, and it yields the compatibility $f_{\lambda^t}=z_\lambda^2f_\lambda$ of the scalars by which $Z(m)$ acts on Specht modules (see Lemma~\ref{lem:ftrans}).

In the second part \cite{LNX25} of this series, the authors developed a general theory constructing wreath modules over certain (possibly infinite-dimensional) quantum wreath products.
Such a theory enables us to finally construct the desired 1-cover via the (dual) Specht modules.
These Specht modules in fact play the role of cell modules in the sense that they lead to a classification theorem of simple modules over the Hu algebras and the corresponding Schur algebra (see Theorem \ref{thm:PL}).}

{
Let $\Pi_m$ be the set of partitions of $m$.
These dual Specht modules $S_\bblambda$ for $\cH_q(\Sigma_m)\wr \cH(d)$ are
indexed by the set of all $\Pi_m$-indexed  multipartitions
\begin{equation}
    \Omega := \left\{ \bblambda: \Pi_m \to \bigcup\nolimits_{n \geq 0} \Pi_n ~\middle|~ \sum\nolimits_{\lambda \in \Pi_m} \#\bblambda(\lambda) = d \right\} 
\end{equation}
whose partition sizes add up to $d$.

The construction of $S_\bblambda$ is built on a technical isomorphism (Theorem \ref{thm:dualS}), which lifts the isomorphism of $\cH_q(\Sigma_m)$-modules
\[
(S^\lambda)^\# \cong (S^{\lambda^t})^*
\quad
\textup{where}
\quad
\lambda \vdash m
\]
to the quantum wreath product level.
{For the Hu algebras, the hypotheses of Theorem \ref{thm:dualS} are verified using the explicit formula from \cite{LNX24} discussed above.}

We further show that the (dual) Specht modules satisfy the desired homological properties, 
consequently proving that the seemingly artificial partial order (Definition \ref{def:KPO}) on $\Omega$ arises naturally from the K\"unneth formula (see Lemma \ref{lem:gHuExt1} and \cite[Lemma 5.5.1]{LNX25}).}
Our detailed analysis leads to the following  results:

\begin{mainthm}[Propositions \ref{prop:WSF}, \ref{prop:gHuDS}, Theorems \ref{thm:R1G}, \ref{thm:PL}, \ref{thm:SDHWC}] \label{thmC}
Let $\cH_q(\Sigma_m)\wr \cH(d)$ be the generalized Hu algebra ($d\geq 2$) over a field $K$ such that $d!\in K^\times$ and the separation condition \eqref{eq:separate} holds.  
Assume that there is an embedding $\iota: \cH_q(\Sigma_m)\wr \cH(d) \hookrightarrow \cH_q(\Sigma_{md})$ such that $\iota|_{\cH_q(\Sigma_m)^{\otimes d}}$ is the identity map. 
Let $\bS(n,m \wr d)$ be its corresponding wreath Schur algebra. 
\begin{enua}
\item The Schur functor $\cF$ and the inverse Schur functor $\cG$ are constructed, provided $n \geq md$.
\item If $\ell \geq 4$, then $\cH_q(\Sigma_m)\wr \cH(d)$ admits dual Specht modules $S_\bblambda$.
In particular, the Young modules $Y^\bblambda$ are constructed.
\item If $n\geq md$ and $\ell \geq 4$, then $R^1\cG(S_\bblambda) = 0$ for all $\bblambda$.
Moreover, $\cH_q(\Sigma_m)\wr \cH(d)$ admits a dual Specht theory. 
Consequently, the Schur functor $\cF$ is a 1-cover.
Additionally, the sets of pairwise non-isomorphic projective indecomposable and irreducible $\bS(n, m\wr d)$-modules are given by $\{P(\bblambda) := \cG(Y^\bblambda)\}_\bblambda$ and $\{L(\bblambda) := \hd P(\bblambda)\}_\bblambda$, respectively.
\end{enua}
\end{mainthm}
Since these assumptions all hold for the case $d=2$, as a corollary, 
the Hecke algebra $\cH_q(W^\D_{2m})$ has a 1-cover (see Theorem \ref{thm:SDHWC}) given by 
the Schur functor $\mathcal{F}_\D$ for the algebra 
\[
S^\D_q(n,2m):= \bS(n,m\wr 2) \oplus \Big(
\bigoplus_{i < m} S_q(n,i) \otimes S_q(n, 2m-i)
\Big),
\]
provided $n \geq 2m$ and
\begin{equation}\label{eq:finalell}
\ell = \infty
\quad
\textup{or}
\quad
\ell \geq\begin{cases}
 5 &\tif q \textup{ is an odd root of unity};
\\
 4m &\tif q \textup{ is an even root of unity}.
\end{cases}
\end{equation}
Finally, we apply our general framework to resolve the GGOR conjecture for type D in the separated case. 
In particular, a detailed analysis of the highest weight cover of the rational Cherednik algebra provided by the KZ functor is required,
where the standard objects are ordered via the {$\bc$-function of \cite{GGOR03}, which is expressed in terms of Lusztig's $\ba$- and $\bA$-functions.}

On the other hand, our wreath product construction allows us to establish a modified dominance order for our dual Specht modules $S_\bblambda$.
By demonstrating the Schur-concavity of the $\ba$-function, 
we prove that the $\bc$-function order $\leq_1$ refines our modified dominance order $\leq_2$ (see Proposition \ref{prop:refinePO}). 
This critical lemma enables us to apply Rouquier's unicity theorem, yielding the following result:
%----------------------
\begin{mainthm}[Theorem \ref{thm:GGORD}]
Suppose that $K = \bbC$, $q=\exp(2\pi\sqrt{-1}c)$ with $c\in\QQ_{>0}$ or $c\notin\QQ$. If $n\geq 2m$ and \eqref{eq:finalell} holds,
then the following diagram commutes:
\begin{equation}\label{cd:GGORD1}
\begin{tikzcd}
\cO(\mathbb{H}(W^\D_{2m})) \arrow{rd}[swap]{\KZ_{\D}} \arrow{rr}{\simeq} &{}& S_q^\D(n,2m)\textup{-mod} \arrow{ld}{\mathcal{F}_\D}
\\
{}&\cH_{q}(W^\D_{2m})\textup{-mod}
\end{tikzcd}
\end{equation}
Moreover, $(\cO(\bbH), \KZ_\D)$ and $(S_q^\D(n,2m)\textup{-mod}, \cF_\D)$ are equivalent highest weight covers.
\end{mainthm}
Therefore, the GGOR problem for type $\D_{2m} = G(2,2,2m)$ is solved for those choices of parameters for which both 1-covers exist.
Note that Theorem \ref{thmC} is not restricted to the case $d=2$.
We anticipate that the construction for $d>2$ will be used in resolving the GGOR conjecture for $G(d,d,dm)$.

\vskip .25cm
\noindent {\bf Acknowledgements.}
{
The authors thank Ben Webster for pointing out an oversight in an earlier version of \cite{LNX25}. 
Section \ref{sec:GGOR} is a byproduct of our discussions with him.
Moreover, we acknowledge Meinolf Geck and Jun Hu for helpful comments.
}
%=============================================
\section{Dual Specht Theory}\label{sec:dST}
%=============================================
Throughout this paper, let $K$ be a commutative ring. 
All algebras are assumed to be free of finite rank as $K$-modules, %All algebras are free algebras over $K$, 
and all tensor products are considered to be over $K$, unless specified otherwise.
For this section, assume that $K$ is a field, and let $\bS, \bH$ be arbitrary finite-dimensional $K$-algebras.

{
Suppose that $e\bS e\cong \bH$  for some idempotent $e\in \bS$.
The goal of this section is to develop a {\em dual Specht theory} so that the Schur functor $\cF:\bS$-mod $\to \bH$-mod, $\cF(M) := eM$ is a quasi-hereditary 1-cover. See Theorem \ref{thm:1-cover}.
}
%==========================================
\subsection{Quasi-Hereditary 1-Covers}  \label{sec:QHA}
%==========================================
The formal definition of a highest weight category is due to Cline, Parshall, and Scott \cite{CPS88}. For the category of finite-dimensional modules over a finite-dimensional $K$-algebra $\bS$,  
Rouquier  \cite{Ro08} has the following equivalent formulation. 
A category $\mathcal{C} = \bS$-mod is a {\em highest weight category} if there is a collection $\{\Delta(\lambda)\}_{\lambda \in \Lambda}$ of standard modules indexed by a finite poset $\Lambda$ such that the following hold:
\enu[\textup{(}H1\textup{)}]
\item $\End_{\mathcal{C}}(\Delta(\lambda)) = K$ for all $\lambda \in \Lambda$;
\item If $\Hom_{\mathcal{C}}(\Delta(\lambda) , \Delta(\mu)) \neq 0$ then $\lambda \leq \mu$;
\item If $\Hom_{\mathcal{C}}(\Delta(\lambda), M) = 0$ for all $\lambda \in \Lambda$ then $M = 0$;
\item For each $\lambda \in \Lambda$ there is a projective module $P(\lambda) \in \mathcal{C}$ and a surjection $P(\lambda) \twoheadrightarrow \Delta(\lambda)$ whose kernel $\mathrm{ker}(P(\lambda) \to \Delta(\lambda))$ has a $\Delta$-filtration with subquotients of the form $\Delta(\mu)$ with $\mu > \lambda$.
\endenu
Recall that the finite-dimensional algebra $\bS$ is called {\em quasi-hereditary} if and only if $\bS$-mod is a highest weight category.

For now, suppose that $\bS$-mod is a highest weight category.
A {\em 1-cover} (i.e., a 1-faithful quasi-hereditary cover following \cite{Ro08}) of $\bH$ is an exact functor $\cF: \bS$-mod $\to \bH$-mod  satisfying the following:
\enu[\textup{(}C1\textup{)}]
\item $\bH = \End_\bS(P)^{\textup{op}}$ for some projective $\bS$-module $P$ and $\cF\cong \Hom_\bS(P,-)$;
\item for projective $\bS$-modules $M,N$,
$
\Hom_\bS(M,N) \cong  \Hom_\bH(\cF(M),\cF(N));
$
\item for modules $M,N \in \bS$-mod admitting $\Delta$-filtrations and $0\leq i \leq 1$,
\[
\Ext^i_\bS(M,N) \cong  \Ext^i_\bH(\cF(M),\cF(N)).
\]
\endenu
{
The existence of a 1-cover is not guaranteed in general.
When one exists, there is a uniqueness result (see Proposition \ref{Ro449}).
}

%=============================================
\subsection{Dual Specht modules and Young modules}
%=============================================
{To develop dual Specht theory, we begin by introducing a new notion of dual Specht modules for an arbitrary finite-dimensional algebra.
Note that in quantum characteristic two, the Hecke algebra $\cH_q(\Sigma_d)$ does not admit dual Specht modules under this new definition. 
}
%--------------------------
\begin{Def}\label{def:dS}
%--------------------------
 
The algebra $\bH$  admits {\em dual Specht modules} if
there is a family of non-isomorphic finite-dimensional $\bH$-modules (called dual Specht modules) $\{S_\lambda\}_{\lambda\in \Lambda}$ indexed by a finite poset $(\Lambda, \leq)$ such that 
\[
\Hom_{\bH}(S_\lambda , S_\mu) \neq 0 \Rightarrow \lambda \leq \mu,
\quad
\Ext^1_\bH(S_\lambda, S_\mu) \neq 0 \Rightarrow \lambda < \mu,
\quad
\End_{\bH}(S_\lambda) = K \textup{ for all }\lambda \in \Lambda.
\]
\end{Def}
%--------------------------
{
The construction of Young modules for $\cH_q(\Sigma_m)$ appears in \cite[Sections 4.3--4.4]{HN04} by building a Specht filtration from bottom to top.
An alternative construction using a dual Specht filtration from top to bottom was sketched in \cite[Section 4.7]{HN04}.
We generalize this alternative approach to construct a finite-dimensional indecomposable $\bH$-module $Y^{\lambda}$ for each $\lambda\in \Lambda$. 
These $\bH$-modules are called the Young modules.}
%--------------------------
\begin{Def}[Construction of Young modules]\label{def:YoungMod}
%--------------------------
Suppose that $\bH$ admits dual Specht modules.
Set $\lambda_0:= \lambda$ and $Y^0 := S_{\lambda_0}$. 
If $\Ext^1_{\bH}(Y^0, S_\nu) = 0$ for all $\nu \in \Lambda$, set $Y^\lambda := S_\lambda$.
Otherwise, we construct $Y^i$ inductively for $i \geq 1$ as follows:
let $\lambda_i$ be a minimal element in the poset $(\Lambda, \leq)$ such that $\Ext^1_{\bH}(Y^{i-1}, S_{\lambda_i}) \neq 0$.
It follows from Definition \ref{def:dS} that 
%$\lambda = \lambda_0 < \lambda_1 < \dots < \lambda_i$. 
\[
\lambda < \lambda_j ~(\textup{for all } 1\leq j \leq i),
\quad
\lambda_j\textup{'s are all distinct, and}
\quad
\lambda_j \not< \lambda_k \textup{ for } j>k.
\]
In other words, there is a non-split exact sequence
\begin{equation}\label{eq:YSES}
0 \to S_{\lambda_i}^{\oplus t_i} \to Y^i \rightarrow Y^{i-1} \to 0,
\end{equation}
for some $\bH$-module $Y^i$, where $t_i:=\dim \Ext^1_{\bH}(Y^{i-1}, S_{\lambda_i})$. 
Since $\Lambda$ is finite and the choice of $\lambda_i$'s is irredundant, the process must terminate, say, at the $N$th step.
Let $Y^\lambda := Y^N$, which is finite-dimensional.
\end{Def}
%--------------------------
Note that it appears that $Y^\lambda$ depends on the choice of minimal elements $(\lambda_1, \dots, \lambda_N)$. 
However, once a dual Specht theory is established, $Y^\lambda$ does not depend on this choice.
Indeed, in that case $Y^\lambda\cong\cF(P(\lambda))$, where $P(\lambda)$ is the projective cover of $L(\lambda)=\hd\Delta(\lambda)$ (see Proposition~\ref{prop:PDfil}).
%--------------------------------------------------------------------
\begin{theorem} \label{thm:Young}
Suppose that $\bH$ admits dual Specht modules. Then, for each $\lambda \in \Lambda$:
\begin{enua}
    \item The Young module $Y^\lambda$ is indecomposable.
    \item $Y^\lambda$ has a dual Specht filtration 
$\{0\} = F_{0} \subsetneq F_{1} \subsetneq \dots \subsetneq F_{t} = Y^{\lambda}$ 
such that 
\[
F_{i}/F_{i-1}\cong S_{\mu_i}
\quad
\textup{where}
\quad
\mu_t = \lambda < \mu_i 
\quad\textup{for}
\quad
i < t. 
\]
\end{enua}
\end{theorem} 
%--------------------------------------------------------------------
\begin{proof} 
Part (a) can be proved using the argument in \cite[Proposition 4.4.2]{HN04}. 
That is, one can show that each $Y^i$ is indecomposable by induction on $i$.
For the initial case, $Y^0 := S_\lambda$ is indecomposable since $\End_\bH(S_\lambda) \cong K$.
The inductive case is proved by showing that the only nonzero idempotent in $\End_\bH(Y^i)$ is the identity.
Part (b) follows directly from the construction.
\end{proof} 
%=============================================
\subsection{Schur Functors}\label{sec:SF}
%=============================================
{A key ingredient of our dual Specht theory concerns the behavior of dual Specht modules under the right derived functor $R^1\cG$ of the inverse Schur functor $\cG$.
This functor exists in the following framework (cf.~\cite{DEN04}).}

From now on, assume that there is an idempotent $e = e^2 \in \bS$ such that $e \bS e \cong \bH$.
Thus, one can define the Schur functor $\cF$ and the (right adjoint) inverse Schur functor $\cG = \cG_{\textup{Hom}}$, respectively, by
\begin{equation}
\begin{array}{rlcrl}
\cF: \bS\textup{-mod}&\to \bH\textup{-mod},&
&\cG: \bH\textup{-mod}&\to \bS\textup{-mod}.
\\
M&\mapsto  eM&
&M&\mapsto \Hom_{\bH}(e\bS , M)
\end{array}
\end{equation}
Here, 
$\cF$ is covariant and exact,
$\cG$ is left exact, and
$\cF\circ\cG \cong \id$.

In general, it is not true that  $\cG \circ \cF \cong \id$. 
In the special case when $\bS = S_q(n,d)$, $\bH = \cH_q(\Sigma_d)$, and $\ell  \geq 4$, one has $\cG \circ \cF(M) \cong M$  for those $M$ having a Weyl filtration, provided $n \geq d$.
%-----------------------------------------------------------------------------------------------
%\begin{prop}[]

Suppose that $M \in \bS$-mod, and $N\in \bH$-mod. %, and $\ell \geq 3$. 
%If there is an idempotent $e\in \bS$ such that $e \bS e \cong \bH$,
It follows from \cite[Section 2.2]{DEN04} that there is a first quadrant Grothendieck spectral sequence
\[
E_2^{i,j} = \Ext^i_\bS\big( M, R^j\cG(N) \big) \Rightarrow 
\Ext^{i+j}_\bH\big(\cF(M), N\big).
\]
%\end{prop}
%-----------------------------------------------------------------------------------------------
As a consequence, $E_2^{0,0} \cong E_0$.
Moreover,  there is a five-term exact sequence
$0 \to E_2^{1,0} \to E_1 \to E_2^{0,1} \to E_2^{2,0} \to E_2$.
The following isomorphisms can be verified using the five-term exact sequence:
%-----------------------------------------------------------------------------------------------
\begin{thm}\label{thm:exactseq}
Suppose that $M, N \in \bS$-mod and $N \cong \cG \circ \cF(N)$. Then,
\begin{enua}
\item 
$\Hom_\bS\big( M, N \big) \cong\Hom_\bH\big(\cF(M), \cF(N)\big)$.
\item
There is an exact sequence
\[
0 
\to \Ext^1_\bS\big( M, N \big) 
\to \Ext^{1}_\bH\big(\cF(M), \cF(N)\big)
\to \Hom_\bS\big( M, R^1\cG(\cF(N)) \big).
\]
In particular, $\Ext^1_\bS\big( M, N \big) \cong\Ext^1_\bH\big(\cF(M), \cF(N)\big)$ if $R^1\cG(\cF(N)) = 0$.
\end{enua}
\end{thm}
In other words, Theorem \ref{thm:exactseq}(b) suggests that a vanishing condition for $R^1\cG$ is needed in order to identify the $\Ext^1$ groups.
One can show that the inverse Schur functor takes a Young module $Y^\lambda$ to an indecomposable $P(\lambda) := \cG(Y^\lambda)$ as follows. 
The axioms for $\{P(\lambda)\}_{\lambda}$ to form a complete set of non-isomorphic projective indecomposables will be discussed in Section \ref{sec:dStheory}.

%--------------------------
\begin{prop}\label{prop:indec}
Suppose that $\bH$ admits dual Specht modules $\{S_\lambda\}_{\lambda \in \Lambda}$.
  Then $P(\lambda) := \cG(Y^\lambda)$ is indecomposable for all $\lambda \in \Lambda$.
\end{prop}
%--------------------------
\begin{proof}
Suppose that  $\cG(Y^{\lambda})\cong M_{1}\oplus M_{2}$ with $M_{j}\neq 0$. 
Applying the Schur functor, it follows that 
\begin{equation}
Y^{\lambda}\cong e \cG(Y^{\lambda}) \cong e{M_{1}}\oplus e{M_{2}}.
\end{equation}
Since the socle of $\cG(Y^{\lambda})$ consists of simple modules not killed by $e$:
\[
\Hom_\bS(V, \cG(N)) \cong \Hom_\bH(eV, N) = 0
\quad
\textup{whenever}
\quad
eV = 0,
\]
one has $eM_{j}\neq 0$ for $j=1,2$. 
This contradicts the indecomposability of $Y^{\lambda}$. 
\end{proof}
%=============================================
\subsection{Dual Specht theory}\label{sec:dStheory}
%=============================================
Let $\bH$ and $\bS$ be as in Section \ref{sec:SF}.
In this section, we establish a general theory for determining whether the Schur functor $\cF$ is a 1-cover.
%--------------------------
\begin{Def} \label{def:dST}
%--------------------------
The algebra $\bH$ {\em admits a dual Specht theory} if:
\begin{enumerate}[\textup{(}S1\textup{)}]
\item $\bH$ admits dual Specht modules $\{S_\lambda\}_{\lambda \in \Lambda}$ for some poset $(\Lambda, \leq)$.
\item $P(\lambda) := \cG(Y^\lambda)$ is projective for all $\lambda \in \Lambda$.
\item $R^1\cG(S_\lambda) = 0$ for all $\lambda \in \Lambda$.
\item Let $\lambda, \mu \in \Lambda$. Then, $P(\lambda)\cong P(\mu)$ if and only if $\lambda=\mu$. 
\item 
If $P \in \bS$-mod is projective indecomposable, then $P\cong P(\lambda)$ for some $\lambda \in \Lambda$. 
%\red{The set $\Lambda$ is in bijection with the set of isomorphism classes of simple modules of $\bH_\KK := \KK \otimes_\cA \bH$, where $\KK$ is a splitting field for $\bH$.}
\end{enumerate} 
\end{Def}
%--------------------------
%and $\cF(\nabla(\lambda)) = S^\lambda$.
%An $\bS$-module is said to have a {\em $\Delta$-filtration} (or $\nabla$-filtration, resp.) if it is filtered by modules of the form $\Delta(\lambda)$ (or $\nabla(\lambda)$, resp.) for some $\lambda \in \Lambda$.
Let $\Delta(\lambda) := \cG(S_{\lambda})$.
It is immediate that $\cF(\Delta(\lambda)) = \cF \circ \cG (S_{\lambda}) = S_{\lambda}$.
{Having dual Specht theory available, 
one can give a criterion for a certain $\bH$-module to have a $\Delta$-filtration as in the lemma below.}
%--------------------------------------------------------------------
\begin{lem}\label{lem:WeylFil}
Suppose that $\bH$ admits a dual Specht theory, and $N \in \bH$-mod.
Then,
\begin{enua}
\item The $\bS$-module $\cG(N)$ has a $\Delta$-filtration
if and only if
$N$ has a dual Specht filtration.  
\item $\cG \circ \cF(M) \cong M$ for any  $\bS$-module $M$ admitting a $\Delta$-filtration.
\end{enua}
\end{lem}
%--------------------------------------------------------------------
\begin{proof}
Since $\cF$ is exact, a $\Delta$-filtration of $\cG(N)$ induces a filtration of $N$ by modules of the form $\cF(\Delta(\lambda)) = S_\lambda$, i.e., a dual Specht filtration.
While $\cG$ is only left exact, it still takes a short exact sequence of the form $0 \to S_\lambda \to N \to Q \to 0$ to another short exact sequence $0 \to \Delta(\lambda) \to \cG(N) \to \cG(Q) \to 0$, since $R^1\cG(S_\lambda) = 0$.
Hence, part (a) holds by induction on the length of the filtration.

For (b), let $\tau(M)$ be the largest submodule of $M$ whose composition factors are killed by $\cF$. 
We first show that $\tau(\Delta(\nu)) = 0$ for every $\nu \in \Lambda$. 
Recall that $\Delta(\nu) = \Hom_\bH(e \bS, S_\nu)$. 
Suppose that $U \subseteq \Delta(\nu)$ is a submodule with $eU = 0$. 
Then for any $f \in U \subseteq \Hom_\bH(e \bS, S_\nu)$,
\[
f(es) = f(e \cdot es) = (e\cdot (s \cdot f))(e) = 0,
\quad
\textup{for all}
\quad
s \in  \bS,
\]
which implies $f = 0$ since $s\cdot f \in U$, and hence $U = 0$. 
That is, $\Delta(\nu)$ contains no nonzero submodule whose composition factors are killed by $\cF$, and hence $\tau(\Delta(\nu)) = 0$.
Since $M$ has a $\Delta$-filtration, $\tau(M)$ must be zero, and hence we can apply \cite[Theorem 3.3]{DEN04} to show that the unit map is an embedding $M \hookrightarrow \cG \circ \cF(M)$ (note that $\tau(M)$ is denoted by $S(M)$ therein).

Finally, since $M$ has a $\Delta$-filtration, $\cF(M)$ has a dual Specht filtration.
Moreover,  $\cG\circ\cF(M)$ has a $\Delta$-filtration by (a). 
This means that both $M$ and $\cG\circ\cF(M)$ have $\Delta$-filtrations, with the same multiplicities.  
Therefore, $M \cong \cG\circ\cF(M)$.
\end{proof}
{Additionally, the dual Specht theory allows us to determine all projective indecomposable $\bS$-modules as follows.}

\begin{prop}\label{prop:PDfil}
Suppose that $\bH$ admits dual Specht modules. 
\begin{enua}
\item
If (S3) holds, then
$P(\lambda)$ has a $\Delta$-filtration for all $\lambda \in \Lambda$.
Moreover, the top layer of the filtration is $\Delta(\lambda)$, and all other subquotients are of the form $\Delta(\mu)$ with $\mu > \lambda$.
\item 
If (S2) and (S3) hold, then
the head of $P(\lambda)$, denoted by $L(\lambda) := \hd P(\lambda)$, coincides with $\hd \Delta(\lambda)$.
\item 
If $\bH$ admits a dual Specht theory, then
$\{P(\lambda)~|~ \lambda\in \Lambda\}$ forms a complete set of non-isomorphic projective indecomposable $\bS$-modules. 
Consequently, $\{L(\lambda) ~|~ \lambda\in \Lambda\}$ forms a complete set of non-isomorphic simple $\bS$-modules. 
\end{enua}
\end{prop}
%--------------------------------------------------------------------
\begin{proof}
For (a), by applying $\cG$ to the short exact sequence \eqref{eq:YSES} for each $i$, one obtains a short exact sequence $0 \to \Delta(\lambda_i)^{\oplus t_i} \to \cG(Y^i) \to \cG(Y^{i-1}) \to 0$, since $R^1\cG(S_{\lambda_i}) = 0$.
Note that $\lambda_i > \lambda$ thanks to the construction of the Young module.

For (b), one has a %n injection $S^{\lambda}\hookrightarrow Y^{\lambda}$, which induces a surjection 
%$(Y^{\lambda})^{*} \twoheadrightarrow S_\lambda $.
%Since $Y^{\lambda}\cong (Y^{\lambda})^{*}$, we have a 
surjection 
$Y^{\lambda}\twoheadrightarrow  S_\lambda$.
Now,  $R^{1}\cG(S_\mu) =0$ for all $\mu\in \Lambda$. 
This shows $P(\lambda) \twoheadrightarrow \Delta(\lambda)$. 
Therefore, the head of $\Delta(\lambda)$ is $L(\lambda) = \hd P(\lambda)$ since $P(\lambda)$ is projective indecomposable. 

{For (c), because $P(\lambda)$ is indecomposable by Proposition \ref{prop:indec}, and is projective by (S2), it follows from (S4)--(S5) that the set $\{P(\lambda) \mid \lambda \in \Lambda\}$ contains all such modules without redundancy. 
}
\end{proof}

\begin{cor}\label{cor:simple}
Suppose that $\bH$ admits a dual Specht theory. Then the complete set of pairwise non-isomorphic  simple $\bH$-modules is given by
    \[
    \{D^\lambda := \cF(L(\lambda)) ~|~ \lambda \in \Lambda, D^\lambda \neq 0\}.
    \]
\end{cor}
%--------------------------------------------------------------------
%--------------------------------------------------------------------
\subsection{Dual Specht theory guarantees the existence of 1-covers}

%--------------------------------------------------------------------
Finally, we show that the Schur functor is a 1-cover for any finite-dimensional algebra admitting a dual Specht theory. 

%--------------------------------------------------------------------
\begin{thm} \label{thm:1-cover}
Suppose that $\bH$ admits a dual Specht theory.
Then the Schur functor $\cF: \bS$-mod $\to \bH$-mod, $M \mapsto eM$ is a 1-cover  with standard modules $\{\Delta(\lambda) := \cG(S_\lambda)\}_{\lambda \in \Lambda}$. 
\end{thm} 
%--------------------------------------------------------------------
\begin{proof}
{First, we show that $\bS$-mod is a highest weight category.
By Theorem \ref{thm:exactseq}(a), 
\[
\Hom_\bS\big( \Delta(\lambda), \Delta(\mu) \big) 
\cong\Hom_\bH\big(\cF(\Delta(\lambda)), \cF(\Delta(\mu))\big)
\cong\Hom_\bH\big(S_\lambda, S_\mu\big)
\]
for all $\lambda, \mu \in \Lambda$.
Hence, (H1) and (H2) follow from (S1).

For (H3), recall from Proposition \ref{prop:PDfil}(b) that $L(\lambda) = \hd\Delta(\lambda)$
for any $\lambda \in \Lambda$. 
Thus, there is a short exact sequence
\begin{equation}
0 \to \textup{Ker}(\Delta(\lambda) \to L(\lambda)) \to \Delta(\lambda) \to L(\lambda) \to 0.
\end{equation}
Applying the functor $\Hom_\bS(-,M)$, one obtains that $\Hom_\bS(L(\lambda),M) =0$ since $\Hom_\bS(\Delta(\lambda),M) = 0$. 
Thus, $M$ must be zero, since $\{L(\lambda)\}_{\lambda \in \Lambda}$ is a complete set of simples due to Proposition \ref{prop:PDfil}(c).
Finally, (H4) follows from Proposition \ref{prop:PDfil}(a).}

Next, we show that $\cF$ is a 1-cover.
(C1) follows from the fact that $P = \bS e$ is a projective $\bS$-module.
For (C2), 
it follows from Proposition \ref{prop:PDfil} that the projective module $N$ has a $\Delta$-filtration.
It then follows from Lemma \ref{lem:WeylFil}(b) that $N \cong  \cG \circ \cF(N)$.
Thus, (C2) follows from Theorem \ref{thm:exactseq}(a).
For (C3), 
since $\cF$ is exact, $\cF(N)$ has a dual Specht filtration.
Hence, $R^1\cG (\cF(N)) = 0$ thanks to (S3).
By Theorem~\ref{thm:exactseq}(b) and Lemma \ref{lem:WeylFil}(b), (C3) holds.
Therefore, $\cF$ is a 1-cover of $\bH$.
\end{proof}
%--------------------------------------------------------------------

%=============================================
\section{Quantum Wreath Products}\label{sec:QWP}
%=============================================
{
In this section, basic properties and representation theory of quantum wreath products are reviewed.
}
%=============================================
\subsection{Definition} \label{sec:QWPprop}
%=============================================
Let $V$ be a $K$-vector space, and let $d\geq 2$ be an integer.
For $Z \in V\otimes V$ and $\phi \in \End_K(V\otimes V)$,  define
  $Z_i \in V^{\otimes d}$ and $\phi_i \in \End_K(V^{\otimes d})$ by:
\begin{equation}\label{def:Zphii}
Z_i :=  1^{\otimes i-1} \otimes Z \otimes 1^{\otimes d-i-1},
\quad
\phi_i(v_1 \otimes \dots \otimes v_d) := v_1 \otimes \dots \otimes v_{i-1}\otimes  \phi(v_{i} \otimes v_{i+1}) \otimes v_{i+2} \otimes \dots \otimes v_{d}.
\end{equation} 
Let $B$ be an associative $K$-algebra.
Let $Q = (R,S,\rho,\sigma)$ where $R, S \in B\otimes B$ and $\rho, \sigma \in  \End_K(B\otimes B)$ such that $\sigma$ is an automorphism. 

\begin{Def}
The {\em quantum wreath product} is the associative $K$-algebra generated by the algebra $B^{\otimes d}$ and $H_1, \dots, H_{d-1}$ such that the following conditions hold, for $1\leq k \leq d-2, 1\leq i \leq d-1, |j-i|\geq 2$:
\begin{align}
&\textup{(braid relations) } \label{def:BR}
& H_k H_{k+ 1} H_k = H_{k+ 1} H_k H_{k + 1},  \quad H_i H_j = H_j H_i,
\\
&\textup{(quadratic relations) }\label{def:QR}
&H_i^2 = {S_i}H_i + {R_i},
\\
&\label{def:WR}\textup{(wreath relations) }
& H_ib = \sigma_i(b)H_i + \rho_i(b) \quad ( b\in B^{\otimes d}).
\end{align}
\noindent Often,  this algebra is denoted by $B \wr_Q \cH(d)$, or $B \wr \cH(d)$ whenever it is convenient.
\end{Def}
It follows from \eqref{def:BR} that $H_w := H_{i_1} \dots H_{i_N} \in B\wr \cH(d)$ for any $w\in \Sigma_d$ is well-defined, where $s_{i_1} \dots s_{i_N}$ is a reduced expression of $w$. 

In this paper, it will be assumed that
$B\wr \cH(d)$ has a PBW basis, i.e., a $K$-basis of $B \wr \cH(d)$  of one of the following forms:
\begin{equation}\label{def:PBWbases}
\{ (b_{i_1} \otimes \dots\otimes b_{i_d}) H_w ~|~ i_j \in I, w\in \Sigma_d\},
\quad\textup{or}\quad
\{H_w (b_{i_1} \otimes \dots\otimes b_{i_d})  ~|~ i_j \in I, w\in \Sigma_d\},
\end{equation}
where $\{b_i ~|~ i\in I\}$ is a $K$-basis of $B$.
This assumption is equivalent to the conditions \cite[(P1)--(P9)]{LNX25} on $Q$.

\begin{Def}\label{def:Hpar}
For a composition $\mu = (\mu_1, \dots, \mu_r) \vDash d$, denote the associated Young subgroup by
$\Sigma_\mu := \Sigma_{\mu_1} \times \dots \times \Sigma_{\mu_r}$.
The {\em parabolic subalgebra} of $B \wr \cH(d)$ corresponding to $\mu \vDash d$ is given by
\[
B \wr \cH(\mu) := \< B^{\otimes d}, H_w~|~ w \in \Sigma_\mu\> \subseteq B \wr \cH(d).
\]
The tensor product $(B \wr \cH(\mu_1)) \otimes \dots \otimes (B\wr \cH(\mu_r))$ among quantum wreath products with the same parameter $Q$ is canonically isomorphic to the parabolic subalgebra $B\wr \cH(\mu)$ via $(b_1 H_{w_1})\otimes \dots \otimes (b_r H_{w_r}) \mapsto (b_1\otimes \dots \otimes b_r) H_{w_1 \dots w_r}$.
\end{Def}
%=================================
\subsection{Wreath Modules}
%=================================
{
Suppose that $M$ is a $B$-module, and that the elements $S, R \in B\otimes B$ act on $M\otimes M$ by the following maps $S^M$ and $R^M$, respectively:
\begin{equation}\label{eq:SMRM}
S^M(m) = \chi(S)m,
\quad
R^M(m) = \chi(R)m,
\end{equation}
for some map $\chi: B\otimes B \to K$.
Let $\tcH = \tcH^\chi_d$ be the twisted Hecke algebra generated by $h_1, \dots, h_{d-1}$ subject to the following relations:
\begin{align}
&\textup{(braid relations) } 
& h_k h_{k+ 1} h_k = h_{k+ 1} h_k h_{k + 1},  \quad h_i h_j = h_j h_i,
\\
&\textup{(quadratic relations) }
&h_i^2 = {\chi(S)}h_i + {\chi(R)},
\end{align}
where $1\leq k \leq d-2, 1\leq i \leq d-1, |j-i|\geq 2$. 
Suppose that $N$ is an $\tcH$-module.
Let $\sigma^M, \rho^M \in \End_K(M\otimes M)$, and let 
\begin{equation}\label{def:QWPAction}
    H_i \cdot (m\otimes n) := \sigma^M_i(m) \otimes (h_i\cdot n) + \rho^M_i(m) \otimes n.
\end{equation}
In \cite{LNX25}, the following conditions are introduced in order to characterize precisely when the $H_i$-action \eqref{def:QWPAction} induces a  $B\wr\cH(d)$-action on the vector space $M \wr N := M^{\otimes d}\otimes N$:
\begin{align}
&\label{def:M2}\tag{M2}
\sigma^M_i(b\cdot m) = \sigma_i(b)\cdot \sigma^M_i(m),
\quad
\rho_i^M(b\cdot m) =  \sigma_i(b) \cdot \rho^M_i(m) + \rho_i(b) \cdot m,
\\
&\label{def:M4}\tag{M4}
\chi(S)(\sigma^M_i)^2+\rho^M_i\sigma^M_i+\sigma^M_i\rho^M_i = l^M_{S_i}\sigma^M_i,
\quad
\chi(R) (\sigma^M_i)^2+(\rho^M_i)^2 = l^M_{S_i} \rho^M_i + l^M_{R_i},
\end{align}
where $l^M_X\in \End_K(M^{\otimes d})$ means the left $B^{\otimes d}$-action by $X \in {B^{\otimes d}}$.
\begin{align}    
&\label{def:M5}\tag{M5}
\sigma^M_i\sigma^M_j\sigma^M_i = \sigma^M_j \sigma^M_i \sigma^M_j,
\quad
\rho^M_i\sigma^M_j\sigma^M_i = \sigma^M_j \sigma^M_i \rho^M_j,
\\
&\label{def:M6}\tag{M6}
\rho^M_i\sigma^M_j\rho^M_i = \chi(S) \sigma^M_j \rho^M_i \sigma_j^M + \rho_j^M\rho_i^M\sigma_j^M + \sigma_j^M \rho_i^M \rho_j^M,
\\
&\label{def:M7}\tag{M7}
\rho_i^M\rho_j^M\rho_i^M +\chi(R) \sigma_i^M \rho_j^M \sigma_i^M
= \rho_j^M\rho_i^M\rho_j^M + \chi(R) \sigma_j^M \rho_i^M \sigma_j^M.
\end{align}
These explicit conditions on the tuple $(M, \sigma^M, \rho^M, S^M, R^M)$ are described as follows.
\begin{thm}[{\cite[Theorem A]{LNX25}}]\label{thmWM}
The action \eqref{def:QWPAction} induces a $B\wr\cH(d)$-module structure on $M \wr N$ if and only if the following are true, for all $b\in B^{\otimes d}, m \in M^{\otimes d}, 1\leq i < d$:
\begin{enumerate}[\textup{(}i\textup{)}]
\item \eqref{def:M2}, \eqref{def:M4} hold;
\item \eqref{def:M5}--\eqref{def:M7} hold additionally for $|i-j|=1$, if $d\geq 3$.
\end{enumerate}
In this case, we call $M$ a {\bf ground} module and $M\wr N$ a {\bf wreath module}.
\end{thm}
}

%=============================================
\subsection{Wreath Index Set and Two Partial Orders} \label{sec:indexset}\label{sec:transpose}
%=============================================
The goal of this section is to lift the notion of Specht modules  over the base algebra $B$, if they exist, to the notion of Specht modules over the quantum wreath product $B\wr\cH(d)$. 

Let $I_B$ be a set. 
Recall that $\Pi_m$ is the set of partitions of $m$. 
Write $\Pi = \bigcup_{m \geq 0} \Pi_m$ and $\#\lambda := \sum_{i} \lambda_i$.
Denote the set of $I_B$-indexed multipartitions of $d$ by
\begin{equation}
\textstyle
\Pi_{I_B}(d) := \left\{
\bblambda: I_B \to \Pi
~\middle|~
\sum_{\nu \in I_B} \#\bblambda(\nu) = d
\right\}.
\end{equation}
Write $\Omega:=  \Pi_{I_B}(d) $ for brevity.
It is useful to express these multipartitions $\bblambda$  as a sum of characteristic functions via
$\bblambda \equiv \sum_{\lambda \in I_B} e_\lambda^{\bblambda(\lambda)}$ where
\begin{equation}\label{def:charfcn}
e_\lambda^\nu(\lambda) := \nu,
\quad
e_\lambda^\nu(\gamma) := \varnothing 
\quad
\textup{for all }
\quad \gamma \in I_B \setminus \{\lambda \}.
%e_\lambda^\nu:  \mu \mapsto \begin{cases}
%\nu &\tif \lambda = \mu;
%\\
%\varnothing &\textup{otherwise}.
%\end{cases}
\end{equation}
For example, if $m = d = 2$, then
\[
\Pi_m = \{(2), (1,1)\}, 
\quad
\Omega = \{ e_{(2)}^{(2)}, \ e_{(2)}^{(1,1)},\ e_{(1,1)}^{(2)},\ e_{(1,1)}^{(1,1)},\ e_{(2)}^{(1)} + e_{(1,1)}^{(1)}\}.
\] 
In particular, $e_{(2)}^{(1)} + e_{(1,1)}^{(1)} = e_{(1,1)}^{(1)} + e_{(2)}^{(1)}$ is the map taking both $(2)$ and $(1,1)$ to the  partition $(1)$.

{
Suppose that $(I_B,\leq)$ is a poset from now on.
The following is a natural poset structure on $\Omega$ which specializes to the usual dominance order on multipartitions by letting $I_B = \{1, \dots, m\}$ be the totally ordered set.}
\begin{Def}[Dominance Order]
\label{def:domPO}
Impose a partial order on $\Omega$ induced from $(I_B,\leq)$  by  the following: $\bblambda \unlhd \bblambda'$ if and only if 
\begin{equation}\label{def:POIA}
\sum_{j \leq k}\bblambda(\nu)_j
+
\sum_{\gamma \in I_B; \gamma > \nu} \#\bblambda(\gamma)  
\ \leq \ 
\sum_{j \leq k}\bblambda'(\nu)_j
+
\sum_{\gamma \in I_B; \gamma > \nu} \#\bblambda'(\gamma)
\qquad
(\textup{for all}
~ \nu\in I_B, 
~k\geq 1).
\end{equation}
\end{Def}
Therefore, for $\lambda, \lambda' \in I_B$ and $\nu, \nu' \in \Pi_d$, the following implications are immediate:
\begin{equation}\label{eq:POimplications}
 \lambda < \lambda' \in I_B \Rightarrow  e_\lambda^\nu \lhd e_{\lambda'}^{\nu'},
\quad
\textup{and}
\quad
\nu \lhd \nu' \in \Pi_d \Rightarrow e_\lambda^\nu \lhd e_{\lambda}^{\nu'},
\end{equation}
where $(\Pi_d, \unlhd)$ is the dominance order.
%------------------------------------------------------------
{The following observation describes certain binary relations in $(\Omega,\unlhd)$ when the supports of two $I_B$-indexed multipartitions have the same cardinality.}
%------------------------------------------------------------
\begin{prop}\label{prop:PO}
Let $\mu = (\mu_1, \dots, \mu_r) \vDash d$,
and let $\lambda^{(1)}, \dots, \lambda^{(r)}, \gamma^{(1)}, \dots,  \gamma^{(r)}$ be elements in $I_B$ such that
$\lambda^{(i)} \neq \lambda^{(j)}$,
$\gamma^{(i)}\neq \gamma^{(j)}$ unless $j=i$.
Suppose that for all $1\leq i \leq r$,
$\lambda^{(i)} \leq \gamma^{(i)}$ in $I_B$,
$\nu^{(i)}, \delta^{(i)} \vdash \mu_i$,
and
\[
\nu^{(i)} \lhd \delta^{(i)}
\quad
\textup{if}
\quad
\lambda^{(i)} = \gamma^{(i)}.
\]
Then, $\sum_i e_{\lambda^{(i)}}^{\nu^{(i)}} \lhd \sum_i e_{\gamma^{(i)}}^{\delta^{(i)}}$ in $\Omega$.
%If there is a $w \in \Sigma_d$ such that $\overline{\bblambda}(i) < \overline{\bbmu}(w(i))$ for all $i$, then $\bblambda < \bbmu$ in $\Omega$.
\end{prop}
%------------------------------------------------------------
\begin{proof}
{We omit the proof as it is elementary and will not be used in later sections.
}
\end{proof}
%------------------------------------------------------------
{
Motivated by the combinatorial description in Proposition~\ref{prop:PO}, 
a new partial order on $\Omega$ is imposed as follows. 
As we will show in Section~\ref{sec:1-cover}, this seemingly artificial partial order in fact arises naturally from the K\"unneth formula.}

%-----------------------------
In order to define this partial order, it is useful to fix an enumeration $\{\lambda^{(1)}, \dots, \lambda^{(s)} \}$ of $I_B$ that is compatible with the partial order $\leq$, i.e., if $\lambda^{(i)}<\lambda^{(j)}$ in $I_B$ then $i<j \in \ZZ$.
Given $\bblambda:I_B\to \Pi$ with support supp$(\bblambda) = \{\lambda^{(k_1)}, \dots, \lambda^{(k_r)}\}$ where $1\leq k_1 < \dots < k_r \leq s$,
define a composition $\mu = \mu(\bblambda) \vDash d$ by $\mu_i := \#\bblambda(\lambda^{(k_i)})$, and a tuple
\begin{equation}\label{eq:bblambdabar}
\overline{\bblambda} := (
\underset{\mu_1\textup{ times}}{\underbrace{\lambda^{(k_1)}, \dots, \lambda^{(k_1)}}}, 
\underset{\mu_2\textup{ times}}{\underbrace{\lambda^{(k_2)}, \dots, \lambda^{(k_2)}}}, 
\dots, 
\underset{\mu_r\textup{ times}}{\underbrace{\lambda^{(k_r)}, \dots, \lambda^{(k_r)}}} 
).
\end{equation}
Equivalently, $\overline{\bblambda}$ is a map $\{1,\dots,d\} \to I_B$ such that $\overline{\bblambda}(i) := \lambda^{(k_j)}$ where $\mu_1 + \dots + \mu_{j-1} < i \leq \mu_1 + \dots + \mu_j$.

%------------------------------------------------------------
\begin{Def}[K\"unneth Order]
\label{def:KPO}
{
Let $\bblambda, \bbmu \in \Omega$. 
Write $\bbmu \uparrow_K \bblambda$ if there exists $w\in \Sigma_d$ such that
\begin{enu}[(K1)]
    \item $\overline{\bbmu}(i) \leq \overline{\bblambda}(w(i)) \in I_B$ for all $1\leq i \leq d$;
    \item For each $\nu \in I_B,$ there is a partition $\beta \vdash n_{\nu}$ that is contained in both $\bblambda(\nu)$ and $\bbmu(\nu)$, 
    where
    $n_{\nu} = n_\nu(w) := \#\{i~|~\overline{\bbmu}(i) = \nu = \overline{\bblambda}(w(i))\}$.
\end{enu}
The {\em K\"unneth order} $\leq_K$ on $\Omega$ is the transitive closure of the relation $\uparrow_K$.
}
\end{Def}
{
Note that $\leq_K$ is indeed a partial order.
It follows from construction that $\leq_K$ is reflexive  and transitive.
For antisymmetry, set 
\begin{equation}
||\bblambda ||:=\sum\nolimits_{i=1}^d j_i,
\quad
\textup{where}
\quad
\overline{\bblambda}(i)=\lambda^{(j_i)}.
\end{equation}
If $\bbmu\uparrow_K\bblambda$ via $w$, then (K1) implies that $||\bbmu||\leq ||\bblambda ||$, with equality only if $\overline{\bbmu}(i)=\overline{\bblambda}(w(i))$ for all $i$. 
In that case, $n_\nu(w)=\#\bbmu(\nu)=\#\bblambda(\nu)$ for all $\nu\in I_B$, and hence (K2) forces $\bbmu=\bblambda$.
That is, $||-||$ strictly increases along any chain of relations $\bbmu\uparrow_K\bblambda$ with $\bbmu\neq\bblambda$, and so $\leq_K$ is antisymmetric.
}

%------------------------------------------------------------
Now, assume that the poset $I_B$ admits the notion of transposition, i.e., an order-reversing involution $t: \nu \mapsto \nu^t$.
Then, define the transpose map on $\Omega$ via $\bblambda \mapsto \bblambda^t$, where
\begin{equation}\label{def:t}
\bblambda^t : I_B \to \Pi,
\quad
\nu \mapsto \bblambda(\nu^t)^t.
\end{equation}
Indeed, $(\bblambda^t)^t = \bblambda$ since 
$
(\bblambda^t)^t(\nu) 
= (\bblambda((\nu^t)^t)^t)^t
= \bblambda(\nu)$
for all $\nu \in I_B$.
%===============================
\subsection{Parabolically Induced Wreath Modules}\label{sec:newSpecht}
%===============================
We are now in a position to define certain induced modules over $B\wr \cH(d)$ with respect to an index $\bblambda$ in the poset $\Omega$.
From now on, we use the following shorthand notations for any composition $\mu = (\mu_1, \dots, \mu_r) \vDash d$:
\begin{align}
\label{eq:ind}
&\ind_\mu^d =\ind_{\mu_1, \mu_2, \dots, \mu_r}^d := \ind_{B \wr \cH(\mu)}^{B \wr \cH(d)},
\\
\label{eq:HEshort}
&
{\Hom_{\mu} :=\Hom_{B\wr \cH(\mu)}
\quad
\textup{and}
\quad
\Ext^n_{\mu} :=\Ext^n_{B\wr \cH(\mu)}
\quad
(n\geq 0).
}
\end{align}

\begin{Def}\label{def:wreathSpecht}
Let $\{S^\lambda_B ~|~ \lambda \in I_B\}$ be a family of ground $B$-modules with associated $\chi_\lambda$, labeled by a poset $I_B$ such that $\End_B(S_B^\lambda) \cong K$.
    Denote the support of $\bblambda \in \Omega$ by $\{\lambda^{(1)}, \lambda^{(2)}, \dots, \lambda^{(r)}\} \subseteq I_B$.
Write $\mu_i := \# \bblambda(\lambda^{(i)})$, and hence $\mu = (\mu_1, \dots, \mu_r) \vDash d$. 
For $\nu \vdash k \leq d$, let $S_{\chi_\lambda}^{\nu}$ be the Specht module of the twisted Hecke algebra $\tcH_{k}^{\chi_\lambda}$ (see \cite[Section 3.3]{LNX25}).
Since each $S^{\lambda^{(i)}}_B \wr S_{\chi_{\lambda^{(i)}}}^{\bblambda(\lambda^{(i)})}$ is a module over $B \wr\cH(\mu_i)$,
their tensor product is a module over the parabolic subalgebra $B\wr \cH(\mu)$,
which induces to  
\[
S^\bblambda = 
\ind_{\mu}^{d} \big( \bigotimes\nolimits_{i=1}^r (S_{B}^{\lambda^{(i)}} \wr S_{\chi_{\lambda^{(i)}}}^{\bblambda(\lambda^{(i)})}) \big).
\]
\end{Def}
When $S^\lambda_B$ are the Specht modules over $B$, we call $S^\bblambda$ the Specht modules over $B \wr \cH(d)$.
%-----------------------------
\begin{prop}\label{prop:indexact}
    The parabolic induction functor $\ind_\mu^d: B\wr \cH(\mu)$-mod $\to B\wr\cH(d)$-mod is exact.
\end{prop}
%-----------------------------
\begin{proof}
 Note that $B \wr \cH(d)$ is a free right $B \wr \cH(\mu)$-module with a free basis 
 $\{ H_w  ~|~ w \in \Sigma^\mu \}$ 
 where $\Sigma^\mu$ is the set of minimal length coset representatives for $\Sigma_d / \Sigma_\mu$.
Since tensoring with a free module preserves all short exact sequences, the parabolic induction functor is exact.
\end{proof}

%=============================================
\subsection{Mackey Theory for Quantum Wreath Products}\label{sec:Mackey}
%=============================================
Throughout this article, the notation of the restriction functor is suppressed.
That is, a $B\wr\cH(d)$-module $M$ is regarded as a $B \wr \cH(\mu)$-module without further comment.

In this section, we develop Mackey theory for quantum wreath products admitting PBW bases.
Assume further that $\sigma$ is the flip map.
These results will be used in Section \ref{sec:DSgHu} for the verification of (S1).

Let $\mu = (\mu_1, \dots, \mu_r) \vDash d$ and $\eta = (\eta_1, \dots, \eta_{r'}) \vDash d$ be compositions of $d$, 
with corresponding Young subgroups $\Sigma_\mu, \Sigma_\eta \le \Sigma_d$.
The {\em $i$th interval of} $\mu$ is the set
\begin{equation}\label{def:interval}
    \ZZ^\mu_i := \{\mu_1 + \dots + \mu_{i-1} +k ~|~ 1\leq k \leq \mu_i\}.
\end{equation}
Let ${}^\mu\Sigma^\eta$ denote the set of minimal length double coset representatives in $\Sigma_\mu\backslash\Sigma_d/\Sigma_\eta$.
For each $w \in {}^\mu\Sigma^\eta$, it is well-known that the intersection $\Sigma_\mu \cap {}^w \Sigma_\eta = \Sigma_\mu \cap w \Sigma_\eta w^{-1}$ is a Young subgroup $\Sigma_{\delta(\mu,w,\eta)}$ corresponding to the composition $\delta(\mu,w,\eta) \vDash d$ whose entries are given by the nonzero sizes 
\[
\#(
\ZZ^\mu_i
\cap 
w\ZZ^\eta_j)
\]
following this order: 
$(i,j) = (1,1), (1,2)$, $\dots$, $(1,r')$, $(2,1)$, $(2,2)$, $\dots$, $(r,r')$.
For $w \in \Sigma_d$, denote the place permutation automorphism on $\cB:= B^{\otimes d}$ by
\[
w(b_1 \otimes \dots \otimes b_d) := b_{w^{-1}(1)} \otimes \dots \otimes b_{w^{-1}(d)}.
\]
For $w \in {}^\mu\Sigma^\eta$ and $Y \in B\wr\cH(\eta)$-mod, let ${}^w Y$ be the $B\wr\cH(\delta(\mu,w,\eta))$-module with underlying space $Y$ via the new action $\circ$ given by, for $b \in \cB$ and $y \in Y$:
\[
b \circ y := w^{-1}(b) \cdot y, 
\quad
H_i \circ y := H_{w^{-1}(i)} \cdot y.
\]

\begin{prop}[Mackey Filtration]\label{prop:Mackey}
Let $\mu, \eta \vDash d$, and let $Y \in B\wr\cH(\eta)$-mod.
Then $\ind_\eta^d Y$ admits a filtration as a $B\wr\cH(\mu)$-module:
\[
\{0\} = V_0 \subseteq V_1 \subseteq \dots \subseteq V_k = %\res^d_\mu
\ind_\eta^d Y, \qquad (k = |{}^\mu\Sigma^\eta|),
\]
whose subquotients are given by
\begin{equation}\label{eq:MackeySubquot}
V_i / V_{i-1} \cong 
\ind_{\delta(\mu,w_i,\eta)}^{\mu} ( 
{}^{w_i} 
%\big(\res^{\eta}_{\delta(\eta, w_i^{-1}, \mu)} 
Y 
%\big) 
), \qquad w_i \in {}^\mu\Sigma^\eta.
\end{equation}
In particular, as a $\cB$-module, $\ind_\eta^d Y$ has a filtration with subquotients isomorphic to ${}^w Y$ for $w \in \Sigma^\eta$.
\end{prop}

\begin{proof}
By the PBW theorem for quantum wreath products (see \eqref{def:PBWbases}), $B\wr\cH(d)$ is free as a left $B\wr\cH(\mu)$-module with basis $\{H_w \mid w \in {}^\mu\Sigma\}$
and as a right $B\wr\cH(\eta)$-module with basis $\{H_w \mid w \in \Sigma^\eta\}$. 
Decomposing $\Sigma_d = \bigsqcup_{w \in {}^\mu\Sigma^\eta} \Sigma_\mu w \Sigma_\eta$, 
each double coset defines a $(B\wr\cH(\mu))$--$(B\wr\cH(\eta))$-bimodule 
\[
(B\wr\cH(\mu)) H_w (B\wr\cH(\eta)) \subseteq B \wr \cH(d).
\]
Filtering $B\wr\cH(d)$ compatibly with the Bruhat order on ${}^\mu\Sigma^\eta$ and tensoring with $Y$ yields the filtration with subquotients \eqref{eq:MackeySubquot}.
\end{proof}
For any $w \in {}^\mu \Sigma^\eta$, write $\delta_w := \delta(\mu,w,\eta)$.
{The corollary below allows us to detect non-vanishing homomorphisms and extensions between parabolically induced modules,
under certain conditions on the parameters.
}
\begin{cor}\label{cor:MackeyHomExt}
{Suppose that $S= 0, \rho = 0$ and $R$ is invertible.
Let $X$ be a $B\wr\cH(\mu)$-module and $Y$ a $B\wr\cH(\eta)$-module.
\begin{enua}
\item 
One has 
\[
\dim \Hom_{B\wr\cH(d)}(\ind_\mu^d X, \ind_\eta^d Y) \leq \sum\nolimits_{w \in {}^\mu\Sigma^\eta} \dim \Hom_{B\wr\cH(\delta_w)}(X, {}^wY).
\]
\item If $\Ext^1_{B\wr\cH(d)}(\ind_\mu^d X, \ind_\eta^d Y) \neq 0$, then $\Ext^1_{B\wr\cH(\delta_w)}(X, {}^wY) \neq 0$ for some $w \in {}^\mu\Sigma^\eta$.
\end{enua}
}
\end{cor}
\begin{proof}
The assumptions on the parameters $(S, R, \sigma, \rho)$ imply that each $H_x$ $(x \in \Sigma_\mu)$ is invertible and
\[
H_x H_y \in \cB^\times H_{xy}
\quad
\textup{for all}
\quad
x, y \in \Sigma_\mu.
\]
In other words, $B \wr \cH(\mu)$ and $B \wr \cH(\delta_w)$ are the crossed products $\cB *\Sigma_\mu$ and $\cB * \Sigma_{\delta_w}$, respectively.
Hence,  
$B\wr \cH(\mu)$ is a free Frobenius extension of 
$B\wr\cH(\delta_w)$,
so $\ind^\mu_{\delta_w}$ is also right adjoint to restriction.
Therefore,
Frobenius reciprocity yields
\begin{equation}
    \label{eq:XZ}
 \Ext^n_{\mu}(X,\ind_{\delta_w}^\mu Z) \cong \Ext^n_{\delta_w}(X,Z),
\end{equation}
for any $w \in {}^\mu\Sigma^\eta$, $B\wr\cH(\delta_w)$-module $Z$, and $n \geq 0$.

Next, by Proposition \ref{prop:indexact} and Frobenius reciprocity,
\[
\Ext^n_{B\wr \cH(d)}(\ind_\mu^d X, \ind_\eta^d Y) 
\cong
\Ext^n_{\mu}(X, \ind_\eta^d Y), 
\]
for all $n \geq 0$.
Recall the filtration $(V_i)_i$ of $\ind_\eta^d Y$ in Proposition \ref{prop:Mackey}.
Let $n=0$.
Applying $\Hom_{\mu}(X, -)$ to each exact sequence $0 \to V_{i-1} \to V_i \to V_i/V_{i-1} \to 0$, one gets
\[
\dim \Hom_{\mu}(X,V_i) 
-\dim\Hom_{\mu}(X,V_{i-1})
\leq 
\dim\Hom_{\mu}(X,V_i/V_{i-1})
\quad
\textup{for all}
\quad
i,
\]
which sums up to 
$\dim \Hom_{\mu}(X, \ind_\eta^d Y) 
\leq 
\sum_i \dim\Hom_{\mu}(X,V_i/V_{i-1})$.
By Proposition \ref{prop:Mackey},
 $V_i/V_{i-1} \cong \ind^\mu_{\delta_{w_i}} ({}^{w_i}Y)$ for some $w_i$
 where the $w_i$ enumerate ${}^\mu\Sigma^\eta$.
 Hence,
 part (a) follows from \eqref{eq:XZ}.
 
For (b), consider $n=1$. One obtains the following exact sequences 
\[\Ext^1_{\mu}(X,V_{i-1}) 
\to \Ext^1_{\mu}(X,V_{i})
\to \Ext^1_{\mu}(X,V_i/V_{i-1}) \]
 for all $i$. 
By a routine induction on $i$, one can show that
if $\Ext^1_{\mu}(X, V_i/V_{i-1}) = 0$ for all $i$, then $\Ext^1_\mu(X, \ind_\eta^d Y) = 0$.
In other words, (b) follows from \eqref{eq:XZ}.
The corollary is proved.
\end{proof}
%=============================================
\section{Duality for Wreath Modules}\label{sec:duality}
%=============================================
{
In this section, sufficient conditions to introduce an automorphism $\#$ and an anti-automorphism $\ast$ on the quantum wreath product $B \wr \cH(d)$ are determined.
These maps should be considered as a generalization of the maps $\#: T_i \mapsto -qT_i^{-1}$ and $\ast: T_i \mapsto T_i$ for the Hecke algebras.
In particular,  there is an isomorphism
\[
(S^\bblambda)^\# \cong (S^{\bblambda^t})^* \in B\wr \cH(d)\hp \mMod
\]
between these parabolically induced wreath modules.
This duality will play a crucial role in establishing a dual Specht theory for the Hu algebras.
}
%=============================================
\subsection{Duality between \texorpdfstring{$\#$}{sharp} and \texorpdfstring{$\ast$}{*} for Hecke algebras}\label{sec:dualSpecht}
%=============================================
Recall that the Hecke algebra $\cH_q(\Sigma_m)$ has an automorphism $\#$ and an anti-automorphism $*$ determined by $T_i^\# := -q T_i\inv$ and $T_i^* := T_i$, respectively. 
Therefore, for any $w \in \Sigma_m$:
\begin{equation}\label{def:mapH}
T_w^\# = (-q)^{\ell(w)} T_{w\inv}\inv,
\quad
\textup{and}
\quad
T_w^* = T_{w\inv}.
\end{equation}
Let $M$ be a finite-dimensional $\cH_q(\Sigma_m)$-module on which the action is denoted by $b\cdot m$ for $b\in \cH_q(\Sigma_m)$, $m\in M$.
There are two associated modules $M^\#$ and $M^*$
whose underlying spaces are $M$ and $\Hom_K(M,K)$, respectively.
The actions $\circ$ of $b\in \cH_q(\Sigma_m)$ on these modules are given by, respectively, 
\begin{equation}
b\circ m := b^\#\cdot m,
\quad
(b \circ f)(m) := f(b^*\cdot m),
\quad
\textup{for}
\quad m\in M, \ f\in M^*.
\end{equation}
The following isomorphism is well-known (see \cite[Exercise 3.14]{Ma99}):
\begin{equation}\label{eq:dualSpechtHecke}
(S^\lambda_m)^\# \cong (S^{\lambda^t}_m)^*,
\quad
\textup{for any}
\quad
\lambda \vdash m,
\end{equation} 
where $S^\lambda_m$ is the Specht module of $\cH_q(\Sigma_m)$, 
and $(S^{\lambda}_m)^*$ is the dual Specht module with respect to $\lambda$.

Next, we consider the twisted Hecke algebra $\tcH^\chi_d$. 
%with $\chi(S) = 0$ and $\chi(R) = \sqrt{f_\lambda}^2$ for some $\sqrt{f_\lambda} \in K^\times$.
%Therefore, 
One can directly check that there is an automorphism $\#$ and an anti-automorphism $*$ determined by 
$h_i^\# := -\chi(R)h_i^{-1} = \chi(S) - h_i$ and $h_i^* := h_i$, respectively. Therefore,
\begin{equation}\label{def:maptcH}
h_w^\# = (-\chi(R))^{\ell(w)} h_{w^{-1}}^{-1},
\quad
\textup{and}
\quad
h_w^* = h_{w\inv},
\quad
(w \in \Sigma_d).
\end{equation}
{
One can also verify the following analog of \eqref{eq:dualSpechtHecke} for these twisted Hecke algebras.
\begin{lem}\label{lem:tHAdual}
%Suppose that $\{\chi_\lambda~|~\lambda \in \Pi_m\}$ is a collection of characters such that ${\chi_\lambda} = \chi_{\lambda^t}$ for all $\lambda \vdash m$.
Let $\tcH^{\chi_\lambda}_d$ be the twisted Hecke algebra with quadratic relations of the form $h_i^2 = \chi_\lambda(S)h_i + \chi_\lambda(R)$.
If $\chi_\lambda(R)$ is nonzero and $h_i^2 - \chi_\lambda(S)h_i - \chi_\lambda(R)$ splits into linear factors,
then, as $\tcH^{\chi_\lambda}_d$-modules:
\begin{equation}\label{eq:tcHdS}
(S^\nu_{\chi_\lambda})^\# \cong (S^{\nu^t}_{\chi_\lambda})^*
\quad
\textup{for all}
\quad
\nu\vdash d.
%(K^+)^* \cong  K^+,
%\quad
%(K^-)^* \cong  K^-,
%\quad
%(K^+)^\# \cong  K^-,
%\quad
%(K^-)^\# \cong  K^+.
\end{equation}
\end{lem}
}
%=============================================
\subsection{The Automorphism \texorpdfstring{\#}{sharp}}\label{sec:auto}
%=============================================
%Denote the Specht modules over $B \wr \cH(d)$ by $S^\bblambda$ for $\bblambda \in \Omega$, 
%and those over $B$ by $S_B^\lambda$ for $\lambda \in I_B$.
Consider a base algebra $B$ with an automorphism $\#$.
Within this section, assume
\begin{equation}\label{eq:weakFrob}    
\sigma =\textup{flip},
\quad
\sigma(S) = S,
\quad
\sigma(R) = R,
\quad
S b = \sigma(b) S
\quad
\textup{for all}
\quad b\in B\otimes B.
\end{equation}
These assumptions \eqref{eq:weakFrob} arise naturally in \cite[(A2--A3)]{LM26} in the context of polynomial quantum wreath products. 
In particular, \eqref{eq:weakFrob} is satisfied for the Hu algebras and the Yokonuma--Hecke algebras.

Let $A$ be the free associative algebra generated by $B^{\otimes d}$ and $H_i$ $(1\leq i \leq d-1)$.
Let $\theta\in Z(B\otimes B)^{\times}$ be an element satisfying 
$\rho(\theta)=0$.
Consider the algebra homomorphism $\# = \#(\theta)$ uniquely determined by the following:
\begin{equation}\label{def:autA}
\#: A \to B \wr \cH(d), \quad 
b_1 \otimes \dots \otimes b_d \mapsto b_1^{\#} \otimes \dots \otimes b_d^{\#}, \quad
H_i \mapsto (- H_i + S_i){\theta_i}.
\end{equation}
Consider the following conditions:
\begin{align}
\label{eq:sharpB}
&R^\# = %\sigma(R) 
R{\sigma(\theta)\theta },
\quad
S^\# = S{\theta},
\\ \label{eq:sharpW}
&\sigma(b)^\# = \sigma(b^\#), 
\quad
\rho(b)^\# =  -  \rho(b^\#){\theta},
\\ \label{eq:sharptheta}
&{\sigma_i\sigma_j(\theta_i) \sigma_i(\theta_j)\theta_i
=\sigma_j\sigma_i(\theta_j)\sigma_j(\theta_i)\theta_j},
\quad\textup{for all}
\quad |i-j| = 1,
\\
\label{eq:sharpBB}
&S_i S_j S_i + \sigma_i(S_j) R_i = S_j S_i S_j + \sigma_j(S_i) R_j,
\quad\textup{for all}
\quad |i-j| = 1.
\end{align}
\begin{thm}\label{thm:sharp}
Suppose that \eqref{eq:weakFrob} holds.
The map $\#$ in \eqref{def:autA} induces a homomorphism $B \wr \cH(d) \to B \wr \cH(d)$ if and only if \eqref{eq:sharpB}--\eqref{eq:sharpBB} hold.
In this case, $\#$ is an automorphism on $B \wr \cH(d)$.
\end{thm}
%-----------------
\begin{proof}
To show that $\#$ is a homomorphism, it suffices to prove that \eqref{eq:sharpB}--\eqref{eq:sharpBB} are equivalent to the following:
\begin{align}
\label{equivrel1}
&(H^2 - S H - R)^\# = 0,
\quad
(H b - \sigma(b) H - \rho(b))^\# = 0, 
\\
\label{equivrel2}
&(H_i H_j H_i - H_j H_i H_j)^\# = 0
\quad\textup{for all}
\quad |i-j| = 1,
\\
\label{equivrel3}
&(H_i H_j - H_j H_i)^\# = 0
\quad\textup{for all}\quad
|i-j| >1.
\end{align}
Details of the equivalences are provided in Section \ref{app:equiv}.
The verification of the adjacent braid relation \eqref{equivrel2} is surprisingly non-trivial.

It remains to show that $\#$ is bijective. 
First, filter $B\wr \cH(d)$ by the length in $\Sigma_d$. 
Set
\begin{equation}\label{eq:filtrationBHd}
F_p := \operatorname{Span}_K \{ b H_w \mid b \in B^{\otimes d}, \ 
w \in \Sigma_d, \ \ell(w) \leq p \}, \qquad (p \geq 0).
\end{equation}
For any $w \in \Sigma_d$ with a reduced expression $w = s_{i_1} \dots s_{i_l}$, 
\[
H_w^\# = (-H_{i_1} + S_{i_1}){\theta_{i_1}} \dots (-H_{i_l} + S_{i_l}){\theta_{i_l}} \in  {Z(B^{\otimes d})^\times} H_w + F_{\ell(w)-1},
\]
and hence for an arbitrary basis element, one has
$(b H_w)^\# \in {\zeta_w} b^\# H_w  + F_{\ell(w)-1}$
{for some $\zeta_w \in Z(B^{\otimes d})^\times$}.
Since $\#$ is an automorphism of $B$, the associated graded map 
\[
\operatorname{gr}(\#): \operatorname{gr}(B\wr \cH(d)) \to \operatorname{gr}(B\wr \cH(d)), \qquad \overline{b H_w} \mapsto  \overline{ {\zeta_w} b^\# H_w}
\]
sends the PBW basis to a basis. Hence, $\#$ is a $K$-module isomorphism. 
\end{proof}
\begin{rmk}\label{rmk:sharpHu}
If $\sigma = \textup{flip}$, $S=0$ and $\rho=0$,
the conditions \eqref{eq:sharpB}--\eqref{eq:sharpBB} can be reduced to a single condition 
\begin{eq}
R^\# = R \theta \sigma(\theta),
\end{eq}
since the condition $R = \sigma(R)$ follows from $S=0$:
\[
0 = H^3 - H^3 = H R - RH = (\sigma(R) - R)H + \rho(R).
\]
Thus, the assumptions $\rho(\theta)=0$ and \eqref{eq:weakFrob} are automatically true. 
This will be used in the proof of Proposition \ref{prop:genSpechtiso}.
\end{rmk}

Both $\cH_q(\Sigma_d)$ and $\tcH^\chi_d$ can be regarded as quantum wreath products of the form $K \wr \cH(d)$ with $\rho =0$. 
In this case, \eqref{eq:sharpB}--\eqref{eq:sharpBB} automatically hold with $\theta =1$, and hence the maps $\#$ are well-defined in both cases.  
Then, such automorphisms $\#$ on $K \wr \cH(d)$ agree with the map $\#$ on $\cH_q(\Sigma_d)$ or on $\tcH^\chi_d$ defined in \eqref{def:mapH} and \eqref{def:maptcH}, respectively, since $T_i^\# = -qT_i\inv = -T_i + (q-1)$ and
 $h_i^\# = \chi(S)-h_i$.    
%==================
\subsection{Conditions on the Wreath Functor \texorpdfstring{$M^\#\wr-$}{M\# wr -}}
%=================
Having established the conditions under which the map $\#$ defines an automorphism on the algebra $B \wr \cH(d)$, we now turn our attention to its representations. 
{
Recall that we assumed $\theta \in Z(B\otimes B)^\times$, and hence
it acts on $M \otimes M$ by a nonzero scalar (provided $\End_B(M) \cong K$).
Denote this scalar by $\chi(\theta)$ in accordance with \eqref{eq:SMRM}, where
$\chi: B\otimes B \to K$ is a map satisfying 
$S^M(m) = \chi(S)m$ and $R^M(m) = \chi(R)m$ for all $m \in M\otimes M$.
The conditions \eqref{eq:sharpB} suggest that $S^\#$ and $R^\#$ should act 
on $M\otimes M$ by the scalars $\chi(S)\chi(\theta)$ and $\chi(R)\chi(\theta)^2$, respectively.

Equivalently, if there is a ground module structure for $M^\#$, 
the parameters $S^{M^\#}$ and $R^{M^\#}$ must be given by $S^{M^\#} = l_{S^\#}^M$ and 
$R^{M^\#} = l_{R^\#}^M.$
The remaining parameters $\sigma^{M^\#}$ and $\rho^{M^\#}$ are constructed 
by a careful analysis of \eqref{def:M2}--\eqref{def:M7}.
}
The following lemma provides sufficient conditions for $M^\#$ to be a ground module for a given ground module $M$ {under two scenarios.
One requires $S = 0$ and is the case that is relevant to this paper.
The other scenario, which is useful in studying the representation theory of Yokonuma--Hecke algebras, is recorded for completeness.
}
\begin{lem}\label{lem:Msharpground}
{
Suppose that \eqref{eq:weakFrob} holds, and that
$(M, \sigma^M, \rho^M, S^M, R^M)$ is a ground module satisfying $\End_B(M) \cong K$.
Let 
\[
\sigma^{M^\#} = \sigma^M,
\quad
\rho^{M^\#} = l_{S{\theta}}^M -{l^M_\theta}\rho^{M},
\quad
S^{M^\#} = l_{S^\#}^M,
\quad
R^{M^\#} = l_{R^\#}^M.
\]
Then $(M^\#, \sigma^{M^\#}, \rho^{M^\#}, S^{M^\#}, R^{M^\#})$ is also a ground module if
either $S=0$ or 
%\begin{equation}\label{eq:Msharpground}
$\sigma^M = 0$.
%\end{equation}
}
\end{lem}
\begin{proof}
Abbreviate by (Mi)\# the corresponding condition (Mi) for $(M^\#, \sigma^{M^\#}, \rho^{M^\#}, S^{M^\#}, R^{M^\#})$.
By Theorem \ref{thmWM}, it suffices to show that \eqref{def:M2}--\eqref{def:M7}  imply \eqref{def:M2}\#--\eqref{def:M7}\#.
The details of the verifications will be provided in Section \ref{app:Msharpground}.
\end{proof}
In the following theorem, it is shown that taking the automorphism $\#$ is compatible with the wreath module construction. 
{
Note first that the twist $\theta$ rescales the parameters. 
Let $\chi^\#$ be the map determined by
\[
\chi^\#(S) := \chi(S^\#) = \chi(S)\chi(\theta),
\quad
\chi^\#(R) := \chi(R^\#) = \chi(R)\chi(\theta)^2.
\]
Thus, the assignment $h_i \mapsto h_i^\# := \chi(\theta)(-h_i + \chi(S))$ defines an isomorphism $\tcH^{\chi^\#}_d \to \tcH^{\chi}_d$.
Accordingly, for an $\tcH^\chi_d$-module $N$, write $N^\#$ for the $\tcH^{\chi^\#}_d$-module whose underlying space is $N$ and whose action is given by $h_i \circ n := h_i^\# \cdot n$.
}
\begin{theorem}\label{thm:MNsharp}
{
Suppose that \eqref{eq:weakFrob} holds.
Let $M, M^\#$ be ground modules as described in Lemma \ref{lem:Msharpground} such that $\End_B(M) \cong K$, and {$\sigma^M = 0$ unless $S = 0$}.
Let $N$ be a module over the twisted Hecke algebra $\tcH^\chi_d$ with quadratic relations of the form $h_i^2 = \chi(S) h_i + \chi(R)$.
Then, $(M \wr N)^\# \cong M^\# \wr N^\#$ as $B \wr \cH(d)$-modules.
}
\end{theorem}
\begin{proof}
{
    For $m \otimes n \in M \wr N$, recall that the actions of $H_i$ and $b\in B^{\otimes d}$ are given by, respectively:
    \[
    b \cdot (m\otimes n) = (b \cdot m) \otimes n,
    \quad
    H_i \cdot (m\otimes n) = \sigma^M_i(m) \otimes (h_i\cdot n) + \rho^M_i(m) \otimes n.
    \]
    Hence, the corresponding actions $\circ$ on $(M\wr N)^\#$ are given by, respectively:
    \begin{equation}
    \begin{split}
    b \circ (m\otimes n) = (b^\# \cdot m) \otimes n,
    \quad
    H_i &\circ (m\otimes n) 
    =(S_i - H_i){\theta_i} \cdot (m\otimes n)
    \\
    &= (l^M_{S_i{\theta_i} } - {l^M_{\theta_i}} \rho^M_i)(m) \otimes n -  {\chi(\theta)}\sigma^M_i(m) \otimes (h_i \cdot n).
    \end{split}    
    \end{equation}
Since $l^M_{S_i{\theta_i} } - {l^M_{\theta_i}} \rho^M_i = \rho^{M^\#}_i$, 
$\sigma^M_i = \sigma^{M^\#}_i$,
and $h_i^\# = {\chi(\theta)}({\chi(S)}-h_i)$, the two actions agree as long as $\chi(\theta)\chi(S)\sigma^M_i(m)\otimes n = 0$, which follows from the assumption that $\sigma^M = 0$ unless $S=0$.
}
\end{proof}
%=============================================
\subsection{The Anti-Automorphism *}\label{sec:anti-auto}
%=============================================
%Denote the Specht modules over $B \wr \cH(d)$ by $S^\bblambda$ for $\bblambda \in \Omega$, 
%and those over $B$ by $S_B^\lambda$ for $\lambda \in I_B$.
Assume for now that $B$ has an anti-automorphism $*$.
%Let $A$ be the free associative algebra generated by $B^{\otimes d}$ and $H_i$. 
Consider the algebra anti-homomorphism $*$ uniquely determined by the following:
\begin{equation}\label{def:antiautA}
\ast: A \to B \wr \cH(d), \quad 
b_1 \otimes \dots \otimes b_d \mapsto b_1^{\ast} \otimes \dots \otimes b_d^{\ast}, \quad
H_i \mapsto H_i.
\end{equation}
We consider the following conditions:
\begin{align}\label{eq:starB}
&\sigma(S^\ast) = S,
\quad
\rho(S^\ast) + R^\ast = R,
\\
\label{eq:starW}
&\sigma(\sigma(b)^\ast) = b^\ast,
\quad
\rho(\sigma(b)^\ast) + \rho(b)^\ast = 0.
\end{align}
\begin{thm}\label{thm:*}
The map $\ast$ in \eqref{def:antiautA} induces an anti-homomorphism $B \wr \cH(d) \to B \wr \cH(d)$ if and only if \eqref{eq:starB}--\eqref{eq:starW} hold.
In this case, $\ast$ is an anti-automorphism on $B \wr \cH(d)$.
\end{thm}

\begin{proof}
To show that $\ast$ is an anti-homomorphism, it suffices to prove that \eqref{eq:starB}--\eqref{eq:starW} are equivalent to the following: 
\begin{align}
\label{eq:*1}
&(H^2 - S H - R)^\ast = 0, 
\quad
(H b - \sigma(b) H - \rho(b))^\ast = 0, 
\\
\label{eq:*2}
&(H_i H_j H_i - H_j H_i H_j)^\ast = 0
\quad\textup{for all}
\quad |i-j| = 1,
\\
\label{eq:*3}
&(H_i H_j - H_j H_i)^\ast = 0
\quad\textup{for all}\quad
|i-j| >1.
\end{align}
In contrast to the proof of Theorem \ref{thm:sharp}, here the braid relations \eqref{eq:*2}--\eqref{eq:*3} hold immediately. 
The equivalence for \eqref{eq:*1} is deduced from the following equalities:
\[ \begin{split}
(H^2 - S H - R)^\ast 
%&=  (H^\ast)^2 - H^\ast S^\ast - R^\ast 
%\\
&=  H^2 - H S^\ast  - R^\ast 
\\
&=  (S H + R) - (\sigma(S^\ast) H + \rho(S^\ast))  - R^\ast 
\\
&=  (S - \sigma(S^\ast)) H + R - R^\ast - \rho(S^\ast),
\end{split}
\]
\[ \begin{split}
(Hb - \sigma(b)H - \rho(b))^* &= b^* H 
- H \sigma(b)^* 
   -\rho(b)^* 
   \\
&= (b^*- \sigma(\sigma(b)^*) ) H 
- (\rho(\sigma(b)^*)+\rho(b)^*).
\end{split}
\]
Recall the filtration on $B\wr \cH(d)$ from \eqref{eq:filtrationBHd}.
Since 
\[
(b H_w)^* 
= H_w^* b^* 
= H_{w^{-1}} b^* 
\in \sigma_{w^{-1}}(b^*) H_{w^{-1}} + F_{\ell(w)-1},
\]
the bijectivity of $\ast$ can be proved by an argument similar to that for the bijectivity of $\#$.
\end{proof}  
Again, both $\cH_q(\Sigma_d)$ and $\tcH^\chi_d$ can be regarded as quantum wreath products of the form $K \wr \cH(d)$ with $\rho =0$. 
In this case, \eqref{eq:starB}--\eqref{eq:starW} automatically hold, and hence the maps $*$ are well-defined in both cases.  
Then, such anti-automorphisms $*$ on $K \wr \cH(d)$ agree with the map $*$ on $\cH_q(\Sigma_d)$ or on $\tcH^\chi_d$ defined in \eqref{def:mapH} and \eqref{def:maptcH}, respectively, since $T_i^* = T_i$ and
 $h_i^* = h_i$.
%=====================
\subsection{Conditions on the Wreath Functor \texorpdfstring{$M^\ast\wr-$}{M* wr -}}
%=================
 For finite-dimensional $K$-spaces $U$ and $V$, one has the
canonical identification $\Phi: \Hom_K(U,K) \otimes \Hom_K(V,K) \to \Hom_K(U\otimes V,K)$ given by $\Phi(f \otimes g)(u \otimes v) = f(u)g(v)$ for all $u \in U, v\in V$.
 {
 The next lemma specifies exactly when the dual $B$-module $M^*$ inherits a ground module structure.
 }
\begin{lem}\label{lem:Mstarground}
{
Suppose that \eqref{eq:starB}--\eqref{eq:starW} hold, and
$(M, \sigma^M, \rho^M, S^M, R^M)$ is a finite-dimensional ground module.
Let 
\[
\sigma^{M^\ast}(f) = f\sigma^M,
\quad
\rho^{M^\ast}(f) = f\rho^M,
\quad
S^{M^\ast} = l^{M^\ast}_{\chi(S)\id},
\quad
R^{M^\ast}= l^{M^\ast}_{\chi(R)\id}
\]
for all $f \in \Hom_K(M^{\otimes 2},K)$.
Then, $(M^\ast, \sigma^{M^\ast}, \rho^{M^\ast}, S^{M^\ast}, R^{M^\ast})$ is also a ground module if $\sigma^2 = \id$.
}
\end{lem}
\begin{proof}
{Abbreviate by (Mi)* the corresponding condition (Mi) for $(M^\ast, \sigma^{M^\ast}, \rho^{M^\ast}, S^{M^\ast}, R^{M^\ast})$.
By Theorem \ref{thmWM}, it suffices to show that \eqref{def:M2}*--\eqref{def:M7}* follow from \eqref{def:M2}--\eqref{def:M7}.
The details of the verifications will be provided in Section \ref{app:Mstarground}.
}
\end{proof}
The following theorem shows that taking the anti-automorphism $\ast$ is compatible  with the wreath module construction. 
{
\begin{theorem}\label{thm:MNstar}
Suppose that \eqref{eq:starB}--\eqref{eq:starW} hold,
and $\sigma^2 = \id$.
Let $M, M^\ast$ be finite-dimensional ground modules as described in Lemma \ref{lem:Mstarground}, and 
let $N$ be a finite-dimensional module over the twisted Hecke algebra $\tcH^{\chi}_d$ with quadratic relations of the form $h_i^2 = \chi(S) h_i + \chi(R)$.
Then, $(M \wr N)^\ast \cong M^\ast \wr N^\ast$ as $B \wr \cH(d)$-modules.
\end{theorem}
\begin{proof}
Pick an element in $(M\wr N)^\ast$ of the form $f\otimes g$ where $f \in \Hom_K(M^{\otimes d},K)$ and $g \in \Hom_K(N,K)$.
Let $m \otimes n \in M \wr N$. 
Thus, for any $b\in B^{\otimes d}$, the $b$-action on $f\otimes g \in (M\wr N)^\ast$ is given by
\[
b \circ (f\otimes g) : M \wr N \to K,
\quad
m\otimes n 
    \mapsto 
(f\otimes g)(b^\ast \cdot (m \otimes n))
    = f(b^\ast \cdot m) g( n) = (b \circ f)(m) g(n),
\]
which agrees with the evaluation of the map $b \circ (f\otimes g) \in M^\ast \wr N^\ast$ at $m\otimes n$.
Next, consider the $H_i$-action on $f\otimes g \in (M\wr N)^\ast$ given by
    \begin{equation}
    \begin{split}
    H_i \circ (f\otimes g) : M \wr N \to K,
    \quad 
    (m\otimes n) 
    &\mapsto (f\otimes g)(H_i^\ast \cdot (m \otimes n))
    \\
    &= (f\otimes g)(\sigma^M_i(m)\otimes(h_i\cdot n) + \rho^M_i(m)\otimes n)
    \\
    &= f\sigma^M_i(m) g(h_i \cdot n) + f\rho^M_i(m)g(n),
    \end{split}    
    \end{equation}
which agrees with the evaluation of the map $H_i \circ (f\otimes g) \in M^\ast \wr N^\ast$ at $m\otimes n$
since $f\sigma^M_i = \sigma_i^{M^\ast}(f)$ and $f\rho^M_i = \rho_i^{M^\ast}(f)$.
\end{proof}
}
%=======================================
\subsection{Duality between \texorpdfstring{$\#$}{sharp} and \texorpdfstring{$\ast$}{*} for quantum wreath products} \label{sec:newdSpecht}
%=======================================
{
Let 
$\{S^\lambda_B ~|~ \lambda \in I_B\}$ be a collection of finite-dimensional ground $B$-modules such that
\begin{equation}\label{eq:Biso}
\End_B(S_B^\lambda) \cong K,
\quad
(S^\lambda_B)^\# \cong (S^{\lambda^t}_B)^*
\quad
\textup{as ground modules (for all }
\lambda\in I_B).
\end{equation}
Let $\tcH^{\chi_\lambda}_d$ be the twisted Hecke algebra corresponding to the ground module $S^\lambda_B$.
Assume that the duality isomorphism holds for Specht modules over $\tcH^{\chi_\lambda}_d$:
\begin{equation}\label{eq:Hiso}
(S^\nu_{\chi_\lambda})^\# \cong (S^{\nu^t}_{\chi_\lambda})^\ast
\quad
\textup{for all}
\quad
\nu \vdash d.
\end{equation}
In this section, we will specify sufficient conditions to lift the isomorphism \eqref{eq:Biso} to the level of modules over the quantum wreath product $B\wr \cH(d)$.
}
%-----------------------------------------------------------------------
{
\begin{thm}\label{thm:dualS}
Suppose that there is an automorphism $\#$ on $B$, an anti-automorphism $\ast$ on $B$,  and an order-reversing involution $t$ on a poset $I_B$ such that \eqref{eq:Biso} holds.
Assume that \eqref{eq:Hiso} holds for all $\tcH^{\chi_\lambda}_d$.
Then, as $B\wr \cH(d)$-modules:
\[
(S^\bblambda)^\# \cong (S^{\bblambda^t} )^\ast
\quad
\textup{for any}  
\quad
\bblambda \in \Pi_{I_B}(d)
\]
if the following hold:
\begin{enumerate}[\textup{(}i\textup{)}]
    \item The maps $\#$ in  \eqref{def:autA} and $\ast$ in \eqref{def:antiautA} are well-defined for $B\wr \cH(d)$. That is, \eqref{eq:weakFrob}, \eqref{eq:sharpB}--\eqref{eq:sharpBB}, \eqref{eq:starB}--\eqref{eq:starW} all hold.
    \item For all $\lambda \in I_B$, 
    {there is $\alpha_\lambda \in K^\times$ such that 
    \[
    \chi_{\lambda^t}(S) = \alpha_\lambda\chi_\lambda(S),
    \quad
    \chi_{\lambda^t}(R) = \alpha_\lambda^2\chi_{\lambda}(R).
    \]}
    \item Either $S = 0$ or $\sigma^{S^\lambda_B} = 0$ for all $\lambda \in I_B$.
\end{enumerate}
%Here, $\chi_\lambda$ records the action of $S$ and $R$ on $S^\lambda_B\otimes S^\lambda_B$ and $\tcH^{\chi_\lambda}_d$ is the twisted Hecke algebra with quadratic relations of the form $h_i^2 = \chi_\lambda(S)h_i + \chi_\lambda(R)$.
{Note that ({\it ii}) implies that $h_i \mapsto \alpha_\lambda h_i$ defines an isomorphism $\tcH^{\chi_{\lambda^t}}_d \to \tcH^{\chi_{\lambda}}_d$.
}
%$\sigma^{M}=\textup{flip}$, $\rho^M = 0$ for all Specht module $M = S^\lambda_m$,}
\end{thm}
%-----------------------------------------------------------------------
\begin{proof}
It suffices to verify the case when $S^\bblambda = S^\lambda_B \wr S^\nu_{\chi_\lambda}$ for some $\lambda \in I_B, \nu \vdash d$.
Then, 
\begin{equation}\label{eq:dualS}
(S^\lambda_B \wr S^\nu_{\chi_\lambda})^\#
\cong
(S^\lambda_B)^\# \wr (S^\nu_{\chi_\lambda})^\#
\cong
(S^{\lambda^t}_B)^* \wr (S^{\nu^t}_{\chi_{\lambda^t}})^*
\cong
(S^{\lambda^t}_B \wr S^{\nu^t}_{\chi_{\lambda^t}})^*
=
(S^{\bblambda^t})^*.
\end{equation}
Here, the first isomorphism in \eqref{eq:dualS} follows from Theorem \ref{thm:MNsharp}, which is available since ({\it iii}) holds.
The second isomorphism in \eqref{eq:dualS} follows from combining the assumptions \eqref{eq:Biso}, ({\it ii}) and Lemma \ref{lem:tHAdual}.
The third isomorphism follows from Theorem \ref{thm:MNstar} since ({\it i}) implies the condition $\sigma^2 = \id$.
Finally, $S^{\lambda^t}_B \wr S^{\nu^t}_{\chi_{\lambda^t}} = S^{\bblambda^t}$ since $(e_\lambda^\nu)^t = e_{\lambda^t}^{\nu^t}$ from \eqref{def:t}.

The general case follows since the functors $(-)^{\#}$ and $(-)^*$ commute with the parabolic induction functor $\ind_\mu^d$.
\end{proof}
}
%=============================================
\section{The Hu Algebras}\label{sec:Hu}
%=============================================
\subsection{Hu Algebras as Quantum Wreath Products}
%=============================================
From now on, assume that $K$ is a field of characteristic not equal to 2 containing $q = v^2$, and that \eqref{eq:separate} holds.
Let $B = \cH_q(\Sigma_m)$ be the Hecke algebra for $\Sigma_m$ with generators $I_1, \dots, I_{m-1}$ subject to the braid relations and the following normalized quadratic relations:
\begin{equation}
(I_i + v\inv)(I_i - v) = 0,
\quad
\textup{or,}
\quad
I_i^2 = (v-v\inv) I_i +1.
\end{equation}
That is, $I_i = T_i/v$ where $T_i$ is the usual generator satisfying $(T_i+1)(T_i - v^2) = 0$.
In particular, the bar-involution on $B$ sends each $I_{w}$ to $\overline{I}_w := I_{w\inv}\inv$.

For positive integers $a \leq b < 2m$, define the following elements in $\cH_q(\Sigma_{2m})$:
\begin{equation}
I_{a\to b}^{\epsilon_i} := 
\begin{cases}
I_a I_{a+1} \dots I_b &\tif \epsilon_i = +;
\\
I\inv_a I\inv_{a+1}\dots I\inv_b &\tif \epsilon_i = -,
\end{cases}
\quad
C_\epsilon = \prod_{i;\epsilon_i=-}  (I_{2m-i} \dots  I_{m+1}) (I_{m+1}  \dots I_{2m-i}).
\end{equation}
Let $\Hh =\Hh(m) \in \cH_q(\Sigma_{2m})$ be the following element 
\begin{equation}\label{eq:preciseH1}
\Hh := v^{m(2m-1)}{\textstyle \sum_{\epsilon = (\epsilon_1, \dots, \epsilon_m) \in \{+,-\}^m}  
I_{m\to 2m-1}^{\epsilon_1}
\dots
I_{1\to m}^{\epsilon_m} C_\epsilon},
\end{equation}
which will later be identified with the generator $H_1$ in a quantum wreath product.
{
In particular, for the extremal summand corresponding to $\epsilon = (-,\dots,-)$, the coefficient $C_\epsilon = \gamma^2$ is the full twist on the second half $\{m+1, \dots, 2m\}$.
To be precise, the full twist is the square of the following half twist:
\begin{equation}
\gamma := \iota^{-1}(1\otimes I_{w_\circ})
=
I_{m+1} (I_{m+2} I_{m+1}) \dots (I_{2m-1} \dots I_{m+1}) \in \cH_q(\Sigma_m\times \Sigma_m),
\end{equation}
where $w_\circ \in \Sigma_m$ is the longest element, and $\iota$ is the identification map $\cH_q(\Sigma_m \times \Sigma_m) \to B\otimes B$.
}

{
The element $\Hh \in \cH_q(\Sigma_{2m})$ has two crucial properties: 
First, it squares to an element $z = z(m)$ that is central in $\cH_q(\Sigma_m\times \Sigma_m)$. 
Second, it behaves extremely well under the bar involution (and hence the automorphism $\#$), 
leading to the bar-invariant element 
\begin{equation}
    b_1(m):= v^{-m(2m-1)} \Hh \overline{\gamma},
\end{equation}
which is useful for constructing a bar-invariant basis \cite{LNX24}.
In turn, such a property motivated the definition of the automorphism $\#$ in Section \ref{sec:auto}.
}
%----------------------------------------------------------
\begin{prop}[{\cite[Proposition~6.2.2, Theorem 6.3.3, Proposition 6.4.2]{LNX24}}]\label{prop:H1}
Suppose that  $f_{2m}(q):= \prod_{i=0}^{2m-1}(1+q^i)$ is nonzero.
Then, 
\begin{enua}
    \item $z:= \Hh^2 \in Z(\cH_q(\Sigma_m \times \Sigma_m))$ and hence $Z = Z(m) := \iota(z) \in Z(B \otimes B)$.
    Moreover, $z$ and $Z$ are invertible.
    \item The element $\Hh$ intertwines $\cH_q(\Sigma_m\times \Sigma_m)$, i.e., 
    \[
    \Hh I_{(x,y)} = I_{(y,x)} \Hh
    \quad
    \textup{for all}
    \quad
    (x,y) \in \Sigma_m\times \Sigma_m.
    \]
    \item {
    The element $b_1$ is bar-invariant.
    Equivalently,
    \[
    \overline{\Hh} = v^{-2m(2m-1)} \Hh \gamma^{-2}.
    \]}
\end{enua}
\end{prop}
%----------------------------------------------------------
{
One can use Proposition \ref{prop:H1} to compute the relation between $Z$ and $Z^\#$ as follows:
}
%----------------------------------------------------------
\begin{prop}\label{prop:Zsharp}
{
Let $\Gamma := I_{w_\circ} \otimes I_{w_\circ} \in B\otimes B$. Then
    $
    Z^\# = Z \Gamma^{-2}.
    $    
}
\end{prop}
%----------------------------------------------------------
\begin{proof}
It follows from $\overline{I}_i = I_i^{-1}$ and $I_i^\# = -\overline{I}_i$ 
that $\overline{I}_i^\# = -I_i$.
One can then combine \eqref{eq:preciseH1} and Proposition \ref{prop:H1}(c)  to show that $\Hh^\#$ differs from $\Hh$ by a full twist $\gamma^2$ (up to a sign), i.e.,
$\Hh^\# = (-1)^m \Hh \gamma^{-2}$.
Therefore, by Proposition \ref{prop:H1}(b),
\[
Z^\# 
= \iota(\Hh^2)^\#
= \iota((\Hh^\#)^2)
= \iota(\Hh \gamma^{-2} \Hh \gamma^{-2})
= \iota(\Hh^2) (I_{w_\circ}^{-2} \otimes I_{w_\circ}^{-2})
= Z \Gamma^{-2}.
\]
\end{proof}
%----------------------------------------------------------
For $\lambda\vdash m$, denote by $S_m^\lambda$ the Specht module over $\cH_q(\Sigma_m)$, and by $D_m^\lambda$ the simple quotient of $S_m^\lambda$ by its radical.
Let $\lambda', \lambda'' \vdash m$. It follows from \cite[Lemmas 3.2 and 3.4]{Hu02} that $Z$ acts on $S_m^{\lambda'} \otimes S_m^{\lambda''}$ by multiplication by a scalar which we call $f_{\lambda', \lambda''} \in \QQ(v)$. If, additionally, $\lambda' = \lambda = \lambda''$, then $\sqrt{f_{\lambda}} := \sqrt{f_{\lambda,\lambda}}$ is well-defined (up to a sign) and lies in $\ZZ[v, v\inv]$ thanks to \cite[Theorem 3.5]{Hu02}.
%----------------------------------------------------------
\begin{exa} 
Let $d=2$.
For $m=1$, the central element $Z$ is given by
\[
Z = (v(I_1 + I_1\inv))^2 = v^2 (v+v\inv)^2 = (v^2+1)^2 ,
\]
which is the square of a quantum integer two  $[2]_{v^2}$ where $[2]_u := (u+1)$. 
{In this case,
there is only one $\sqrt{f_\lambda}^2$ given by $\sqrt{f_{(1)}} = \pm (q+1)$.
On the other hand, $\gamma = 1$, $\Gamma = 1\otimes 1$, and hence indeed $Z^\# = (q+1)^2 = Z \Gamma^{-2}$.
}
For $m=2$, one computes that
\begin{equation}
\begin{split} 
 Z =v^{12}&\big(
 (v-v\inv)^2 (\phi(v) - 2(v^2 + v^{-2})) (I_{1}\otimes I_1)
 \\
&+(v-v\inv) \phi(v) (I_1\otimes 1 +1\otimes I_1)
+2\phi(v)
\big),
\end{split}
\end{equation}
where $\phi(v) := v^4+{4}v^2-2+{4}v^{-2}+v^{-4}$. 
{It is worth mentioning that for the one-row and the one-column partitions, their corresponding $f_\lambda$ are (normalized) products of squares of various quantum integer twos:}
\[
f_{(m)}(v) = v^{4m(m-1)}
[2]_{v^2}^2 [2]_{v^4}^2 \dots [2]_{v^{2m}}^2
= v^{2m(m-1)} f_{(1^m)}(v).
\]
{
Then, since $\gamma = \iota^{-1}(1\otimes I_1) = I_3$, one has $\Gamma =I_1\otimes I_1$. 
It then follows from \eqref{eq:preciseH1} that
\[
\begin{split}
\Hh^\# 
&= v^6\big(
(I_2I_3I_1I_2 +I_2I_3\overline{I_1I_2})^\#
+ (\overline{I_2I_3}I_1I_2 +\overline{I_2I_3I_1I_2})^\#(I_3^2)^\#
\big)
\\
&= 
v^6
\big(
(\overline{I_2I_3I_1I_2}+\overline{I_2I_3}I_1I_2)
+ (I_2I_3\overline{I_1I_2} +I_2I_3I_1I_2)(\overline{I}_3^2)
\big)
\\
&= \Hh I_3^{-2}.
\end{split}
\]
Therefore, 
$Z^\# 
= \iota((\Hh^\#)^2)
= \iota(\Hh I_3^{-2} \Hh I_3^{-2})
= \iota(\Hh^2 I_1^{-2} I_3^{-2})
= Z (I_1^{-2}\otimes I_1^{-2})
= Z \Gamma^{-2}$.
}
\end{exa}
%----------------------------------------------------------
Following \cite{Hu02}, let $\cA(m)$ be the subalgebra of $\cH_q(\Sigma_{2m})$ generated by $\cH_q(\Sigma_m\times\Sigma_m)$ and $\Hh$.
Combining $\iota:\cH_q(\Sigma_m\times\Sigma_m) \to B\otimes B$ and the assignment $\Hh \mapsto H_1$, the Hu algebra $\cA(m)$ is isomorphic to the quantum wreath product $B\wr \cH(2)$ with the following choices:
\begin{equation}\label{par:Hu}
S=0,
\quad
R = Z(m),
\quad 
\sigma: a\otimes b \mapsto b\otimes a,
\quad
\rho = 0.
\end{equation}
In other words, $\cA(m)$ is generated by $B^{\otimes 2}$ and $H_1$ subject to the following relations:
\begin{equation}
H_1^2 = Z(m),
\quad H_1(a \otimes b) = (b \otimes a) H_1,
\qquad
(a,b \in B).
\end{equation}
{From now on, by a generalized Hu algebra we mean a quantum wreath product $B \wr \cH(d)$ for $d\geq 2$ with respect to the choice of parameters \eqref{par:Hu}.}
%=============================================
\subsection{Representation Theory of Hu Algebras}\label{sec:RTHu}
%=============================================
The complete set of Specht and simple modules over the Hu algebra is best described using the set $\Omega = \Pi_{\Pi_m}(2)$ of $\Pi_m$-indexed multipartitions of $d=2$.
Let $\chi_\lambda$ be given by  
$\chi_\lambda(S) = 0$, $\chi_\lambda(R) = f_{\lambda,\lambda}$.
For each $\bblambda \in \Omega$, the assignment $\lambda \mapsto \bblambda(\lambda)$ indicates that the wreath module $S^\lambda_m \wr S_{\chi_\lambda}^{\bblambda(\lambda)}$ will be used to construct $S^\bblambda$.
When $d=1$, it is understood that $\tcH^{\chi}_1 \equiv K$ with a single Specht module $S_{\chi}^{(1)} \equiv K$.
Hence, the wreath module $S^\lambda_m \wr S_{\chi_\lambda}^{(1)} \equiv S^\lambda_m$ is the Specht module over $B$.

The Specht modules $S(\lambda,\lambda)_\pm$ and $S(\lambda',\lambda'')$ in \cite{Hu02} are identified as in the table below:
\[
\begin{array}{|c||c|c|c|}
\hline
\textup{Specht module in \cite{Hu02}}&S(\lambda, \lambda)_+ &S(\lambda, \lambda)_- &S(\lambda', \lambda'')
\\
\hline
\textup{Multipartition }\bblambda&  e_\lambda^{(2)} & e_\lambda^{(1,1)} & e_{\lambda'}^{(1)} + e_{\lambda''}^{(1)}
\\
\hline
\textup{Our Specht module }S^\bblambda & S_m^\lambda \wr S^{(2)}_{\chi_\lambda} & S_m^\lambda \wr S^{(1,1)}_{\chi_\lambda} 
& \ind_{(1,1)}^2\big((S_m^{\lambda'} \wr S_{\chi_{\lambda'}}^{(1)})\otimes (S_m^{\lambda''}\wr S_{\chi_{\lambda''}}^{(1)}) \big)
\\
\hline
\end{array}
\]
In other words, this yields the following characterization of Specht modules for the Hu algebras.
%--------------------------------------
\begin{prop}[{\cite[Section 4]{Hu02}}]\label{prop:gHuIrrep}
The complete set of Specht modules, up to isomorphism, of the Hu algebra $\cA(m)$ is given by $\{S^\bblambda~|~ \bblambda\in\Omega\}$.
\end{prop}
%--------------------------------------

The following identifications in the examples below will be useful in Example \ref{ex:POD} when various partial orders are being compared.

%------------------------------------------------------------
\begin{exa}\label{ex:Huindex}
For the Hu algebra $\cA(m)$ or the group algebra of $\Sigma_m \wr \Sigma_2$, 
%the base algebras are the Hecke algebra $B= \cH_q(\Sigma_m)$ and the group algebra $K\Sigma_m$, respectively.
the Specht modules over $B$ are in bijection with the set $I_B =\Pi_m$ of partitions of $m$. 
In the special case $m=2$,  the  poset structure $I_B =\{(1,1) < (2)\}$  induces a partial order on $\Omega$ given by
\begin{equation}\label{eq:HuPoset}
\left(\begin{array}{l}  (2)\mapsto\varnothing \\ (1,1) \mapsto (1,1)\end{array}\right)
\unlhd
\left(\begin{array}{l} (2) \mapsto \varnothing\\ (1,1)\mapsto (2)\end{array}\right)
\unlhd
\left(\begin{array}{l} (2) \mapsto (1)\\ (1,1)\mapsto (1)\end{array}\right)
\unlhd
\left(\begin{array}{l} (2) \mapsto (1,1)\\ (1,1)\mapsto \varnothing\end{array}\right)
\unlhd
\left(\begin{array}{l} (2) \mapsto (2) \\ (1,1)\mapsto \varnothing\end{array}\right).
\end{equation}
Note that the decorated orbits in \cite{Hu02} are identified as
$[\lambda, \lambda]_+ \equiv e_\lambda^{(2)}$, $[\lambda, \lambda]_- \equiv e_\lambda^{(1,1)}$, and $[\lambda, \mu] \equiv e_\lambda^{(1)} + e_\mu^{(1)}$, and hence 
 \eqref{eq:HuPoset} becomes 
\[
 [(1,1), (1,1)]_- \unlhd [(1,1), (1,1)]_+ \unlhd [(1,1), (2)] \unlhd [(2), (2)]_- \unlhd [(2), (2)]_+.
\]
\end{exa}
\begin{exa}
Consider the Hu algebra $B \wr \cH(2)$ with  base algebra $B = \cH_q(\Sigma_m)$.
The Specht modules of the form $S^\bblambda = S(\lambda^{(1)}, \lambda^{(2)})$ for $\lambda^{(1)} \neq \lambda^{(2)} \in \Pi_m$ correspond to the multipartition $\bblambda$ supported on $\{\lambda^{(1)}, \lambda^{(2)}\}$ with the same value $(1)$.
Thus, $\bblambda^t$ is supported on $\{(\lambda^{(1)})^t, (\lambda^{(2)})^t\}$ with the same value $(1)$.
That is, the  induced transpose map $S^\bblambda \mapsto S^{\bblambda^t}$ takes $S(\lambda^{(1)}, \lambda^{(2)})$ to $S((\lambda^{(1)})^t, (\lambda^{(2)})^t)$ for $\lambda^{(1)} \neq \lambda^{(2)} \in \Pi_m$.

The other Specht modules $S^\bblambda = S(\lambda, \lambda)_\pm$ correspond to the multipartitions 
\begin{equation}
\bblambda: \nu \mapsto
\begin{cases}
(2) &\tif \nu = \lambda;
\\
\varnothing &\textup{otherwise},
\end{cases}
\quad\textup{and}\quad
\bblambda: \nu \mapsto
\begin{cases}
(1,1) &\tif \nu = \lambda;
\\
\varnothing &\textup{otherwise},
\end{cases}
\end{equation}
respectively. 
Hence, the induced transpose map swaps the two Specht modules $S(\lambda,\lambda)_\pm$ and $S(\lambda^t,\lambda^t)_\mp$ for any $\lambda \in \Pi_m$. 
\end{exa}

%=============================================
\subsection{Wreath Modules over Hu Algebras}\label{sec:WMHu}
%=============================================
Let $B \wr \cH(d)$ be the generalized Hu algebra $(d\geq 2)$.
We obtain the following wreath modules as a consequence of Theorem~\ref{thmWM}.

%--------------------------------------
\begin{prop}\label{prop:MwrNHu}
Suppose that $M, M_1, \dots, M_r \in \cH_q(\Sigma_m)$-mod such that the element $Z \in \cH_q(\Sigma_m) \otimes \cH_q(\Sigma_m)$ acts on $M\otimes M, M_i\otimes M_i$ by  multiplication by $\chi(Z)$ and $\chi_i(Z)$, respectively.
\begin{enua}
\item
If $N$ is a module over a  twisted Hecke algebra $\tcH^{\chi}_d$,
then $M\wr N \in B \wr\cH(d)$-mod, on which the action is given by, for $b\in \cH_q(\Sigma_m)^{\otimes d}, m \in M^{\otimes d}$, $y\in N$:
\[
b \cdot (m \otimes y) = (b \cdot m) \otimes y,
\quad
H_i \cdot (m \otimes y) =  \sigma^M_i(m) \otimes (h_i \cdot y),
\]
where $\sigma^M(m' \otimes m'') = m''\otimes m'$ is the flip map on $M\otimes M$.
\item
Let $\mu = (\mu_1, \dots, \mu_r) \vDash d$.
For $i = 1, \dots, r$, if each $N_i$ is a module over a  twisted Hecke algebra $\tcH^{\chi_i}_{\mu_i}$, then
\[
(M_1 \wr N_1) \otimes (M_2 \wr N_2) \otimes  \dots \otimes (M_r \wr N_r) 
\equiv
M_1^{\otimes \mu_1} \otimes \dots \otimes M_r^{\otimes \mu_r} \otimes N_1 \otimes \dots \otimes N_r
\] 
is a module over the parabolic subalgebra $\cH_q(\Sigma_m) \wr \cH(\mu)$ of the generalized Hu algebra.
\end{enua}
\end{prop}
%--------------------------------------
%We need to be careful when $I_B$ is not a totally ordered set.
We may fix an enumeration $I_B = \{\lambda^{(1)}, \lambda^{(2)}, \lambda^{(3)}, \dots\}$ that is compatible with the dominance order on $I_B = \Pi_m$, i.e., if $\lambda^{(i)} \lhd \lambda^{(j)}$ then $i<j$.
Thus,  each $\bblambda \in \Omega$ has support of the form 
$\{\lambda^{(k_i)} ~|~  k_1 < k_2 < \dots < k_r \}$ for some $r$.
Write 
$|\bblambda| := (\#\bblambda(\lambda^{(k_1)}), \dots , \#\bblambda(\lambda^{(k_r)})) \vDash d$.
%,
%and $|\bblambda|^+ = w |\bblambda|$ the partition of $d$ obtained from $|\bblambda|$ by the shortest place permutation $w\in \Sigma_r$.
For example, if $\bblambda = e_{\lambda^{(1)}}^{(2)} + e_{\lambda^{(3)}}^{(1)} +e_{\lambda^{(4)}}^{(1,1)}$, then $|\bblambda| = (2,1,2)$.
%, $w = s_2$ and $|\bblambda|^+ = (2,2,1)$.
For each $i$, we set $M_i := S^{\lambda^{(k_i)}}_m$ to be the corresponding Specht module for $\cH_q(\Sigma_m)$,
and let $\chi_i$ be the twist determined by $\chi_i(S) = 0$ and $\chi_i(R) = f_{\lambda^{(k_i)}}$.

Next, let
$\nu^{(i)} := \bblambda(\lambda^{(k_i)})$. 
We pick  $N_i = S^{\nu^{(i)}}_{\chi_i}$ to be the Specht module over $\tcH^{\chi_i}_{\mu_i}$.
Consequently, Proposition \ref{prop:MwrNHu} implies that the following vector space, which we call the Specht module, is indeed a module over the Hu algebra $B \wr \cH(d)$ for $d\geq 2$:
\begin{equation}\label{def:SpechtgHu}
S^\bblambda := \ind_{|\bblambda|}^d  \big( 
(S_{m}^{\lambda^{(k_1)}} \wr S_{\chi_1}^{\nu^{(1)}}) 
\otimes
\dots
\otimes
(S_{m}^{\lambda^{(k_r)}} \wr S_{\chi_r}^{\nu^{(r)}})
\big).
\end{equation}
Note that $S^\bblambda$ does not depend on the choice of enumeration of $I_B$.

%=============================================
\subsection{The Dual Specht Modules}
%=============================================
{
Let $\cH_q(\Sigma_m) \wr \cH(d)$ be the generalized Hu algebra $(d\geq 2)$.
We begin with a lemma which verifies the condition (ii) in Theorem~\ref{thm:dualS}, which will produce isomorphisms  $(S^\bblambda)^\# \cong (S^{\bblambda^t})^\ast$ for all $\bblambda \in \Omega$.
%---------------------------
\begin{lemma}\label{lem:ftrans}
For $\lambda\in\Pi_m$,
let $n(\lambda):= \sum_{i\geq 1} (i-1)\lambda_i$,
and let $z_\lambda\in K^{\times}$ be the scalar by which $I_{w_\circ(m)}^{-2}$ acts on $S^{\lambda}_m$. Then:
\begin{enua}
\item $I_{w_\circ}^{2}=L_1\cdots L_m$, where $L_1 :=1 $, $L_{i+1} := I_iL_iI_i$ are the Jucys--Murphy elements of $B$. 
\item $z_\lambda = q^{n(\lambda) - n(\lambda^{t})}$. 
In particular, 
$z_\lambda z_{\lambda^{t}}=1$ 
and 
$z_\lambda=v^{2(n(\lambda)-n(\lambda^{t}))}$ is a square in $\ZZ[v,v^{-1}]^{\times}$.
\item 
For all $\lambda,\mu\in\Pi_m$, 
$f_{\lambda^{t},\mu^{t}} = z_\lambda z_\mu f_{\lambda,\mu}$.  
In particular,
      $f_{\lambda^{t}} = z_\lambda^{2}f_{\lambda}$.
\end{enua}
\end{lemma}
\begin{proof}
{
Part (a) follows by induction on $m$.
For part (b), it follows from part (a) that $I_{w_\circ(m)}^{-2}$ acts by the scalar
$q^{-\sum_{b\in\lambda}\textup{cont}(b)}$ where $\textup{cont}(i,j) := j-i$ is the content of the box $b = (i,j)$ in the Young diagram $\lambda$.
Part (b) then follows from combining \cite[Proposition 3.35]{Ma99} and the content formula $\sum_{b\in\lambda} \textup{cont}(b) = n(\lambda^t) - n(\lambda)$.

For Part (c), note that the central element $Z = Z^\ast$ \cite[Lemma 1.10]{Hu02} acts on
$(S^{\lambda^{t}}_m)^{*}\otimes(S^{\mu^{t}}_m)^{*}$
by $f_{\lambda^{t},\mu^{t}}$. 
On the other hand, $Z$ acts on
$(S^{\lambda}_m\otimes S^{\mu}_m)^{\#}$ by the scalar of $Z^{\#} = Z \Gamma^{-2}$ on
$S^{\lambda}_m\otimes S^{\mu}_m$, i.e., $z_\lambda z_\mu f_{\lambda,\mu}$.
The lemma is proved.
}
\end{proof}
%---------------------------
{
Thanks to Lemma \ref{lem:ftrans}(c),
we may fix the choice of $\{\sqrt{f_\lambda}\}_\lambda$ such that
$\sqrt{f_{\lambda^{t}}} = z_\lambda \sqrt{f_{\lambda}}$ for all $\lambda$.
}
One can then apply Theorem \ref{thm:dualS} to obtain the isomorphism $(S^\bblambda)^\# \cong (S^{\bblambda^t})^\ast$ as follows.
%---------------------------
\begin{prop} \label{prop:genSpechtiso}
Suppose that $\cH_q(\Sigma_m)\wr \cH(d)$ is the generalized Hu algebra $(d\geq 2)$. Then,
\begin{enua}
\item     
Let $\theta := 1\otimes I_{w_\circ}^{-2}$.
Then, $\theta \in Z(B\otimes B)^\times$, and the map $\# = \#(\theta)$ (see \eqref{def:autA}) defines an automorphism for $\cH_q(\Sigma_m) \wr \cH(d)$. In particular,
\[
\begin{split}
H_1^\# &=  - H_1 (1\otimes I_{w_\circ}^{-2}),
\\
(T_x \otimes T_y)^\#
&=(-q)^{\ell(x)+\ell(y)}(T_{x^{-1}}^{-1} \otimes T_{y^{-1}}^{-1})
\quad
\textup{for all}
\quad
x,y \in \Sigma_m.
\end{split}
\]
\item 
The map $*$ in \eqref{def:antiautA} defines an anti-automorphism for $\cH_q(\Sigma_m) \wr \cH(d)$.
In particular,
\[
((T_x \otimes T_y) H_1)^* 
= (T_{y^{-1}} \otimes T_{x^{-1}}) H_1, 
\quad 
\textup{for all } 
x, y \in \Sigma_m.
\]
\item 
For any  $\bblambda \in \Omega$, $(S^\bblambda)^\# \cong (S^{\bblambda^t})^\ast$. 
\end{enua}
\end{prop}
\begin{proof}
For part (a), the first assertion follows from combining Lemma \ref{lem:ftrans} and the fact that $w_\circ = w_\circ^{-1}$.
By Remark \ref{rmk:sharpHu},
the second assertion reduces to the 
verification of $Z^\# = Z \theta \sigma(\theta)$ for some $\theta \in Z(B\otimes B)^\times$, which follows from Proposition \ref{prop:Zsharp}.
Indeed, $\theta \sigma(\theta) = (I_{w_\circ}^{-2} \otimes I_{w_\circ}^{-2}) = \Gamma^{-2}$.

For Part (b), it follows from applying Theorem \ref{thm:*}.
The only nontrivial equality to check is $Z^* = Z$, which follows from \cite[Lemma 1.10]{Hu02}.

Part (c) follows by verifying the assumptions ({\it i})--({\it iii}) in Theorem \ref{thm:dualS}. 
Set $S^\lambda_B = S^\lambda_m$ for $\lambda \in I_B = \Pi_m$. 
Since the twisted Hecke algebra $\tcH^{\chi_\lambda}_d$ has quadratic relations of the form $(h_i + \sqrt{f_\lambda})(h_i - \sqrt{f_\lambda}) = 0$,
Lemma~\ref{lem:tHAdual} applies.
First, we have checked ({\it i}) in Parts (a) and (b) of this proposition.
For ({\it ii}), Lemma \ref{lem:ftrans} gives
$\chi_{\lambda^{t}}(R)=f_{\lambda^{t}}=z_\lambda^{2}f_\lambda=z_\lambda^{2}\chi_\lambda(R)$,
while $\chi_{\lambda^{t}}(S)=0=\chi_\lambda(S)$. 
Hence, ({\it ii}) holds with
$\alpha_\lambda= z_\lambda\in K^{\times}$. 
Finally, ({\it iii}) holds since $S=0$ and $\sigma$ is the flip map.
\end{proof}
}
%=============================================
\subsection{The new trivial and sign modules}\label{sec:trivsgn}
%=============================================
The Hecke algebra $\cH_q(\Sigma_m)$ carries no natural comultiplication, so one cannot
expect the tensor product of two $\cH_q(\Sigma_m)$-modules to be an $\cH_q(\Sigma_m)$-module
again. What the automorphism $\#$ provides is a substitute for tensoring with the sign
module: the analog of Theorem \ref{thm:dualS} over $\cH_q(\Sigma_m)$ in the literature
usually takes the form
\begin{equation}\label{eq:sgntwist}
S^\lambda_m \otimes \sgn \cong (S^{\lambda^t}_m)^\ast \in \cH_q(\Sigma_m)\textup{-mod},
\quad
\textup{where}
\quad
\sgn := S_m^{(1^m)},
\end{equation}
in which the left-hand side is to be read as $(S^\lambda_m)^\#$.

For the Hu algebras, it turns out that $(S^\bblambda)^\#$ cannot be regarded as a tensor
product with a single sign module.
Suppose that $S^\bblambda = \ind_\mu^d ( \bigotimes_i  S^{\lambda^{(i)}}_m  \wr S_{\chi_i}^{\nu^{(i)}} )$
as in \eqref{def:SpechtgHu}.
Since $h_i$ acts on $(S^{\nu^{(i)}}_{\chi_i})^\#$ by $-\alpha_{\lambda^{(i)}}h_i$, it follows
from Theorem \ref{thm:MNsharp} that
\begin{equation}\label{eq:sharpSpechtHu}
(S^\bblambda)^\# \cong \ind_\mu^d  \Big(\bigotimes\nolimits_i
(S^{\lambda^{(i)}}_m \otimes \sgn) \wr (S_{\chi_i}^{\nu^{(i)}} \otimes K_{-\sqrt{f_{(\lambda^{(i)})^t}}} )
\Big).
\end{equation}
%where the $i$-th twisted Hecke factor is a module over $\tcH^{\chi_{(\lambda^{(i)})^t}}_{\mu_i}$
%rather than over $\tcH^{\chi_i}_{\mu_i}$. 
%Note that this is forced: the ground module
%$S^{\lambda^{(i)}}_m \otimes \sgn \cong (S^{(\lambda^{(i)})^t}_m)^\ast$ has $Z$ acting by
%$f_{(\lambda^{(i)})^t}$, and the two parameters must match for the tensor factor to be a
%wreath module.

It is tempting to believe that the following $K$-module (which depends on $\bblambda$) could
be a $B\wr \cH(d)$-module:
\begin{equation}\label{eq:sgnlambda}
\sgn_{\bblambda}:= \ind_\mu^d  \Big( \bigotimes\nolimits_i
(\sgn \wr K_{-\sqrt{f_{(\lambda^{(i)})^t}}})
\Big).
\end{equation}
This fails for two independent reasons. First, by the matching just described, each tensor
factor $\sgn \wr K_{-\sqrt{f_{(\lambda^{(i)})^t}}}$ fails to be a wreath module unless
$f_{(\lambda^{(i)})^t} = f_{(1^m)}$, that is, unless $\lambda^{(i)}$ is a single row.
Second, the one-dimensional factors in \eqref{eq:sgnlambda} live over the pairwise distinct
twisted Hecke algebras $\tcH^{\chi_{(\lambda^{(i)})^t}}_{\mu_i}$, whose parameters vary with
$i$.
Thus, no single sign module is available even formally.
Consequently, we will avoid the notion of tensoring with the sign module, and consider only the
following special cases of (dual) Specht modules.

\begin{Def}\label{def:trivsgnHu}
Let $\id := e_{(m)}^{(d)} \in \Omega$, and hence $\id^t = e_{(1^m)}^{(1^d)}$.
Denote the trivial and sign modules for the generalized Hu algebra $\cH_q(\Sigma_m)\wr \cH(d)$
by, respectively:
\[
S^\id = S_m^{(m)} \wr K_{\sqrt{f_{(m)}}},
\quad
S^{\id^t} = S_m^{(1^m)} \wr K_{-\sqrt{f_{(1^m)}}},
\]
where the parameters of the two twisted Hecke algebras are related by
$\sqrt{f_{(1^m)}} = \alpha_{(m)}\sqrt{f_{(m)}}$ with
$\alpha_{(m)} = v^{-m(m-1)}$.
\end{Def}
These modules behave nicely under the maps $\#$ and $\ast$ as follows.

\begin{prop}\label{prop:trivstar}
    As modules over the generalized Hu algebra,
    \[
    (S^\id)^\ast \cong S^\id,
    \quad
    (S^{\id^t})^\ast \cong S^{\id^t},
    \quad
    (S^\id)^\# \cong S^{\id^t},
    \quad
    (S^{\id^t})^\# \cong S^\id.
    \]
\end{prop}
\begin{proof}
All modules involved are finite-dimensional, so Theorem \ref{thm:MNstar} applies and gives
\[
(S^\id)^\ast \cong (S_m^{(m)})^\ast \wr (K_{\sqrt{f_{(m)}}})^\ast.
\]
Each of these one-dimensional modules is self-dual, since $I_i^\ast = I_i$ and
$h_i^\ast = h_i$,
\[
\begin{split}
(I_i \circ \varphi)(x) = \varphi(I_i \cdot x) = v\,\varphi(x)
\quad&\textup{for all}\quad
\varphi \in (S_m^{(m)})^\ast,
    \\
(h_i \circ \psi)(y) = \psi(h_i \cdot y) = \sqrt{f_{(m)}}\,\psi(y)
\quad&\textup{for all}\quad
\psi \in K_{\sqrt{f_{(m)}}}^\ast.
\end{split}
\]
Therefore $(S^\id)^\ast \cong S_m^{(m)} \wr K_{\sqrt{f_{(m)}}} = S^\id$.
The proof of $(S^{\id^t})^\ast \cong S^{\id^t}$ is identical, using that $I_i$ acts on
$S_m^{(1^m)}$ by $-v\inv$ and $h_i$ acts on $K_{-\sqrt{f_{(1^m)}}}$ by $-\sqrt{f_{(1^m)}}$.

For the remaining two isomorphisms, Theorem \ref{thm:MNsharp} gives
\[
(S^\id)^\# \cong (S_m^{(m)})^\# \wr (K_{\sqrt{f_{(m)}}})^\#
\cong S_m^{(1^m)} \wr K_{-\sqrt{f_{(1^m)}}} = S^{\id^t},
\]
where the second isomorphism uses \eqref{eq:dualSpechtHecke} together with the self-duality of
$S_m^{(1^m)}$ for the first factor, and, for the second factor, the fact that $h_i$ acts on
$(K_{\sqrt{f_{(m)}}})^\#$ by $-\alpha_{(m)}\sqrt{f_{(m)}} = -\sqrt{f_{(1^m)}}$.
The proof of $(S^{\id^t})^\# \cong S^\id$ is similar.
Alternatively, both follow from Theorem \ref{thm:dualS} and the two self-duality statements
above, since $(\id)^t = \id^t$.
\end{proof}
%========================
\subsection{A Transfer Map}
%========================
In this section, we start off by defining a transfer theory for the Hu algebra.
%========================
\begin{lem}\label{lem:transfer}
Let $B= \cH_q(\Sigma_m)$, and $\bH = B\wr \cH(d)$ be the generalized Hu algebra with $d\geq 2$ and $d! \in K^\times$. Let $M$ and $N$ be two $\bH$-modules. Then, ${\rm Tr} : \Hom_{\cB}(M, N) \to \Hom_{\bH}(M, N)$ given by
\[
{\rm Tr}(f)(x):=\frac{1}{d!}\sum_{w\in\Sigma_d} H_w^{-1}f(H_w x)
\]
is well-defined, and restricts to the identity on $\Hom_{\bH}(M, N)$. Moreover, for $\bH$-linear maps $g : N \to N'$ and $h : M' \to M$, ${\rm Tr}(g \circ f \circ h) = g \circ {\rm Tr}(f) \circ h$.
\end{lem}
%------------------------------------------------------------------
\begin{proof}
We first prove that for $f \in \Hom_{\cB}(M, N)$, ${\rm Tr}(f)$ is an $\bH$-module homomorphism. Consider the case when $d=2$ first and let $m \in M$ be an arbitrary element. Note that
\begin{equation}    
{\rm Tr}(f)(H_1 \cdot m) 
= \tfrac{1}{2} ( f(H_1 \cdot m) + H_1^{-1} f(H_1^2 \cdot m))
= \tfrac{1}{2} ( f(H_1 \cdot m) + H_1^{-1} f(Z \cdot m)).
\end{equation}
Because $Z\in \cB$ and $f$ is $\cB$-linear,
$H_1^{-1} f(Z \cdot m) = H_1^{-1}Z f(m) = H_1 f(m)$,
and hence
\begin{equation}
\begin{split}
{\rm Tr}(f)(H_1 \cdot m) 
&= \tfrac{1}{2}(H_1 H_1^{-1}f(H_1 \cdot m) + H_1 f(m))
\\
&= H_1(\tfrac{1}{2}(H_1^{-1}f(H_1 \cdot m) + f(m)))
\\
&= H_1 {\rm Tr}(f)(m).
\end{split}
\end{equation}
That is, ${\rm Tr}(f)$ commutes with $H_1$.
For $b \in \cB$, one has
\begin{equation}
\begin{split}
{\rm Tr}(f)(b \cdot m) 
&= \tfrac{1}{2} ( f(b \cdot m) + H_1^{-1} f(H_1 b \cdot m))
\\
&= \tfrac{1}{2} ( b f( m) + H_1^{-1} \sigma(b) f( H_1 \cdot m))
\\
&= b {\rm Tr}(f)(m).
\end{split} 
\end{equation}
In other words, ${\rm Tr}(f) \in \Hom_{\bH}(M, N)$.
This verification generalizes to the case $d>2$ since each $H_w$ is invertible, $d! \in K^\times$ and $H_w H_x \in \cB^\times H_{wx}$, which follows from the fact that $S=0$ and $R \in (B\otimes B)^\times$.

The other two claims that ${\rm Tr}(g \circ f \circ h) = g \circ {\rm Tr}(f) \circ h$ and ${\rm Tr}$ restricts to the identity follow from a direct calculation.
\end{proof}

\begin{cor} \label{cor:exttransfer}
Under the same assumption as in Lemma~\ref{lem:transfer}, the restriction
\[
  \Ext^1_{\bH}(M, N) \to \Ext^1_{\cB}(M, N)
\]
is injective.
\end{cor}
\begin{proof}
{
Let $0\to N\to E\overset{p}{\to}M\to0$ be a short exact sequence of $\bH$-modules that splits over $\cB$. 
Then there is a $\cB$-linear map $s\colon M\to E$ with $p\circ s=\mathrm{id}_M$. 
By Lemma~\ref{lem:transfer}, applied to the $\bH$-modules $M$ and $E$, 
the map ${\rm Tr}(s)\colon M\to E$ is $\bH$-linear. 
Since $p$ and $\mathrm{id}_M$ are $\bH$-linear, the same lemma gives
\[
p\circ{\rm Tr}(s)={\rm Tr}(p\circ s)={\rm Tr}(\mathrm{id}_M)=\mathrm{id}_M .
\]
Hence, the sequence splits over $\bH$.
}
\end{proof}
%=============================================
\section{Dual Specht Theory for Hu Algebras}\label{sec:1-cover}
%=============================================
In this section, 
we will apply the dual Specht theory developed in Section \ref{sec:dST}
to the generalized Hu algebra $\cH_q(\Sigma_m)\wr \cH(d)$ and the wreath Schur functor $\cF_\bH$ constructed in \cite{LNX24} (see Section \ref{sec:WSA}).
The axioms (S1)--(S5) in Definition \ref{def:dST} will be verified in Sections \ref{sec:DSgHu}--\ref{sec:S45}, under an embedding assumption \eqref{eq:goodH} which holds for $d=2$.

For the Hu algebra  $\cH_q(\Sigma_m)\wr \cH(2)$ (and hence  the Hecke algebra of type $\D_{2m}$), its 1-cover is obtained in Section \ref{sec:1coverThm}.
For the generalized Hu algebra $\cH_q(\Sigma_m)\wr \cH(d)$ with $d > 2$, 
it is anticipated that our dual Specht theory applies to the Schur functor associated with the Schur algebras appearing in \cite{HM12},
which will lead to a potential construction of a 1-cover for the Hecke algebra of type $G(d,d,md)$.

%=============================================
\subsection{Construction of wreath Schur functors}\label{sec:WSA}
%=============================================
To verify the axioms (S1)--(S5), one must first construct the inverse Schur functor $\mathcal{G}$. 
Let  $\bH = B \wr \cH(d)$ be the generalized Hu algebra where $B =\cH_q(\Sigma_m)$ and $d \geq 2$.
Recall from Proposition~\ref{prop:H1} that $R = Z \in B\otimes B$ is central and invertible.
This section details the construction of $\cG$ for $\bH$ using the compatibility of Schur dualities for $\bH$ and for type A.

Fix positive integers $m$ and $n$.
Let $V = V(n)$ be the free right $K$-module with basis $\{v_i\}_{1\leq i \leq n}$.
Then, the $K$-module $V^{\otimes m}$ has a basis 
\[
\{v_f:= v_{f(1)} \otimes \cdots \otimes v_{f(m)}~|~ 
f:\{1, \dots, m\} \to \{1, \dots, n\}\},
\]
on which $\Sigma_m$ acts from the right by permuting tensor factors.
There is a unique right module structure over $B=\cH_q(\Sigma_m)$ on $V^{\otimes m}$ determined by the following:
\begin{equation}\label{eq:fundaction}
  v_f \cdot T_i = \begin{cases}
    v_{f}\cdot s_i%(b\cdot s_k) 
    %+ v_{f} \partial^{\beta}_k(b) 
    &\textup{if } f(i) < f(i+1),
    \\
     qv_{f} %(b\cdot s_k) 
     %+ v_{f} \partial^{\beta}_k(b) 
     &\textup{if }  f(i) = f(i+1),
    \\
    qv_{f}\cdot s_i %(b\cdot s_k) 
    + (q-1)v_{f} %(\partial^\beta_k(b) + S_k(b\cdot s_k) 
    &\textup{if }  f(i) > f(i+1).
    \end{cases}
\end{equation}
{
The $B$-module  $V(n)^{\otimes m}$ induces a natural action of $\cB := B^{\otimes d}$ on $(V(n)^{\otimes m})^{\otimes d}$. 
Write $\bT = \bT({n,m,d}) := V(n)^{\otimes md}$ as an $\cH_q(\Sigma_{md})$-module.
Let $S_q(n,m) = \End_{B}(V(n)^{\otimes m})$ be the $q$-Schur algebra of type A.
Denote 
$S^\cB :=  \End_{\cB}(\bT)$
and the {\em wreath Schur algebra} by 
\begin{eq}\label{def:goodWS}
\bS = \bS(n,m\wr d):=\End_{\bH}(\bT(n,m,d)).
\end{eq}    
}
%------------------------------------------------------------------
\begin{lem}\label{lem:summand}
Let $B= \cH_q(\Sigma_m)$, and $\bH = B\wr \cH(d)$ be the generalized Hu algebra with $d\geq 2$.
Suppose that the following condition holds:
\begin{eq}\label{eq:goodH}
\textup{There is an embedding }\iota: \bH \hookrightarrow \cH_q(\Sigma_{md}) \textup{ such that } \iota|_{\cB}\textup{ is the identity map.} 
\end{eq}
Then, one has the following compatible towers of algebras acting on the same tensor space:
\begin{equation}\label{eq:STHdiag}
\begin{array}{cccccccc}
S_q(n,dm)&\crr& V(n)^{\otimes dm} &\crl& \cH_{q}(\Sigma_{dm})
\\
\dsubseteq &&\veq&&\usubseteq
\\
\bS(n,m\wr d) &\crr& \bT(n,m,d) & \crl& \bH
\\
\dsubseteq &&\veq&&\usubseteq
\\
S^\cB \cong S_q(n,m)^{\otimes d}&\crr& (V(n)^{\otimes m})^{\otimes d}& \crl& \cB = \cH_q(\Sigma_m)^{\otimes d}
\end{array}
\end{equation}
Moreover, if in addition $d! \in K^\times$, then $\bS$ is a direct summand of $S^\cB$ as an $\bS$-module.
\end{lem}
%------------------------------------------------------------------
\begin{proof}
Note that \eqref{eq:goodH} implies the following:
\begin{eq}\label{eq:goodT}
\textup{The subalgebra }\cB \subseteq \bH \textup{ acts on }\bT\textup{ factor-wise by the }B\textup{-action on }V(n)^{\otimes m},
\end{eq}
and hence 
$S^\cB \cong S_q(n,m)^{\otimes d}$. The fact that $\bS$ is a direct summand follows from Lemma~\ref{lem:transfer}.
\end{proof}
%------------------------------------------------------------------
 
For $d=2$, \eqref{eq:goodH} holds due to Proposition \ref{prop:H1}. 
Assuming \eqref{eq:goodH} holds for $d > 2$, the following proposition describes the construction of the wreath Schur functor, generalizing \cite[Section 7.6]{LNX24}.
%------------------------------------------------------------------
\begin{prop}\label{prop:WSF}
Let $B= \cH_q(\Sigma_m)$, and $\bH = B\wr \cH(d)$ be the generalized Hu algebra with $d\geq 2$ and $d! \in K^\times$.
Suppose that \eqref{eq:goodH} holds. 
If $n \geq md$, then there is a Schur functor $\cF_\bH : \bS$-mod $\to B\wr \cH(d)$-mod, $M \mapsto eM$ for some idempotent $e = e^2 \in \bS$ such that $e\bS e \cong B\wr \cH(d)$.

Consequently, the inverse Schur functor $\cG: B\wr \cH(d)$-mod $\to \bS(n,m\wr d)$-mod,  $M \mapsto \Hom_\bH(e \bS, M)$ is well-defined. 
\end{prop}
%------------------------------------------------------------------
\begin{proof}
In order to obtain the (inverse) Schur functor, it suffices to check that \cite[Theorem 7.3.1]{LNX24} applies.
In other words, we need to check (i) $\bH$ is a self-injective algebra, and (ii) $\bT$ is a faithful $\bH$-module.

Note that the parameters for the generalized Hu algebras are chosen so that $\sigma = {}$flip, $\rho = 0$,  and $R = Z(m)$ is invertible.
It follows from \cite[Proposition 3.4.1(b)]{LNX24} that the quantum wreath product $\bH$ is symmetric, and hence self-injective (see \cite{CR81}). That is, (i) holds.

For (ii), it suffices to show that $\bT$ has a splitting $\bT \cong \bH \oplus Q'$ as a right $\bH$-module.
By \cite[Lemma 7.6.1]{LNX24}, the existence of such a splitting follows from a splitting $\bT \cong \cH_q(\Sigma_{md}) \oplus Q$ as a right $\cH_q(\Sigma_{md})$-module.
Since $n \geq md$, the submodule of $\bT$ generated by $v_1 \otimes \dots \otimes v_{dm}$ is isomorphic to $\cH_q(\Sigma_{md})$ as an $ \cH_q(\Sigma_{md})$-module.
Therefore, both desired splittings of $\bT$ exist. The condition (ii) holds. 
\end{proof}
%------------------------------------------------------------------
To sum up, we have constructed the Schur functor for the Hu algebra $B\wr \cH(2)$. 
Such a construction generalizes to $d>2$ for any embedding of the form \eqref{eq:goodH}, if one exists.

%=============================================
\subsection{Verification of (S1): Dual Specht Modules}\label{sec:DSgHu}
%=============================================
Within this section, assume that $d \geq 2$ and $d! \in K^\times$.
We will prove that the generalized Hu algebra $B\wr \cH(d)$ admits dual Specht modules, provided that $\ell \geq 4$.

Recall for any $\bblambda \in \Omega$, one can assign a Specht module
\[
S^\bblambda = \ind_{|\bblambda|}^d \big(
\bigotimes\nolimits_{i=1}^r (S^i \wr N^i)
\big),
\]
where 
$S^i := S^{\lambda^{(k_i)}}_m$, $N^i := S^{\nu^{(i)}}_{\chi_i}$,
$\{\lambda^{(k_1)}, \dots, \lambda^{(k_r)}\}$ is the support of $\bblambda$, 
$\chi_i = \chi_{\lambda^{(k_i)}}$, 
and $|\bblambda|$ is the composition of $d$ such that $\nu^{(i)} \vdash |\bblambda|_i$.
Recall $\overline{\bblambda}$ from \eqref{eq:bblambdabar}.
Let $X = \bigotimes_{i=1}^{r} (S^i \wr N^i)$ as above and let $w \in \Sigma_d$.
Then, as a $\cB$-module:
\begin{equation}\label{eq:wTwistB}
{}^{w} X \cong 
\big(
S_m^{\overline{\bblambda}(w^{-1}(1))} \otimes S_m^{\overline{\bblambda}(w^{-1}(2))} \otimes \dots \otimes S_m^{\overline{\bblambda}(w^{-1}(d))}
\big)^{\oplus \prod_i \dim N^i}
\quad
(w \in \Sigma_d).
\end{equation}
Recall $\delta_w := \delta(\mu, w, \eta)$ for some $w \in {}^\mu\Sigma^\eta$.
From \eqref{def:interval}, each interval $\ZZ^{\delta_w}_k$ of $\delta_w$ is of the form $\ZZ^\mu_a \cap w \ZZ^\eta_b$ for some $a$ and $b$.  
Hence, there exists $\lambda' = \lambda'(k)$ and $\lambda'' = \lambda''(k) \vdash m$ such that
\[
\overline{\bblambda}(i) = \lambda'
\quad
\textup{and}
\quad
\overline{\bbmu}(w^{-1}(i)) = \lambda''
\quad
\textup{for all}
\quad
i \in \ZZ^{\delta_w}_k.
\]
By \eqref{eq:wTwistB} and its analog for $Y$, 
for any $i \in \ZZ^{\delta_w}_k$, the $i$th tensor factor of $X$ (and ${}^wY$, resp.) is $S_m^{\lambda'}$ (and $S_m^{\lambda''}$, resp.). 
For $\nu \vdash m$, we have
\begin{eq}\label{def:nnu}
n_\nu(w) = \#\{ i ~|~ \overline{\bbmu}(i) = \nu = \overline{\bblambda}(w(i)) \} = \sum\nolimits_{k; \lambda'(k) = \nu = \lambda''(k)} \#\ZZ^{\delta_w}_k,
\end{eq}
and at most one interval $\ZZ^{\delta_w}_k$ contributes to this sum.
%------------------------
\begin{Def}\label{def:HudS}
For $\bblambda \in \Omega$, let
$S_\bblambda := (S^\bblambda)^\ast \cong (S^{\bblambda^t})^\#$.
\end{Def}
This implies that $S_{\bblambda^t} \cong (S^\bblambda)^\#$.
We will show that $\{S_\bblambda\}_\bblambda$ forms a set of dual Specht modules in the sense of Definition \ref{def:dS}.
Let us first establish the following lemma regarding the extension groups between Specht modules.
%--------------------------------------------------------------------
\begin{lem}\label{lem:gHuExt1}
Let $\bblambda, \bbmu \in \Omega$. Assume that $d! \in K^\times$.
\begin{enua}
\item For all $n \geq 0$, $\Ext^n_\bH(S^\bblambda, S^\bbmu) \cong \Ext^n_\bH(S_{\bblambda^t}, S_{\bbmu^t})$.
\item For all $n \geq 0$, 
\begin{equation}\label{eq:ExtHSlSm}
\Ext^n_\bH(S^\bblambda, S^\bbmu)
\cong 
\Ext^n_{|\bblambda|}\big( 
\bigotimes\nolimits_i (S^i\wr N^i),  \ind^d_{|\bbmu|} (\bigotimes\nolimits_j ({}'S^j\wr {}'N^j))
\big),
\end{equation}
where $S^\bbmu = \ind_{|\bbmu|}^d \big(
\bigotimes\nolimits_{j=1}^{r'} ({}'S^j \wr {}'N^j)
\big)$,
 with ${}'S^j:= S^{\lambda^{(k'_j)}}_m$,
 ${}'N^j :=  S^{\bbmu(\lambda^{(k'_j)})}_{\chi'_j}$, 
 $\{k'_1, \dots, k'_{r'}\}$ is the support of $\bbmu$, 
 $\chi'_j = \chi_{\lambda^{(k'_j)}}$,
 and $|\bbmu|$ is the composition of $d$ such that $\bbmu(\lambda^{(k'_j)}) \vdash |\bbmu|_j$.
\item Suppose that $\ell\geq 3$, and let $w \in {}^\mu\Sigma^\eta$. If $\Hom_{\cB}(X, {}^wY) \neq 0$, then $w$ satisfies \textup{(K1)}. 
If $\Hom_{\delta_w}(X, {}^wY) \neq 0$, then $w$ satisfies \textup{(K1)} and \textup{(K2)}. 
In particular, $\bbmu \leq_K \bblambda$.
\item If $\ell \geq 3$ then $\End_\bH(S_\bblambda) \cong K$.
\item Suppose that $\ell \geq 3$. 
Then $\Hom_{\bH}(S^\bblambda, S^\bbmu) \neq 0$ implies that $\bblambda \geq_K  \bbmu$. 
\item Suppose that $\ell \geq 4$. 
Then $\Ext^1_{\bH}(S^\bblambda, S^\bbmu) \neq 0$ implies that $\bblambda >_K \bbmu$. 
\end{enua}
\end{lem}
%--------------------------------------------------------------------
\begin{proof}
Part (a) follows from the fact that 
$\Ext^n_\bH(S^\bblambda, S^\bbmu) 
\cong \Ext^n_\bH((S^\bblambda)^\#, (S^\bbmu)^\#)
\cong \Ext^n_\bH(S_{\bblambda^t}, S_{\bbmu^t})$.
Part (b) follows from Frobenius reciprocity.

For part (c), 
note that
the assumption $d! \in K^\times$
implies that $\tcH^{\chi}_c \cong K\Sigma_c$ is split semisimple for any $c \leq d$. 
Hence, one can apply  the branching rule to obtain that 
$X$,
as a $B\wr\cH(\delta_w)$-module, is  a direct sum of modules of the form 
\begin{equation}\label{eq:Xdecomp}
\bigotimes\nolimits_k 
(S^{\lambda'(k)}_m \wr S^{\beta'_k}_{\chi_{\lambda'(k)}})
\quad
\textup{with}
\quad\beta'_k \subseteq \bblambda(\lambda'(k)).
\end{equation}
Similarly, ${}^wY$ is a direct sum of modules of the form $\bigotimes_k (S^{\lambda''(k)}_m \wr S^{\beta''_k}_{\chi_{\lambda''_k}})$ with $\beta''_k  \subseteq \bbmu(\lambda''(k))$.

Since $\End_B(S^\nu) \cong K$ for all $\nu \vdash m$ for $\ell \geq 3$, it follows from \cite[Proposition 5.3.1]{LNX25}
that
\[
\Hom_{B\wr\cH(c)}(S^\nu_m \wr S^{\beta'}_{\chi_\nu}, S^\nu_m \wr S^{\beta''}_{\chi_\nu}) 
\cong 
\Hom_{\tcH^{\chi_\nu}_c}(S^{\beta'}_{\chi_\nu}, S^{\beta''}_{\chi_\nu}),
\]
for all $\nu \vdash m$, $c \geq 1$ and $\beta', \beta'' \vdash c$.
Since $\tcH^{\chi_\kappa}_c$ is split semisimple, this space is nonzero only if $\beta' = \beta''$.

Now, suppose that $\Hom_{\cB}(X, {}^wY) \neq 0$. 
By \eqref{eq:wTwistB} and the K\"unneth formula over $\cB$, 
\[
\Hom_B(S^{\overline{\bblambda}(p)}_m, S^{\overline{\bbmu}(w^{-1}(p))}_m) \neq 0
\quad
\textup{for all}
\quad
1 \leq p \leq d.
\]
Since $\ell \geq 3$, \cite[Lemma 4.4.1]{HN04} gives $\overline{\bblambda}(p) \unrhd \overline{\bbmu}(w^{-1}(p))$ for all $p$. 
Substituting $p = w(i)$ yields
\[
\overline{\bbmu}(i) \unlhd \overline{\bblambda}(w(i)) \quad\textup{for all}\quad 1 \leq i \leq d,
\]
which is (K1).
Next, suppose that $\Hom_{\delta_w}(X, {}^wY) \neq 0$. 
Then $\Hom_\cB(X, {}^wY) \neq 0$, and hence (K1) holds. 
Let $\nu \vdash m$. 
If there is no interval $\ZZ^{\delta_w}_k$  satisfying $\lambda'(k) = \nu = \lambda''(k)$, 
then $n_\nu(w) = 0$ and (K2) holds for $\nu$ with $\beta' = \varnothing$. 
Otherwise, let $\ZZ^{\delta_w}_k$ be the unique interval contributing to $n_\nu(w) = |\ZZ^{\delta_w}_k|$. 
By the K\"unneth formula over $B\wr\cH(\delta_w)$, there are $\beta'$, $\beta''$ satisfying
\[
\Hom_{B\wr\cH(|\ZZ^{\delta_w}_k|)}(S^\nu_m \wr S^{\beta'(k)}_{\chi_\nu}, 
S^\nu_m \wr S^{\beta''}_{\chi_\nu}) \neq 0,
\]
and hence $\beta' = \beta'' \vdash n_\nu(w)$ is contained in both $\bblambda(\nu)$ and $\bbmu(\nu)$, which is (K2). 
This proves (c).

For part (d), setting $n=0$ and $\bbmu = \bblambda$ (so $\eta = \mu = |\bblambda|$) in (b), we apply the Mackey filtration (Proposition \ref{prop:Mackey}) to $S^\bblambda = \ind_\mu^d X$.
  Suppose that $\Hom_\cB(X, {}^wX) \neq 0$ for some $w \in {}^\mu\Sigma^\mu$.
  By (K1) in (c), $\overline{\bblambda}(k) \unlhd \overline{\bblambda}(w(k))$ for all $1\leq k \leq d$, and hence
  \[
  \overline{\bblambda}(k) \unlhd \overline{\bblambda}(w(k)) \unlhd \dots \unlhd \overline{\bblambda}(w^{N_w}(k)) = \overline{\bblambda}(k),
  \]
  where $N_w$ is the order of $w$. 
  That is, $\overline{\bblambda}(k) = \overline{\bblambda}(w^{-1}(k))$ for all $k$, 
  which implies that $w$ lies in the stabilizer $\Sigma_\mu$ of $\overline{\bblambda}$. 
  It then follows from ${}^\mu\Sigma^\mu \cap \Sigma_\mu = \{1\}$ that $w = 1$.

%---
Thus $\Hom_\cB(X, {}^w X) = 0$ for all $w \neq 1$, and only the identity double coset $w = 1$ contributes to $\End_\bH(S^\bblambda)$,
by Corollary \ref{cor:MackeyHomExt}(a).
Therefore,
\[
\End_\bH(S^\bblambda) 
\cong \End_{B\wr\cH(\mu)}(X) 
\cong \bigotimes\nolimits_{i=1}^r \End_{B\wr\cH(\mu_i)}(S^i \wr N^i).
\]
For each $i$, applying \cite[Proposition 5.3.1]{LNX25} yields $\End_{B\wr\cH(\mu_i)}(S^i \wr N^i) \cong K$ since $\End_B(S_m^\lambda) \cong K$ for $\ell \geq 3$.
Hence, $\End_\bH(S_\bblambda) \cong \End_\bH(S^\bblambda) \cong K$.

For part (e), suppose that $\Hom_\bH(S^\bblambda, S^\bbmu) \neq 0$.
By Corollary \ref{cor:MackeyHomExt}(a), there exists $w \in {}^\mu\Sigma^\eta$ such that $\Hom_{\delta_w}(X, {}^w Y) \neq 0$.
By (c), $w$ satisfies (K1) and (K2) and hence  $\bblambda \geq_K \bbmu$.

For part (f), suppose that $\Ext^1_\bH(S^\bblambda, S^\bbmu) \neq 0$. 
By (b) for $n = 1$ and Proposition \ref{prop:Mackey}, there exists $w \in {}^\mu\Sigma^\eta$ such that $\Ext^1_{\mu}(X, \ind^\mu_{\delta_w}({}^wY)) \neq 0$, and hence $\Ext^1_{\delta_w}(X, {}^wY) \neq 0$.
We analyze the non-vanishing extensions through the following two foundational cases before treating the general case.
In each case, we verify (K1) for $w$ with at least one strict dominance, and (K2) for $w$.
In other words, this gives $\bbmu \leq_K \bblambda$. 
Moreover, one can show that $\bbmu \neq \bblambda$. Otherwise, 
(K1) would give $\overline{\bblambda}(i) \unlhd \overline{\bblambda}(w(i))$ for all $i$ with a strict dominance for some $i$, contradicting the cycle argument in part (d). 

\textbf{Case 1: Purely base-induced modules ($\bblambda = \sum_{i=1}^d e_{\lambda^{(i)}}^{(1)}$ for distinct $\lambda^{(i)} \in I_B$).}

We may label the $\lambda^{(i)}$ so that $\overline{\bblambda}(i) = \lambda^{(i)}$.
Here $\mu = (1^d)$, so $\delta_w = (1^d)$, $B\wr\cH(\delta_w) = \cB$ and $X = S_m^{\lambda^{(1)}} \otimes \dots \otimes S_m^{\lambda^{(d)}}$. 
By \eqref{eq:wTwistB}, as a $\cB$-module, ${}^wY$ is a direct sum of copies of $S_m^{\overline{\bbmu}(w^{-1}(1))} \otimes \dots \otimes S_m^{\overline{\bbmu}(w^{-1}(d))}$.
Since $\Ext^1_{\cB}(X, {}^wY) \neq 0$, the K\"unneth formula gives an index $j \in \{1, \dots, d\}$ such that
\[
\Ext^1_B(S_m^{\lambda^{(j)}}, S_m^{\overline{\bbmu}(w^{-1}(j))}) \neq 0,
\quad\text{and}\quad
\Hom_B(S_m^{\lambda^{(p)}}, S_m^{\overline{\bbmu}(w^{-1}(p))}) \neq 0 \quad\text{for all } p \neq j.
\]
Because $\ell \geq 3$, \cite[Lemma 4.4.1]{HN04} implies that $\lambda^{(p)} \unrhd \overline{\bbmu}(w^{-1}(p))$ in $I_B$ for all $p \neq j$.
Furthermore, because $\ell \geq 4$, \cite[Proposition 4.2.1]{HN04} implies that $\Ext^1_B(S_m^{\lambda^{(j)}}, S_m^{\overline{\bbmu}(w^{-1}(j))}) \neq 0$ forces the strict dominance $\lambda^{(j)} \rhd \overline{\bbmu}(w^{-1}(j))$.
Substituting $p = w(i)$ yields (K1) with a strict dominance for $i = w^{-1}(j)$.
As for (K2), note that $n_\nu(w) \leq \#\bblambda(\nu) \leq 1$ for all $\nu \vdash m$, and if $n_\nu(w) = 1$, then $\beta = (1)$ is contained in both $\bblambda(\nu) = (1)$ and $\bbmu(\nu) \neq \varnothing$. 
Thus, (K2) holds.

\textbf{Case 2: Single wreath blocks ($\bblambda = e_\lambda^\nu$ and $\bbmu = e_{\lambda'}^{\nu'}$ with $\lambda, \lambda' \vdash m$ and $\nu, \nu' \vdash d$).}

In this case, $\mu = \eta = (d)$, so $w = 1$, $S^\bblambda = S_m^\lambda \wr S_{\chi_\lambda}^\nu$ and $S^\bbmu = S_m^{\lambda'} \wr S_{\chi_{\lambda'}}^{\nu'}$.
One can apply \cite[Lemma 5.5.1]{LNX25} and obtain that either 
\begin{itemize}
    \item[(i)] $\Ext^1_B(S_m^\lambda, S_m^{\lambda'}) \neq 0$; or
    \item[(ii)] $\Hom_B(S_m^\lambda, S_m^{\lambda'}) \neq 0$ and
$\Ext^1_{\tcH}(S_\chi^\nu, S_\chi^{\nu'}) \neq 0$.
\end{itemize} 
Note that the twisted Hecke algebra $\tcH$ is split semisimple.
Hence, $\Ext^1_{\tcH}(S_{\chi_\lambda}^\nu, S_{\chi_{\lambda'}}^{\nu'}) = 0$.
Thus, the condition (ii) cannot occur.
It follows from (i) that $\lambda \rhd \lambda'$ by \cite[Proposition 4.2.1]{HN04} since $\ell \geq 4$.
Therefore, (K1) holds for $w = 1$, with $\overline{\bbmu}(i) = \lambda' \lhd \lambda = \overline{\bblambda}(i)$ for all $i$. 
Moreover, (K2) holds trivially, since $\lambda \neq \lambda'$ forces $n_\nu(1) = 0$ for all $\nu \vdash m$.

\textbf{General Case: Arbitrary $\bblambda, \bbmu \in \Omega$.}
Decomposing $X$ and ${}^wY$ as in \eqref{eq:Xdecomp}, 
the K\"unneth formula over $B\wr\cH(\delta_w)$ yields, for some interval $\ZZ^{\delta_w}_k$ of $\delta_w$:
\[
\begin{split}
\Ext^1_{B\wr\cH(|\ZZ^{\delta_w}_k|)}\big(
S^{\lambda'(k)}_m \wr S^{\beta'(k)}_{\chi_{\lambda'(k)}}, S^{\lambda''(k)}_m \wr S^{\beta''(k)}_{\chi_{\lambda''(k)}}
\big) &\neq 0,
\\
\Hom_{B\wr\cH(|\ZZ^{\delta_w}_l|)}\big(
S^{\lambda'(l)}_m \wr S^{\beta'(l)}_{\chi_{\lambda'(l)}}, S^{\lambda''(l)}_m \wr S^{\beta''(l)}_{\chi_{\lambda''(l)}}
\big) &\neq 0
\quad
\textup{for all}
\quad
l\neq k.
\end{split}
\]
If every $|\ZZ^{\delta_w}_k|$ is 1, then the statement can be deduced via Case 1.

For each $\ZZ^{\delta_w}_k$ with $|\ZZ^{\delta_w}_k|\geq 2$, one applies Case 2 to $B\wr\cH(|\ZZ^{\delta_w}_k|)$ in place of $B\wr\cH(d)$.
One obtains $\lambda'(k) \rhd \lambda''(k)$. 
For $l \neq k$, restricting to $B^{\otimes |\ZZ^{\delta_w}_l|}$ and applying \cite[Lemma 4.4.1]{HN04} give $\lambda'(l) \unrhd \lambda''(l)$.
Therefore, $\overline{\bbmu}(i) \unlhd \overline{\bblambda}(w(i))$ for all $i$, 
with a strict dominance whenever $w(i) \in \ZZ^{\delta_w}_k$.
In other words, (K1) holds.
Since $\lambda'(k) \neq \lambda''(k)$, the interval $\ZZ^{\delta_w}_k$ does not contribute to any $n_\nu(w)$. 
If $\ZZ^{\delta_w}_l$ where $l\neq k$ is an interval with $\lambda'(l) = \nu = \lambda''(l)$, then $\beta'(l) = \beta''(l)$, which has size $|\ZZ^{\delta_w}_l| = n_\nu(w)$ and is contained in both $\bblambda(\nu)$ and $\bbmu(\nu)$. 
Consequently, (K2) holds exactly as in (c).
\end{proof}
The lemma can be summarized via this proposition, which assures (S1).
%--------------------------------------------------------------------
\begin{prop}\label{prop:gHuDS}
Suppose that $\ell  \geq 4$ and $d!\in K^\times$.
Then, the generalized Hu algebra $\bH = \cH_q(\Sigma_m) \wr \cH(d)$ admits dual Specht modules $\{S_\bblambda\}_{\bblambda \in \Omega}$ as defined in Definition \ref{def:HudS}.
\end{prop}
%--------------------------------------------------------------------
\begin{proof}
It follows from Proposition \ref{prop:MwrNHu} that the Specht modules 
$\{S^\bblambda\}_{\bblambda \in \Omega}$ constructed in \eqref{def:SpechtgHu} are indeed modules over $\bH$.
Thanks to Proposition \ref{prop:genSpechtiso}, the dual Specht modules $\{S_\bblambda\}_{\bblambda \in \Omega}$ are well-defined.
The poset structure $(\Omega, \leq_K)$ is described in Definition \ref{def:KPO}.
The axioms are all verified in Lemma \ref{lem:gHuExt1}.
\end{proof}
%--------------------------------------------------------------------
\subsection{Verification of (S2): Projectivity} 
%--------------------------------------------------------------------
%Let $\bblambda \in \Omega$.
%Theorem \ref{thm:Young} leads to a construction of the Young module $Y^\bblambda$.
%Note that $Y^\bblambda$ is indecomposable provided $S^\bblambda$ is indecomposable,
%which reduces to the requirement that $S^\lambda_m$ is also indecomposable, or, $\ell \geq 3$.
From now on,
assume that \eqref{eq:goodH} holds.
In other words,
the inverse Schur functor $\cG: \bH$-mod $\to \bS$-mod is well-defined, with respect to the wreath Schur algebra $\bS$ as in \eqref{def:goodWS}.

Since $\{S_\bblambda\}_\bblambda$ are now dual Specht modules, it follows from Definition \ref{def:YoungMod} that the Young modules $Y^\bblambda$ can be constructed, with a dependence on a choice of minimal elements.
Let $P(\bblambda) := \cG(Y^\bblambda)$.
It follows from Proposition \ref{prop:indec} that  $P(\bblambda)$ is indecomposable.
%--------------------------------------------------------------------
\begin{prop} \label{prop:PIM}
Suppose that  $n \geq md$, $\ell \geq 4$, and $d! \in K^\times$. 
Then the indecomposable $\bS$-module $P(\bblambda)$ is projective for any $\bblambda\in \Omega$. 
\end{prop}
%--------------------------------------------------------------------
\begin{proof} 
Set $\overline{\cG} = \operatorname{Hom}_{\cB}(V(n)^{\otimes dm},-) : \cB\textup{-mod} \to S^\cB$-mod. 
By applying the methods in \cite[Theorem 3.4.1]{HN04}, 
one can show that $\overline{\cG}( (Y^{\bblambda})^\#)$ has both a $\nabla$- and a $\Delta$-filtration, as an $S^\cB$-module for all $\bblambda\in \Omega$. 
Therefore, $\overline{\cG}( (Y^{\bblambda})^\#) \cong \bigoplus_{i} T(\bbmu_{i})$ is a direct sum of tilting $S^\cB$-modules. 
Now apply the Schur functor to get that 
$(Y^{\bblambda})^\# \cong \bigoplus_{i} (Y^{\bbmu_{i}^{t}})^\#$ as a $\cB$-module. 
It follows that $Y^{\bblambda}$ is a direct summand of $\bigoplus V(n)^{\otimes dm}$, 
and $\overline{\cG}(Y^{\bblambda})$ is a direct summand of $\bigoplus \overline{\cG}(V(n)^{\otimes dm})\cong \bigoplus S^\cB$. 
Consequently, $\overline{\cG}(Y^{\bblambda})$ is a projective $S^\cB$-module. 
Next, one has 
\begin{equation}
  \cG(Y^{\bblambda})
  \cong \Hom_{\bH}(V(n)^{\otimes dm}, Y^{\bblambda})
  \cong \Hom_{\cB}(V(n)^{\otimes dm}, Y^{\bblambda})^{\Sigma_{d}}
  \cong \overline{\cG}(Y^{\bblambda})^{\Sigma_{d}} .
\end{equation}
Since $n \geq md$,
a decomposition of the tensor space into permutation modules shows that $\bH\otimes_\cB\bT$ is a direct summand of $\bT^{\oplus N}$ for some $N$.
Hence, $S^\cB = \End_\cB(\bT) \cong \Hom_\bH(\bT, \bH\otimes_\cB \bT)$ is projective as a left $\bS$-module, and $\cG(Y^\bblambda) = {\rm Tr}(\overline{\cG}(Y^\bblambda))$ is an $\bS$-direct summand of a projective $\bS$-module by Lemma~\ref{lem:transfer}.
Therefore, 
$\cG(Y^{\bblambda})$ 
%is a direct summand of $\bigoplus (S^\cB)^{\Sigma_{d}} \cong \bigoplus \bS$ which means that $\cG(Y^{\bblambda})$ 
is a projective $\bS$-module. 
\end{proof}
%------------------------------------------------------------------
\subsection{Verification of (S3): \texorpdfstring{$R^1\cG$}{R1G} Vanishing} 
%--------------------------------------------------------------------
This section is dedicated to the proof of (S3).
%, which relies on the duality
%$S_{\bblambda^t} \cong (S^{\bblambda})^\#$ between Specht and dual Specht modules.
%------------------------------------------------------------------
\begin{thm}\label{thm:R1G}
Suppose that  $n \geq md$.
 If $\ell \geq 4$ and $d! \in K^\times$, then
$R^1\cG(S_\bblambda) = 0$ for all $\bblambda \in \Omega$.
\end{thm}
%------------------------------------------------------------------
\begin{proof}
From the construction of the wreath Schur functor 
(Proposition \ref{prop:WSF}), it suffices to show that $\Ext^1_\bH(e \bS, S_\bblambda) = 0$, where $e \bS \cong \bT = (V(n)^{\otimes m})^{\otimes d}$. Since the restriction map $\Ext^1_\bH(e \bS, S_\bblambda) \to \Ext^1_{\cB}(e \bS, S_\bblambda)$ is injective by Corollary~\ref{cor:exttransfer}, it suffices to prove that $\Ext^1_\cB(e \bS, S_\bblambda) = 0$.

Now, as a $\cB$-module, $\bT \cong (V^{\otimes m})^{\otimes d}$ where $V^{\otimes m} \in B$-mod is the tensor space for $B$ described in \eqref{eq:fundaction}.
On the other hand, the restriction of $S^{\bblambda^t} $ to $\cB$  decomposes into a finite direct sum of outer tensor products of the form:
\[
S^{\nu_1} \otimes S^{\nu_2} \otimes \cdots \otimes S^{\nu_d},
\quad
(\nu_i \vdash m).
\]
Thus, as a $\cB$-module, $S_\bblambda \cong (S^{\bblambda^t})^\#$ is a direct sum of the form $S_{\nu_1} \otimes S_{\nu_2} \otimes \cdots \otimes S_{\nu_d}$, where $S_\nu := (S^\nu_m)^\ast \cong (S^{\nu^t}_m)^\# \in B$-mod.
The K\"unneth formula implies that 
\[
\Ext^1_{\cB}((V^{\otimes m})^{\otimes d}, S_{\nu_1} \otimes S_{\nu_2} \otimes \cdots \otimes S_{\nu_d})
 \cong 
 \bigoplus_{i=1}^d 
 \left( 
 \Ext^1_B(V^{\otimes m}, S_{\nu_i}) 
 \otimes 
 \bigotimes_{j \neq i} 
 \Hom_B(V^{\otimes m}, S_{\nu_j}) 
 \right).
\]
Since $\ell \geq 4$, it follows from \cite{HN04} that every $\Ext^1_B(V^{\otimes m}, S_{\nu_i}) = R^1\cG_{\Hom_B}(S_{\nu_i})$ vanishes, and hence 
$R^1\cG(S_\bblambda) = 0$.
\end{proof}
%--------------------------------------------------------------------
%--------------------------------------------------------------------
\subsection{Verification of (S4)--(S5): Classification of Projective Indecomposables} \label{sec:S45}
%--------------------------------------------------------------------
 Let $\overline{\Omega} = \Pi_m^d$, which indexes the Specht modules over $\cB = \cH_q(\Sigma_m)^{\otimes d}$.
 Any multipartition $\overline{\bblambda} \in \overline{\Omega}$
 can be identified as a map  $\overline{\bblambda}: \{1, \dots, d\} \to \Pi_m$. 
 For $\overline{\bblambda} \in \overline{\Omega}$,  denote the Specht module for $\cB$ and the standard module for $S^\cB$, respectively, by
 \begin{equation}
 \overline{S}^{\overline{\bblambda}} := 
 S^{\overline{\bblambda}(1)}_m \otimes  \dots \otimes S^{\overline{\bblambda}(d)}_m \in \cB\textup{-mod},
 \quad
 \overline{\Delta}(\overline{\bblambda}) 
 := \Delta(\overline{\bblambda}(1)) \otimes \dots \otimes \Delta(\overline{\bblambda}(d)) \in S^\cB\textup{-mod}.
 \end{equation} 
 Recall from \eqref{eq:STHdiag} that $\bS$ is a subalgebra of $S^\cB$. 
 In order to show that $P(\bblambda)$, $\bblambda \in \Omega$ exhaust all projective indecomposables in $\bS$-mod, 
 one needs a lemma that compares the standard module $ \overline{\Delta}(\overline{\bblambda}) 
 \in S^\cB$-mod with the standard modules $\Delta(\bblambda) := \cG(S_\bblambda) \in \bS$-mod.
 
% Note that as a $\cB$-module, 
% $ \overline{S}^{\overline{\bblambda}} \cong  
% S^{\lambda^{w(1)}}_m \otimes  \dots \otimes S^{\lambda^{w(d)}}_m$ for all $w\in \Sigma_d$, 
% and so it is convenient to use the $\Sigma_d$-orbit $[\overline{\bblambda}] := \{(\lambda^{(w(i))})_{1\leq i \leq d}~|~ w\in \Sigma_d\}$ to index Specht (and simple) modules over $\cB$.
Suppose that the support of $\overline{\bblambda}$ is $\{\lambda^{(1)}, \dots, \lambda^{(r)}\}$ with $\#\{ 1\leq j \leq d ~|~ \overline{\bblambda}(j) = \lambda^{(i)}\} = \mu_i > 0$ 
{as in \eqref{eq:bblambdabar}}.
%We may assume that $\mu = (\mu_1, \dots, \mu_r) \vdash d$ is a partition.
Define
\[
\Pi(\overline{\bblambda}) := \big{\{}
\overline{\bbnu} = (\nu^{(1)}, \dots, \nu^{(r)})
~|~
\nu^{(i)} \vdash \mu_i \textup{ for all }i
\big{\}}.
\]
Each $\overline{\bbnu} \in \Pi(\overline{\bblambda})$ is associated with a multipartition $\bbnu = \bbnu(\bblambda) \in \Omega$ given by
\[
\bbnu = e^{\nu^{(1)}}_{\lambda^{(1)}} + e^{\nu^{(2)}}_{\lambda^{(2)}} + \dots + e^{\nu^{(r)}}_{\lambda^{(r)}}.
\]
%--------------------------------------------------------------------
  \begin{lem}\label{lem:WeylCompare}  
  Suppose that  $n \geq md$ and $d!\in K^\times$.
  If $\ell \geq 3$, then, as an $\bS$-module,
  \[
  \overline{\Delta}(\overline{\bblambda}) 
  \cong
\bigoplus\nolimits_{ \overline{\bbnu} \in \Pi(\overline{\bblambda})} 
  \Delta(\bbnu)
  ^{\oplus (\prod_{i=1}^r \dim S_{\mu_i}^{\nu^{(i)}})}.
  \]
\end{lem}
%--------------------------------------------------------------------
\begin{proof}
%For brevity, in this proof we write $\varphi_i:= e^{\nu}_{\lambda^{(i)}}$.
By \cite[Proposition 5.4.1]{LNX25}, the induced module $\ind_\cB^\bH( \overline{S}^{\overline{\bblambda}})$ decomposes into a direct sum of Specht modules as follows:
\begin{equation}
  \ind_\cB^\bH \big( \overline{S}^{\overline{\bblambda}} \big)
  \cong
  \bigoplus\nolimits_{\overline{\bbnu} \in \Pi(\overline{\bblambda})}
  (S^{\bbnu})^{\oplus (\prod_{i=1}^r \dim S^{\nu^{(i)}}_{\mu_i})},
\end{equation}
in which the semisimplicity of the twisted Hecke algebra $\tcH_{\mu_i}^{\chi_{\lambda^{(i)}}}$ (which follows from combining the fact that  $\mu_i \leq d$ and the assumption  $d!\in K^\times$) is used.
Since $(S^\lambda_m)^* \cong (S^{\lambda^t}_m)^\#$, for any $\lambda \vdash m$, 
\begin{equation}
\Delta(\lambda) 
= \cG((S^\lambda_m)^*) \cong \cG((S^{\lambda^t}_m)^\#)
= \Hom_B(V(n)^{\otimes m}, (S^{\lambda^t}_m)^\#).
\end{equation}
Then, one has, as an $\bS$-module,
$\overline{\Delta}(\overline{\bblambda}) 
\cong 
\text{Hom}_{\cB}\big(
V(n)^{\otimes dm}, (\overline{S}^{\overline{\bblambda}^t})^\#
\big)$,
where $\overline{\bblambda}^t(i):=\overline{\bblambda}(i)^t$ for all $i$.
Hence,
\begin{equation}
\begin{split}
\overline{\Delta}(\overline{\bblambda}) 
%&
%\cong \cG(S_{(\overline{\bblambda}(1))^t}) \otimes \cdots \otimes \cG(S_{(\overline{\bblambda}(d))^t})
%\\
&\cong \text{Hom}_{\bH}\big(
V(n)^{\otimes dm}, \ind_{\cB}^{\bH} (\overline{S}^{\overline{\bblambda}^t})^\#
\big) 
\\
&\cong \text{Hom}_{\bH}\big(
V(n)^{\otimes dm}, 
\bigoplus\nolimits_{\overline{\bbnu} \in \Pi(\overline{\bblambda})}
((S^{\bbnu^t})^{\oplus (\prod_{i=1}^r \dim S^{\nu^{(i)}}_{\mu_i})}
)^\#
\big) 
\\
&\cong 
\bigoplus\nolimits_{\overline{\bbnu} \in \Pi(\overline{\bblambda})}
\Delta( \bbnu )^{\oplus (\prod_{i=1}^r \dim S^{\nu^{(i)}}_{\mu_i})}.
\end{split} 
\end{equation}
\end{proof}
One can calculate these decompositions in the  following examples. 
%--------------------------------------------------------------------
\begin{exa}
\begin{enua}
\item
Let $d=5$, and let $\overline{\bblambda} = (\lambda, \lambda, \nu, \lambda, \nu)$ for some $\lambda, \nu \in \Pi_m$ such that $\lambda \neq \nu$. Then, $\mu = (3,2) \vdash 5$, and
\begin{equation}
\begin{split}
\ind_{B^{\otimes 5}}^{B\wr \cH(5)} \big( \overline{S}^{\overline{\bblambda}} \big)
&\cong
  S^{e^{(3)}_{\lambda} + e^{(2)}_{\nu}}
  \oplus S^{e^{(3)}_{\lambda} + e^{(1,1)}_{\nu}}
  \oplus (S^{e^{(2,1)}_{\lambda} + e^{(2)}_{\nu}})^{\oplus 2}
  \\
  &\oplus (S^{e^{(2,1)}_{\lambda} + e^{(1,1)}_{\nu}})^{\oplus 2}
  \oplus S^{e^{(1,1,1)}_{\lambda} + e^{(2)}_{\nu}}
  \oplus S^{e^{(1,1,1)}_{\lambda} + e^{(1,1)}_{\nu}}.
  \end{split}
\end{equation}
\item
Let $d=2$. From the identification in Example \ref{ex:Huindex},  we see that
\begin{equation}
\overline{\Delta}(\lambda, \mu) \cong
\begin{cases}
\Delta([\lambda,\lambda]_+) \oplus \Delta([\lambda,\lambda]_-)&\tif \lambda = \mu;
\\
\Delta([\lambda,\mu])&\tif \lambda \neq \mu.
\end{cases}
\end{equation}
\end{enua}
\end{exa}
%----------------------------------------------------------
\subsection{Projective Indecomposable $\bS$-Modules}
%----------------------------------------------------------
We are now ready to give a complete description of the projective indecomposable $\bS$-modules. 
%--------------------------------------------------------------------
\begin{theorem} \label{thm:PL}
Suppose that $\ell \geq 4$, $n \geq md$ and $d! \in K^\times$.
\begin{enua}
\item If $P \in \bS$-mod is projective indecomposable, then $P\cong P(\bblambda)$ for some $\bblambda \in \Omega$. 
\item Let $\bblambda, \bbmu \in \Omega$. Then, $P(\bblambda)\cong P(\bbmu)$ if and only if $\bblambda=\bbmu$. 
\end{enua} 
Consequently, $\{L(\bblambda) := \hd P(\bblambda)~|~ \bblambda\in \Omega\}$ is a complete set of non-isomorphic irreducible representations of $\bS$. 
\end{theorem} 
%--------------------------------------------------------------------
\begin{proof} (a)
For each $\overline{\bblambda} \in \overline{\Omega}$, let $\overline{P}(\overline{\bblambda})$ be the projective cover of 
$\overline{L}(\overline{\bblambda})
:=L(\overline{\bblambda}(1)) \otimes \dots \otimes L(\overline{\bblambda}(d)) \in S^\cB$-mod. 
Consider the surjective $\bS$-module homomorphism $\overline{\Delta}(\overline{\bblambda})\rightarrow \overline{L}(\overline{\bblambda})$. 
By Lemma \ref{lem:WeylCompare}, one sees that as an $\bS$-module,
 $\overline{L}(\overline{\bblambda}) = \hd \overline{\Delta}(\overline{\bblambda})$ consists of irreducible modules of the form $L(\bbnu)$ with $\bbnu \in \Omega$.

By Lemma \ref{lem:summand}, $\bS$ is a direct summand of $S^\cB$.
Since any projective indecomposable $P$ is a direct summand of $\bS$, 
it follows that $P$ is a direct summand of $\overline{P}(\overline{\bblambda})$ for some $\overline{\bblambda} \in \overline{\Omega}$. 
%Recall that the head of $\overline{\Delta}(\overline{\bblambda})$ is $\overline{L}(\overline{\bblambda})$. 
%We are now ready to verify (a). Let $P$ be a projective indecomposable module for $\bS$. 
%Then, $P$ is a direct summand of $\overline{P}(\overline{\bblambda})$ for some $\overline{\bblambda} \in \overline{\Omega}$. 
If $D$ is in the head of $\overline{P}(\overline{\bblambda})$ as an $\bS$-module, 
then $D$ is in the head of $\overline{L}(\overline{\bbmu})$ where  $\overline{L}(\overline{\bbmu})$ is a composition factor of $\overline{P}(\overline{\bblambda})$. 
Consequently, $\hd P \cong L(\bbmu)$ for some $\bbmu\in \Omega$, 
and $P\cong P(\bbmu)$. 

(b) 
%If $P(\bblambda) \cong P(\bbmu)$, applying the Schur functor $\cF_\bH$ yields $Y^\bblambda \cong Y^\bbmu$.
%By Theorem \ref{thm:Young}(b), $S_\bblambda \cong S_\bbmu$. 
%It then follows that  $\bblambda = \bbmu$.
Suppose that $P(\bblambda)\cong P(\bbmu)$, and hence $L(\bblambda)\cong L(\bbmu)$.
By Proposition~\ref{prop:PDfil}(b), $\Delta(\bbmu)$ is a quotient of $P(\bblambda)$.
Recall from Proposition~\ref{prop:PDfil}(a) that
$P(\bblambda)$ has a $\Delta$-filtration with subquotients $\Delta(\bbnu)$ satisfying $\bbnu\geq\bblambda$. 
Hence, there is $\bbnu\geq\bblambda$ such that $\Hom_\bS(\Delta(\bbnu),\Delta(\bbmu))\neq0$.
It then follows from Theorem~\ref{thm:exactseq}(a) and Proposition~\ref{prop:gHuDS} that $\bblambda\leq\bbnu\leq\bbmu$.
By symmetry, $\bbmu\leq\bblambda$, and hence $\bblambda=\bbmu$.
\end{proof} 

%--------------------------------------------------------------------
\subsection{Wreath Schur Functors are 1-Covers}\label{sec:1coverThm}
%--------------------------------------------------------------------
In this section, our dual Specht theory is applied to the quantum wreath product $\bH = \cH_q(\Sigma_m)\wr \cH(d)$ where $d\geq 2$. 
It is shown that the wreath Schur functor $\cF_\bH$ constructed in Proposition \ref{prop:WSF} is indeed a 1-cover under reasonable conditions.
%--------------------------------------------------------------------
\begin{thm} \label{thm:1cover}
Let $B= \cH_q(\Sigma_m)$, and $\bH = B\wr \cH(d)$ be the generalized Hu algebra with $d\geq 2$.
\begin{enua}
\item    
Suppose that \eqref{eq:goodH} holds. 
If $n \geq md$, $\ell \geq 4$ and $d! \in K^\times$, then the wreath Schur functor  is a $1$-cover of $\bH$. 
\item
If $n \geq 2m$ and \eqref{eq:finalell} holds, the wreath Schur functor for the Hu algebra $ \cH_q(\Sigma_m) \wr \cH(2)$ is a 1-cover.
\end{enua}
\end{thm} 
%--------------------------------------------------------------------
\begin{proof}
Part (a) follows from Theorem \ref{thm:1-cover}.
Part (b) is the $d=2$ special case of part (a) since \eqref{eq:goodH} follows from Proposition \ref{prop:H1}.
\end{proof}
%--------------------------------------------------------------------
Suppose that $n\geq 2m$. In light of the Morita equivalence \eqref{eq:MorD}, consider the following Schur algebra of type D:
\[
S^\D_q(n,2m):= \bS(n,m \wr 2) \oplus \Big(\bigoplus_{i=0}^{m-1} S_q(n,i) \otimes S_q(n, 2m-i)\Big).
\]
Let $e^\D:=e^\bH\oplus\bigoplus_{i=0}^{m-1}(e_i\otimes e_{2m-i})$, 
where $e^\bH\in\bS(n,m\wr 2)$ is the idempotent in Proposition~\ref{prop:WSF} 
and $e_j\in S_q(n,j)$ is the idempotent with $e_jS_q(n,j)e_j\cong\cH_q(\Sigma_j)$.
Then $e^\D S^\D_q(n,2m)e^\D\cong\cA(m)\oplus\bigoplus_{i=0}^{m-1}\cH_q(\Sigma_i\times\Sigma_{2m-i})$, 
and one can define
\[
\cF_\D:=F_M^{-1}\circ(e^\D\cdot-)\colon S^\D_q(n,2m)\textup{-mod}\to\cH_q(W^\D_{2m})\textup{-mod},
\]
where $F_M$ is the Morita equivalence \eqref{eq:MorD}.
%--------------------------------------------------------------------
\begin{thm}\label{thm:SDHWC}
Suppose that $n\geq 2m$ and \eqref{eq:finalell} holds. Then $S^\D_q(n,2m)$-mod is a highest weight category with respect to the poset 
\[
\Omega\sqcup\bigsqcup_{i=0}^{m-1}(\Pi_i\times\Pi_{2m-i}), 
\]
in which elements of different components are incomparable, $\Omega$ carries the order in Definition~\ref{def:KPO}, and $\Pi_i\times\Pi_{2m-i}$ carries the product of the dominance orders. Moreover, $\cF_\D$ is a 1-cover of $\cH_q(W^\D_{2m})$.
\end{thm}
\begin{proof}
The module category of a finite direct sum of algebras is the product of the module categories, and (H1)--(H4) can be verified componentwise. 
For $\bS(n,m\wr 2)$, this is Theorem~\ref{thm:1cover}. 
For each $0\leq i<m$, the $q$-Schur algebras $S_q(n,i)$ and $S_q(n,2m-i)$ are  quasi-hereditary with respect to the dominance order since $n\geq 2m$. 
Set 
\[
\Delta(\alpha,\beta):=\Delta(\alpha)\otimes\Delta(\beta),
\quad
P(\alpha,\beta):=P(\alpha)\otimes P(\beta).
\]
Then, (H1)--(H3) follow from the K\"unneth formula, 
since the simple modules are $L(\alpha)\otimes L(\beta)$. 

For (H4), note that the kernel of $P(\alpha,\beta)\to\Delta(\alpha,\beta)$ is filtered by $\ker(P(\alpha)\to\Delta(\alpha))\otimes P(\beta)$ and $\Delta(\alpha)\otimes\ker(P(\beta)\to\Delta(\beta))$, whose $\Delta$-subquotients have labels strictly larger in the product order. 
Since $\ell\geq 4$, it follows from \cite{HN04} that $R^1\cG$ vanishes on the dual Specht modules of each type A factor, and hence on $S_\alpha\otimes S_\beta$ by the K\"unneth formula.
Therefore, by Theorem~\ref{thm:exactseq}, the tensor product of the two type A Schur functors is again a 1-cover,
and so is the Schur functor $\cF_\D$.
\end{proof}

%--------------------------------------------------------------------
\section{GGOR conjecture for type D}\label{sec:GGOR}
%--------------------------------------------------------------------
{
In this section, we will show that the diagram \eqref{cd:GGOR} commutes for type $\D_{2m}$, 
provided $n \geq 2m$ and \eqref{eq:finalell} holds. 
}
%---------------------------------------------------
\subsection{{Rational Cherednik Algebras}}\label{sec:RCA}
%---------------------------------------------------
Following \cite{GGOR03}, let $\bbH = \bbH(W^\D_{2m})$ be the rational Cherednik algebra of type $\D_{2m}$. 
That is, $\bbH$ is a quotient of $T(\fh \oplus \fh^* ) \rtimes W^\D_{2m}$ subject to certain relations,
where $\fh = \CC^{2m}$ is the standard reflection representation of $W^\D_{2m}$.
In particular, the commutator relation between $y \in \fh$ and $x \in \fh^*$ is
\begin{equation}
    [y,x] = \<y,x\> - c \sum_{s} \<\alpha_s, y\> \< x, \alpha_s^\vee\> s,
\end{equation}
where $s$ runs over all reflections, $\alpha_s$ and $\alpha^\vee_s$ are the corresponding root and coroot, respectively.
Moreover, in the notation of \cite[Remark 5.1]{Ro08}, $c$ corresponds to the parameters $h_{H,0}=0$ and $h_{H,1}=c$ for all reflecting hyperplanes $H$.
Note that $\bbH$ is the rational Cherednik algebra $\bbH_{1,c}$ in Etingof--Ginzburg's notation \cite{EG02}.

\subsection{Rouquier's Unicity Theorem}
Let us first recall the following theorem due to Rouquier:

\begin{prop}[{\cite[Theorem 4.49]{Ro08}}]\label{Ro449}
Let $k$ be a commutative Noetherian domain with field of fractions $\KK$.
Let $\bH_k$ be a $k$-algebra which is finitely generated and projective as a $k$-module.
Suppose that $\bH_\KK := \KK\otimes_k \bH_k$ is split semisimple.
For $i\in\{1,2\}$, let $F_i:C_i\to \bH_k$-mod be a 1-cover for some poset $(\Irr \hp \bH_\KK, \leq_i)$.

If $\leq_1$ is a refinement of $\leq_2$, then $C_1$ and $C_2$ are equivalent as highest weight categories inducing bijections $\Irr \hp (\KK C_1) \cong \Irr \hp \bH_\KK \cong \Irr \hp (\KK C_2)$,
where $\KK \cC$ means $(\KK\otimes_k A)$-mod if $\cC = A$-mod for some $k$-algebra $A$.
\end{prop}
%--------------------------------
\subsection{The \texorpdfstring{$\bc$}{c}-function Partial Order}
%--------------------------------
For the rest of the article, let $K=\bbC$, fix $c\in\bbC$, and let $q:=\exp(2\pi\sqrt{-1}\,c)$. 
If $c\notin\QQ$, then $\ell=\infty$ and hence both $\cO(\bbH(W^\D_{2m}))$ and $S^\D_q(n,2m)$-mod are semisimple.
In other words, Theorem~\ref{thm:GGORD} reduces to the matching of labels in Proposition~\ref{prop:cellmod}.
Hence, we may assume from now on that $c\in\QQ_{>0}$ (see \cite[Proposition 5.8 and Remark 6.9]{Ro08}).
Let
\[
k:=\bbC[[t]],
\quad
\bbK:=\textup{Frac}(k) =\bbC((t)),
\quad
c' := c+t,
\quad
\bbq :=\exp(2\pi\sqrt{-1} c')
=q \exp(2\pi\sqrt{-1}t)\in k^{\times},
\]
and let $\bH_k := \cH_{\bbq}(W^\D_{2m})$ be the corresponding Hecke algebra over $k$.
Note that $k$ is a complete discrete valuation ring and hence a commutative Noetherian domain, with $k/(t) = K$ and $2 \in k^\times$.
The algebra
$\bH_\KK$ is split semisimple since 
%$\bbq^N=q^N\exp(2\pi\sqrt{-1} Nt)\notin\bbC$ for all $N\neq 0$, the element 
$\bbq$ is not a root of unity in $\KK$ and 
$\sqrt{\bbq}=\sqrt{q}\exp(\pi\sqrt{-1} t)\in k$.

Moreover, by the Tits deformation theorem, $\Irr\hp \bH_\KK$ is in bijection with the following set:
\begin{equation}    \label{IrrHK}
\left\{
[\lambda, \mu] 
~\middle|~ 
\lambda \neq \mu, (\lambda,\mu) \in \bigsqcup\nolimits_{i=0}^{m} (\Pi_i\times \Pi_{2m-i}) 
\right\}
\sqcup
\big\{
(\lambda,+), (\lambda,-) ~|~ \lambda \in \Pi_m
\big\},
\end{equation}
where $[\lambda, \mu] = \{(\lambda,\mu), (\mu,\lambda)\}$ is the $\Sigma_2$-orbit containing the bipartition $(\lambda,\mu)$.

We begin with a description of the 1-cover on the Cherednik side.
Following \cite[Section 5.2]{Ro08},
let $\cO_k$ be the category of finitely generated $\bbH$-modules that are $\fh$-locally nilpotent with corresponding Knizhnik--Zamolodchikov functor $\KZ_k$. 
{By \cite[Theorem~2.19]{GGOR03} and \cite[Section~5.2]{Ro08}, $\cO_k$ is a highest weight category with respect to the partial order $\preceq$ on $\Irr\hp\bH_\KK$ generated by $E \prec E'$ whenever $\bc(E)-\bc(E')\in\ZZ_{>0}$, 
where $\bc(E)$ is the scalar by which $\sum_s c(1-s)$ acts on $E$, and $s$ runs over all reflections of $W^\D_{2m}$.}
The standard objects in $\cO_k$ are $\{ S(E) \mid E \in \Irr\hp\bH_\KK \}$. 
{
%---------------------------
\begin{prop}\label{prop:cfcnfacts}
Consider Lusztig's $\ba$-function (resp. $\bA$-function) $\Irr\hp\bH_\KK \to \ZZ$. 
That is,
$\ba(E)$ (resp. $\bA(E)$) is the lowest (resp. highest) power of the indeterminate in the generic degree of $E\in \Irr\hp\bH_\KK$ (see \cite[Section 4.1.1]{Lu84}). Then
\begin{enua}
    \item $\bc(E) = c (\ba(E) + \bA(E))$.
    \item $\bA(E) = 2m(2m-1) - \ba(E\otimes \textup{det})$.
\end{enua}
\end{prop}
%---------------------------
\begin{proof}
    Part (a) follows from \cite[Section 4.21 and Proposition 4.1]{BM97}.
    Part (b) follows from \cite[Section 5.11.5]{Lu84}.
    See also \cite[Section 6.2]{GGOR03}.
    Note that the scalar $k_1$ therein is the scalar $c \in \CC$ in our normalization (see Section~\ref{sec:RCA}).
\end{proof}
%---------------------------
Then, we can combine Proposition \ref{prop:cfcnfacts} and the computations in \cite{Ge00} to give the desired partial order $\leq_1$ as follows.
}
{
\begin{Def}[$\bc$-function order $\leq_1$ on $\Irr\hp\bH_\KK$]
\label{def:PO1}
Recall from \cite[Section 6.2]{Ge00} that the 
$\mathbf{a}$-value of a bipartition $(\lambda,\mu)$ is given by
\begin{equation}\label{def:afcn}
\ba(\lambda,\mu) = \frac{-r(r-1)(2r-1)}{6}
+ \sum_{i=1}^r (i-1) (\lambda_i+\mu_i)
+ \sum_{i,j=1}^r \min\{\lambda_i+r-i, \mu_j+r-j\},
\end{equation}
where $r \geq 1$ is chosen such that both $\lambda$ and $\mu$ have at most $r$ parts.
{
Set
\begin{equation}\label{def:Afcn}
\bA(\lambda,\mu) := 2m(2m-1) - \ba(\lambda^t,\mu^t),
\end{equation}
where $2m(2m-1)$ is the total number of reflections in $W^\D_{2m}$. 
}
Define
\begin{equation}\label{leq1}
E \leq_1 E' \quad\Leftrightarrow\quad E \cong E' 
\quad
\textup{or}
\quad
(\ba+\bA)(E) > (\ba+\bA)(E'),
\end{equation}
where $\ba([\lambda, \mu]) := \ba(\lambda,\mu)$, $\ba(\lambda, \pm) := \ba(\lambda, \lambda)$, and similarly for $\bA$.
\end{Def}
%---------------------------
{
Note that if $E$ is labeled by $[\lambda,\mu]$ (resp.\ $(\lambda,\pm)$),
then $E\otimes\det$ is labeled by $[\lambda^t,\mu^t]$ (resp.\ $(\lambda^t,\pm)$ for a suitable sign).
It then follows from Proposition~\ref{prop:cfcnfacts}(b) that $\bA(E)$ agrees with \eqref{def:Afcn}.

Since $c>0$, the order $\leq_1$ refines the order given by $E \prec E'$.
Hence, $\cO_k$ is a highest weight category with respect to $\leq_1$.
Note that $\leq_1$ coincides with the $\ba$-function order for $m\leq 2$.
}

The following is a consequence of \cite[Proposition 4.6]{CGG12}, when restricted to type D.
%---------------------------
\begin{prop}\label{prop:GGOR68}
The image of the set of standard modules in $\cO_k$ is given by 
\[
\{C(E) ~|~ E \in \Irr\hp\bH_\KK\},
\]
where $C(E)$ is the cell module
%$S(E) := \varpi^*(R\otimes_{\bf{Q}} E)$ is the Kazhdan-Lusztig cell module 
corresponding to $E$ constructed in \cite[Example 4.4]{Ge07}.
\end{prop}
%---------------------------
\begin{proof}
In order to apply \cite[Proposition 4.6]{CGG12}, one needs to check the following two assumptions. The first one given in \cite[Section 4.2]{CGG12} was checked for type D in \cite[Section 4.3]{CGG12}.
The second one is the verification of Lusztig's conditions (P1)--(P15) \cite[Section 14.2]{Lu03}, which is known for type D.
\end{proof}
%---------------------------
}
%=======================
\subsection{The K\"unneth Order}
%=======================
{
Next, we  describe the image of the standard objects on the Schur algebra side. 
Recall the K\"unneth order $(\Omega, \leq_K)$ from   Definition~\ref{def:KPO}.
Also recall from Example \ref{ex:Huindex} the identifications
\[
e_\lambda^{(2)}\mapsto (\lambda, +),
\quad
e_\lambda^{(1,1)}\mapsto (\lambda, -),
\quad
e_\lambda^{(1)} + e_\mu^{(1)} \mapsto [\lambda, \mu],
\]
which induce a bijection between $\Omega$ and the following subset of $\Irr\hp\bH_\KK$:
\[
\Lambda^\D_m:=
\left\{
[\lambda, \mu] 
~\middle|~ 
\lambda \neq \mu, (\lambda,\mu) \in  \Pi_m\times \Pi_{m} 
\right\}
\sqcup
\{(\lambda,+), (\lambda,-) ~|~ \lambda \in \Pi_m\}.
\]
Hence, $(\Omega, \leq_K)$ induces a partial order $(\Irr\hp\bH_\KK, \leq_2)$ that we call the modified dominance order as follows.
%-----------------------------------
\begin{Def}[Modified dominance order $\leq_2$ on $\Irr\hp\bH_\KK$]
\label{def:PO2}
Suppose that $\bblambda, \bbmu \in \Omega$ are mapped to $E, E' \in \Lambda^\D_m$ under the above bijection. Then, $E \leq_2 E'$ in $\Irr\hp\bH_\KK$ if and only if $\bblambda \leq_K \bbmu$ in $\Omega$.
For the other relations, we may assume $E = [\lambda,
\mu] \not\in \Lambda^\D_m$. Thus, we can write the rearranged bipartition as
\begin{equation}
    (\lambda^{(1)}(E), \lambda^{(2)}(E)) = 
    \begin{cases}
        (\lambda,\mu) &\tif \# \lambda < \#\mu;
        \\
        (\mu,\lambda) &\tif \# \lambda > \#\mu,
    \end{cases}
\end{equation}
where $\#\lambda$ means the total number of boxes in $\lambda$.
Then, one imposes the following relations using the usual dominance orders if $E, E' \not\in \Lambda^\D_m$:
\begin{equation}
\begin{split}
E \leq_2 E'
\quad
\Leftrightarrow
\quad
&
\lambda^{(1)}(E) \unlhd \lambda^{(1)}(E')  \in \Pi_i,
\quad
 \lambda^{(2)}(E) \unlhd \lambda^{(2)}(E')
 \in \Pi_{2m-i}
\\
&\textup{ for the same } 0\leq i < m.
\end{split}
\end{equation}
There are no other relations imposed.
\end{Def}
}
%-----------------------------------
\begin{exa}\label{ex:POD}
The following is the Hasse diagram for the set of bipartitions of 4 with respect to the dominance order (cf. Definition \ref{def:domPO}).
\begin{equation*}
\xymatrix{
(\varnothing, (4)) \ar[d]
&
((1), (3)) \ar[d] \ar[l]
&
((2), (2)) \ar[d] \ar@/{}^{3pc}/[dd]
&
((3), (1)) \ar[d] 
&
((4), \varnothing) \ar[d]
\\
(\varnothing, (3, 1)) \ar[d] 
&
((1),(2, 1)) \ar[d] \ar[l]
&
((1,1), (2)) \ar@/{}_{3pc}/[dd] \ar[lu]
&
((2, 1), (1)) \ar[d] \ar[lu]
&
((3, 1), \varnothing) \ar[d] \ar[lu]
\\
(\varnothing, (2, 2)) \ar[d]
&
((1), (1,1, 1)) \ar[ld]
&
((2), (1,1)) \ar[d] 
&
((1,1,1),(1)) \ar[lu]
&
((2, 2), \varnothing) \ar[d]
\\
(\varnothing, (2, 1, 1)) \ar[d]
&
&
((1,1), (1,1)) \ar[luu]
&
&
((2, 1, 1), \varnothing) \ar[d] \ar[luu]
\\
(\varnothing, (1, 1, 1, 1))
&
&
&
&
((1, 1, 1, 1), \varnothing) \ar[luu]
}
\end{equation*}
In other words, the dominance order induces the following partial order $\unlhd$ on $\Irr \hp \bH_\KK$.
We will see that the modified dominance order $\leq_2$ can be obtained from this partial order.
\begin{equation*}
\xymatrix{
[\varnothing, (4)] \ar[d]
&
[(1), (3)] \ar[d] \ar[l]
&
[(2), (2)]_+ \ar[d] 
\\
[\varnothing, (3, 1)] \ar[d] 
&
[(1),(2, 1)] \ar[d] \ar[l]
&
[(2), (2)]_- \ar[d] 
\\
[\varnothing, (2, 2)] \ar[d]
&
[(1), (1,1, 1)] \ar[ld]
&
[(2), (1,1)] \ar[d] \ar[luu]
\\
[\varnothing, (2, 1, 1)] \ar[d]
&
&
[(1,1), (1,1)]_+  \ar[d]
\\
[\varnothing, (1, 1, 1, 1)]
&
&
[(1,1), (1,1)]_- \ar[luuu]
}
\end{equation*}
First, for $(\lambda, \mu) \in \Pi_i \times \Pi_{2m-i}, (\lambda', \mu') \in \Pi_j \times \Pi_{2m-j}$ where $0\leq i,j \leq m-1$,
$[\lambda, \mu] \leq_2 [\lambda', \mu']$ if and only if $i=j$ and $[\lambda, \mu] \unlhd [\lambda', \mu']$.
Second, the only relations for $\Lambda^\D_m$ are exactly those given by $\leq_K$.
Therefore, we obtain the following partial order $\leq_2$ for the case $m=2$:
\begin{equation*}
\xymatrix{
[\varnothing, (4)] \ar[d]
&
[(1), (3)] \ar[d] 
&
[(2), (2)]_+ \ar@/{}^{3pc}/[dd] 
\\
[\varnothing, (3, 1)] \ar[d] 
&
[(1),(2, 1)] \ar[d] 
&
[(2), (2)]_- \ar[d] 
\\
[\varnothing, (2, 2)] \ar[d]
&
[(1), (1,1, 1)] 
&
[(2), (1,1)] \ar[d]  \ar@/{}^{3pc}/[dd]
\\
[\varnothing, (2, 1, 1)] \ar[d]
&
&
[(1,1), (1,1)]_+  
\\
[\varnothing, (1, 1, 1, 1)]
&
&
[(1,1), (1,1)]_- 
}
\end{equation*}
\end{exa}
%-----------------------------------
\subsection{Image of Schur Functor}
Recall from \eqref{eq:MorD} the Morita equivalence
\[
F_M:\cH_{q}(W^\D_{2m})\hp\Mod \to  
(\cH_q(\Sigma_m)\wr \cH(2))\hp\Mod
\oplus
\bigoplus\nolimits_{i=0}^{m-1}
\cH_q(\Sigma_i  \times \Sigma_{2m-i})\hp\Mod.
\]
The following key result describes the image of the Schur functor.
{
\begin{prop}\label{prop:cellmod}
Suppose that $n\geq 2m$ and \eqref{eq:finalell} holds. Then, under the Schur functor $\cF_\D$,
the image of the set of standard modules in $S^\D_q(n,2m)$-mod is given by 
$\{F_M^{-1}(N) ~|~ N \in J\}$ where
\[
\begin{split}
J:=&\{S_\lambda \wr \triv, S_\lambda \wr \sgn ~|~ \lambda \in \Pi_m\}
\sqcup 
\{\ind_{(1,1)}^2(S_{\lambda} \otimes S_{\mu})~|~ \lambda, \mu \in \Pi_m, \ \lambda \neq \mu\}
\\
&~\sqcup
\bigsqcup_{i=0}^{m-1}
\{
S_{\lambda} \otimes S_{\mu} ~|~ 
(\lambda, \mu) \in \Pi_i\times \Pi_{2m-i}
\}.
\end{split}
\]
Moreover, the set 
$\{F_M^{-1}(N) ~|~ N \in J\}$ agrees with the set of cell modules considered in \cite[Example 4.4]{Ge07}, in the semisimple case.
\end{prop}
%-----------------------------------
\begin{proof}
The first assertion follows immediately from the construction.
Recall that $\Irr\hp\bH_\KK$ is in bijection  with the following set: 
\begin{equation}    \label{eq:CEHu}
\{
(\lambda, +), (\lambda, -) ~|~ \lambda \in \Pi_m \} 
\sqcup
\{
[\mu',\mu''] \in \Pi_m^2 ~|~ \mu' \neq \mu''
\}
\sqcup
\bigsqcup_{i=0}^{m-1}
\{
(\mu',\mu'') \in \Pi_i \times \Pi_{2m-i}
\}.
\end{equation}
Denote by $S^{(\lambda,\pm)}_\D$, $S^{[\mu', \mu'']}_\D$, and $S^{(\mu', \mu'')}_\D$ the corresponding cell modules of $\bH_\KK$, respectively.
It follows from \cite[Section 4]{Hu02} that $F_M$ sends these cell modules to the following modules:
\[
\begin{split}
&
F_M(S^{(\lambda,+)}_\D)
%=S(\lambda,\lambda)_+ 
= S_\lambda \wr \triv
= S_{(e_\lambda^{(2)})},
\quad
F_M(S^{(\lambda,-)}_\D)
%=S(\lambda,\lambda)_- 
= S_\lambda \wr \sgn
= S_{(e_\lambda^{(1,1)})},
\\
&
F_M(S^{[\mu',\mu'']}_\D) 
= \ind_{(1,1)}^{2}(S_{\mu'} \otimes S_{\mu''}) 
\cong  \ind_{(1,1)}^{2}(S_{\mu''} \otimes S_{\mu'}) ,
\quad
F_M(S^{(\mu',\mu'')}_\D) = S_{\mu'} \otimes S_{\mu''}.
\end{split}
\]
In our convention, dual Specht modules are those with simple heads that form the complete set of all simple modules.
Therefore, the latter assertion is proved.
\end{proof}
}
%-----------------------------------
\subsection{Schur-concavity}
%-----------------------------------
{
Let us collect some basic definitions and properties of the Schur-concave functions (see \cite{MO79}, for example). 
We will analyze the Schur-concavity of Lusztig's $\ba$-function {and deduce that of the $\bc$-function}, which will be used to prove the assumption needed in order to apply Rouquier's unicity theorem.
\begin{Def}
Consider a function $f : \Pi_n \to \mathbb{R}$.
The function $f$ is called {\em Schur-concave}  if $f(\lambda) \leq f(\lambda')$ whenever $\lambda \unrhd \lambda'$. 
The function $f$ is called {\em strictly Schur-concave}  if $f(\lambda) < f(\lambda')$ whenever $\lambda \rhd \lambda'$.     
\end{Def}
\begin{prop} \label{prop:schurconcave}
Suppose that $g : \mathbb{N} \to \mathbb{R}$ is a concave function.
Then, the function defined below is Schur-concave:
\[
\Pi_n \to \mathbb{R},
\quad\lambda \mapsto \sum_i g(\lambda_i).
\]
\end{prop}
\begin{proof}
We treat the function $\Pi_n \to \mathbb{R}$ as a function $\mathbb{R}^r \to \mathbb{R}$ for some sufficiently large $r$, and the result follows from \cite[Proposition C.1]{MO79}. 
\end{proof}
\begin{prop}\label{prop:SchurConc}
    Let $n, n' \in \NN$. 
    The function $\Pi_n\times \Pi_{n'} \to \RR$, $(\lambda,\mu) \mapsto \ba(\lambda,\mu)$ is strictly Schur-concave in both variables.
\end{prop}
\begin{proof}
Observe that we may rewrite \eqref{def:afcn} as
\begin{equation}\label{eq:alm}
    \ba(\lambda, \mu) = C %+ m (n + n') 
    - \sum_{i = 1}^r \left(\sum_{j = 1}^i \lambda_j\right) 
    - \sum_{i = 1}^r \left(\sum_{j = 1}^i \mu_j\right)  
    +\sum_{j = 1}^r \left( \sum_{i = 1}^r \min ( \lambda_i + r - i, \mu_j + r - j )  \right),
\end{equation}
where $C = - \frac{r (r - 1)(2 r - 1)}{6} + r (n + n')$ is a constant.
Here, the second and third summations are Schur-concave (in both $\lambda$ and $\mu$) by the definition of the dominance order. 

For each $j$, 
set $x_i := \lambda_i + r - i$ for $1 \leq i \leq r$.
A standard majorization argument (see \cite[3.C.1]{MO79}) allows us to assume that the summand is of the form $\sum_i g(x_i)$ with $g(x) = \min(x, \mu_j + r - j)$. 
Hence $\lambda \mapsto \sum_{i=1}^r \min(\lambda_i + r - i, \mu_j + r - j)$ is Schur-concave by Proposition~\ref{prop:schurconcave}.
Therefore, the last summand of the $\ba$-function is Schur-concave, and so is the function $\ba(-,\mu):\Pi_n \to \RR$. By a similar argument, $\ba(\lambda,-):\Pi_{n'} \to \RR$ is Schur-concave as well.     

Moreover, the Schur-concavity is strict. 
If $\lambda \rhd \lambda'$, then the second summation in \eqref{eq:alm} strictly decreases while the remaining terms are weakly Schur-concave. Hence $\ba(\lambda,\mu) < \ba(\lambda',\mu)$. 
A similar argument holds in the second variable.
\end{proof}
}
{
\begin{cor}\label{cor:aASchurConc}
    Let $n, n' \in \NN$. 
    The function $\Pi_n\times \Pi_{n'} \to \RR$, $(\lambda,\mu)\mapsto(\ba+\bA)(\lambda,\mu)$ is symmetric and strictly Schur-concave in both variables.
\end{cor}
\begin{proof}
Symmetry follows from that of $\ba$. 
Transposition reverses the dominance order, 
so if $\lambda\rhd\lambda'$ then $\lambda'^t\rhd\lambda^t$. 
Proposition~\ref{prop:SchurConc} then gives $\ba(\lambda'^t,\mu^t)<\ba(\lambda^t,\mu^t)$, i.e., $\bA(\lambda,\mu)<\bA(\lambda',\mu)$. 
Hence $\bA$, and therefore $\ba+\bA$, is strictly Schur-concave in the first variable. 
A similar argument holds in the second variable.
\end{proof}
}
%-----------------------------------
\subsection{Matching the 1-covers}
%-----------------------------------
{
One can use Corollary~\ref{cor:aASchurConc} above to verify the refinement assumption in Proposition \ref{Ro449}.
%----------------------------
\begin{prop}\label{prop:refinePO}
The $\bc$-function order $\leq_1$ (see Definition \ref{def:PO1})  is a refinement of our modified dominance order $\leq_2$ (see Definition \ref{def:PO2}).
\end{prop}
%----------------------------
\begin{proof}
{
It suffices to consider two elements related by either the dominance order or the K\"unneth order.
For the former, let $E,E'\notin\Lambda^\D_m$ with $E<_2E'$. Then $\lambda^{(1)}(E)\unlhd\lambda^{(1)}(E')$ and $\lambda^{(2)}(E)\unlhd\lambda^{(2)}(E')$ with at least one dominance strict, so Corollary~\ref{cor:aASchurConc}, applied in each variable separately, yields $(\ba+\bA)(E)>(\ba+\bA)(E')$.

For the K\"unneth order, 
it suffices to consider $E,E'\in\Lambda^\D_m$ corresponding to $\bblambda,\bbmu\in\Omega$ with $\bblambda\uparrow_K\bbmu$ and $\bblambda\neq\bbmu$. 
Write $\overline{\bblambda}=(\lambda^{(1)},\lambda^{(2)})$ and $\overline{\bbmu}=(\mu^{(1)},\mu^{(2)})$. 
Then $(\ba+\bA)(E)=(\ba+\bA)(\lambda^{(1)},\lambda^{(2)})$ and $(\ba+\bA)(E')=(\ba+\bA)(\mu^{(1)},\mu^{(2)})$.
%, including the case $E=(\lambda,\pm)$, where $\overline{\bblambda}=(\lambda,\lambda)$. 
Since $\ba+\bA$ is symmetric, one may assume $w=1$ in (K1), so that $\lambda^{(1)}\unlhd\mu^{(1)}$ and $\lambda^{(2)}\unlhd\mu^{(2)}$. 
By Corollary~\ref{cor:aASchurConc},
\[
(\ba+\bA)(\lambda^{(1)},\lambda^{(2)})\geq(\ba+\bA)(\mu^{(1)},\lambda^{(2)})\geq(\ba+\bA)(\mu^{(1)},\mu^{(2)}),
\]
and equality holds throughout only if $\lambda^{(1)}=\mu^{(1)}$ and $\lambda^{(2)}=\mu^{(2)}$. 
In that case $\overline{\bblambda}=\overline{\bbmu}$, 
so 
\[
n_\nu(1)=\#\bblambda(\nu)=\#\bbmu(\nu)
\]
for all $\nu\in\Pi_m$.
Thus, (K2) forces $\bblambda=\bbmu$, a contradiction.
Hence, $(\ba+\bA)(E)>(\ba+\bA)(E')$, i.e., $E<_1E'$. 
Therefore, $\leq_1$ is a refinement of $\leq_2$.
}
\end{proof}
}
Let $\bT(n,r)_k:= (\bigoplus_{i=1}^n k v_i)^{\otimes r}$.
Denote by $S_{\bbq}(n,r) := \End_{\cH_{\bbq}(\Sigma_{r})}(\bT(n,r)_k)$ the Schur algebra over $k = \CC[[t]]$.
Let $\mathcal{A}_k$ be the $k$-subalgebra of $\cH_{\bbq}(\Sigma_{2m})$ generated by $\cB_k:=\cH_{\bbq}(\Sigma_m\times\Sigma_m)$ and the element $\Hh$ (see \eqref{eq:preciseH1}). 
Let 
\[
S^\D_k(n,2m):=\End_{\mathcal{A}_k}(\bT_k)\oplus\bigoplus\nolimits_{i=0}^{m-1}
S_{\bbq}(n,i)\otimes S_{\bbq}(n,2m-i).
\]
By a formal argument, one can obtain the following 1-cover over $k$ on the Schur side.
%----------------------
\begin{lem}\label{lem:Schurk}
Suppose that $n\geq 2m$ and \eqref{eq:finalell} holds.
Then $K\otimes_kS^\D_k(n,2m)\cong S^\D_q(n,2m)$. Moreover, $(S^\D_k(n,2m)\textup{-mod},\cF_{\D,k})$ is a $1$-cover of $\bH_k$ with respect to $\leq_2$, where $\cF_{\D,k}$ is defined as $\cF_\D$.
\end{lem}
\begin{proof}
Since the identities in Proposition~\ref{prop:H1} hold over $\ZZ[v,v^{-1}]$, and $f_{2m}(\bbq)\in k^\times$, $\mathcal{A}_k$ is free over $k$ of rank $2(m!)^2$ and specializes to the Hu algebra over $\KK$.
Moreover, the Morita equivalence \eqref{eq:MorD} holds over $k$.
Since $2$ and $Z(m)$ are invertible over $k$, 
the transfer map (see Lemma~\ref{lem:transfer}) is defined over $k$.
Therefore, one can follow the proof of Theorem~\ref{thm:SDHWC} to show that $S^\D_k(n,2m)$-mod is a highest weight category over $k$.
Finally, $\KK\otimes_kS^\D_k(n,2m)$ is split semisimple, 
and the specialization of $\cF_{\D,k}$ at $t=0$ is the $1$-cover $\cF_\D$ of Theorem~\ref{thm:SDHWC}.
Additionally, $\KK\cF_{\D,k}$ matches the labels of the standard modules with those of $\Irr\hp\bH_\KK$ by Proposition~\ref{prop:cellmod}.
Hence, the assertion follows from \cite[Propositions 4.36 and 4.42]{Ro08}.
\end{proof}
%-------------------------

In summary, the GGOR problem is solved for the following choice of parameters for which 1-covers exist:
%-------------------------
\begin{thm}\label{thm:GGORD}
Let $c\in\QQ_{>0}$ or $c\notin\QQ$, and let $q=\exp(2\pi\sqrt{-1}c)$.
Suppose that $n\geq 2m$ and \eqref{eq:finalell} holds.
Then, the following diagram commutes:
\begin{equation}\label{cd:GGORD}
\begin{tikzcd}
\cO(\mathbb{H}(W^\D_{2m})) \arrow{rd}[swap]{\KZ_{\D}} \arrow{rr}{\simeq} &{}& S_q^\D(n,2m)\textup{-mod} \arrow{ld}{\mathcal{F}_\D}
\\
{}&\cH_{q}(W^\D_{2m})\textup{-mod}
\end{tikzcd}
\end{equation}
Moreover, $(\cO(\bbH), \KZ_\D)$ and $(S_q^\D(n,2m)\textup{-mod}, \cF_\D)$ are equivalent highest weight covers.
\end{thm}
%-------------------------
\begin{proof}
{
It follows from Lemma~\ref{lem:Schurk} that the Schur functor $\mathcal{F}_{\D,k}$ is a 1-cover provided $n\geq 2m$ and \eqref{eq:finalell} holds.
For the KZ functor $\KZ_{\D,k}: \cO_k(\bbH)\to \bH_k$-mod over $k$, it follows from \cite[Lemma 2.8]{RSVV16} that $\KZ_{\D,k}$ is a 1-cover if and only if 
\begin{equation}\label{eq:KZ1a}
\textup{lcd}(L) >2
\quad
\textup{for any simple }
L \in \cO(\bbH)_k
\textup{ with }\KZ_{\D,k}(L) = 0,    
\end{equation}
where
$\textup{lcd}(L):= \min\{i \in \ZZ_{\geq 0} ~|~ \Ext^i_{\cO(\bbH)_k}(L, T) \neq 0 \textup{ for some tilting object }T\in\cO(\bbH)_k\}$.
Denote by $\textup{Ch}(L)$ the characteristic variety (see \cite[Section 4.3.4]{GGOR03}) for any $L \in \cO(\bbH)_k$.
By \cite[Lemma 6.2]{RSVV16}, \eqref{eq:KZ1a} is equivalent to 
\begin{equation}\label{eq:KZ1b}
\dim \textup{Ch}(L) < 2m-2
\quad\textup{for any simple }
L \in \cO(\bbH)_k
\textup{ satisfying }
\dim\textup{Ch}(L) < 2m.     
\end{equation}
Moreover, \eqref{eq:KZ1b} is equivalent to the statement that all parabolic subalgebras of $\cH_q(W^\D_{2m})$ of rank at most two are semisimple, which holds when $\ell \geq 4$.

By Proposition~\ref{prop:refinePO}, 
$\leq_1$ is a refinement of $\leq_2$, and so Proposition~\ref{Ro449} yields an equivalence $\cO_k(\bbH) \cong S^\D_k(n,2m)$-mod of highest weight categories.
Therefore, we obtain the desired equivalence $\cO(\bbH)\cong S^\D_q(n,2m)$-mod of highest weight categories, which proves the theorem.
}
\end{proof}

%-------------------------
\appendix

%=============================================
\section{Technical proofs} \label{sec:app}
%=============================================
\subsection{Proof of Theorem \texorpdfstring{\ref{thm:sharp}}{4.2.1}} \label{app:equiv}
Since $\sigma^2=1$, 
\begin{equation}\label{eq:Hsigtt}
H \sigma(\theta) = \sigma^2(\theta) H = \theta H.
\end{equation}   
First, we show that the condition $(H^2 - S H - R)^\# = 0$ is equivalent to \eqref{eq:sharpB}:
\begin{equation}\begin{split}
(H^2 - S H - R)^\# &=  (H^\#)^2 - S^\# H^\# - R^\# 
\\
&=  (H - S){\theta (H-S) \theta} + S^\# (H - S){\theta} - R^\# 
\\
&=  (H - S) (H {\sigma(\theta)}- {\theta} S) {\theta} + S^\# (H - S){\theta} - R^\# 
\\
&=  R{\sigma(\theta)\theta} - (H{\theta} S - S^2 {\theta}) {\theta} + S^\# (H - S){\theta} - R^\# 
\\
&=  
(R{\sigma(\theta)} - \sigma({\theta}S) H - {\sigma(\theta)}\rho(S) + S^2{\theta}){\theta} 
+ S^\# (H - S){\theta} - R^\# 
\\
&=  (S^\# - \sigma(S{\sigma(\theta)})) (H - S){\theta} 
+ {\sigma(\theta)}(S^2 + R - \rho(S) - \sigma(S) S){\theta} 
- R^\#
\\
&=  (S^\# - {S\theta}) (H - S){\theta} 
+ \sigma(R ) {\theta \sigma(\theta)} - R^\#,
\end{split}
\end{equation}
where the third equality follows from \eqref{eq:Hsigtt},
the fourth equality follows from the paraphrase $R = (H-S)H$ of the quadratic relation, 
the second-to-last equality follows from the assumption that $S \theta = \theta S = \sigma(\theta)S$,
and the
last equality follows from \cite[(P2, P4)]{LNX25} and the assumptions $S = \sigma(S)$ and $R = \sigma(R)$, since
\begin{equation}
\sigma(\theta) R \theta 
= \sigma(\theta \sigma(R)) \theta 
= \sigma(\sigma(R) \theta ) \theta 
= R \sigma(\theta ) \theta.
\end{equation}
Next, we show that the condition $(H b - \sigma(b) H - \rho(b))^\# = 0$ is equivalent to \eqref{eq:sharpW}:
\begin{equation}
\begin{split}
(H b - \sigma(b) H - \rho(b))^\# 
&=  H^\# b^\# - \sigma(b)^\# H^\# - \rho(b)^\# \\
&=  (- H + S){\theta} b^\# - \sigma(b)^\# (- H + S){\theta} - \rho(b)^\# \\
&=  (- \sigma(b^\#) H + S b^\# - \rho(b^\#) + \sigma(b)^\# H - \sigma(b)^\# S){\theta} - \rho(b)^\# \\
&=  (\sigma(b)^\# - \sigma(b^\#)) H {\theta} 
+ (S b^\# - \sigma(b)^\# S - \rho(b^\#)){\theta} - \rho(b)^\#.
\end{split}
\end{equation}
Note that any non-adjacent braid relation
\eqref{equivrel3} is automatic due to the construction of $\#$.

For the adjacent braid relations, we need the following auxiliary results.
\begin{lem}\label{lem:sigSS}~
\begin{enua}
\item
$
\sigma_i(S_j) S_i + \rho_i(S_j) 
= S_i S_j.
$
\item 
$\sigma_i(S_j S_i) 
%= \sigma_i\sigma_j(S_i) \sigma_i(S_j) + \sigma_i\rho_j(S_i) 
= S_j \sigma_i(S_j) + \sigma_i\rho_j(S_i).$
\item 
$
\rho_i(S_j S_i) 
%= \rho_i(\sigma_j(S_i)S_j + \rho_j(S_i)) = \rho_i\sigma_j(S_i)S_j + \sigma_i\sigma_j(S_i) \rho_i(S_j) + \rho_i\rho_j(S_i) 
= S_j \rho_i(S_j) + \rho_i\rho_j(S_i).
$
\end{enua}
\end{lem}
\begin{proof}
Part (a) follows from the fact that the left-hand side is equal to 
$\sigma_i\sigma_i(\sigma_j(S_i)) S_i + \rho_i(S_j) 
= S_i \sigma_i\sigma_j(S_i) - \rho_i\sigma_i\sigma_j(S_i) - \sigma_i\rho_i\sigma_j(S_i) + \rho_i(S_j)$, which is equal to the right-hand side. 

Parts (b) and (c) follow from applying $\sigma_i$ and $\rho_i$ to both sides of part (a), respectively, since
\[
\sigma_i\sigma_j(S_i) \sigma_i(S_j) + \sigma_i\rho_j(S_i) 
= S_j \sigma_i(S_j) + \sigma_i\rho_j(S_i),
\]
and
\[
\rho_i(\sigma_j(S_i)S_j + \rho_j(S_i)) 
= \rho_i\sigma_j(S_i)S_j + \sigma_i\sigma_j(S_i) \rho_i(S_j) + \rho_i\rho_j(S_i) 
= S_j \rho_i(S_j) + \rho_i\rho_j(S_i).
\qedhere
\]
\end{proof}
Finally, we are in a position to show that the condition $(H_i H_j H_i - H_j H_i H_j)^\#=0$ is equivalent to \eqref{eq:sharptheta}--\eqref{eq:sharpBB}.
Note that it follows from the assumption $S_i \theta_j = \sigma_i(\theta_j) S_i$ that
\[
\begin{split}
&(H_i H_j H_i)^\# 
=  (S_i - H_i){\theta_i} (S_j - H_j){\theta_j} ( S_i - H_i) {\theta_i}
\\
&\quad 
=  (S_i S_j S_i - H_i S_j S_i - S_i S_j H_i - S_iH_j S_i + H_iS_jH_i + H_iH_j S_i + S_i H_j H_i - H_iH_jH_i)
{\sigma_i\sigma_j(\theta_i) \sigma_i(\theta_j)\theta_i},
\end{split}
\]
in which each term can be expressed in terms of the PBW basis as follows:
\begin{align}
H_i S_j S_i&= \sigma_i(S_j S_i) H_i + \rho_i(S_jS_i),
\\
S_iH_j S_i&= S_i \sigma_j(S_i) H_j + S_i \rho_j(S_i),
\\
H_iS_jH_i &= (\sigma_i(S_j)S_i + \rho_i(S_j) )H_i + \sigma_i(S_j) R_i,
\\
\label{HHS}H_iH_j S_i &= \sigma_i \sigma_j(S_i) H_i H_j + \rho_i \sigma_j(S_i) H_j + \sigma_i\rho_j(S_i) H_i + \rho_i \rho_j(S_i)
\\
\notag&= S_j H_i H_j + \sigma_i\rho_j(S_i) H_i + \rho_i \rho_j(S_i),
\end{align}
where the second equality of \eqref{HHS} follows from \cite[(P8)]{LNX25}.
Therefore, the coefficient for $H_i$ in $(H_iH_jH_i)^\#$ can be simplified as
\begin{equation}
\sigma_i(S_j) S_i + \rho_i(S_j) + \sigma_i\rho_j(S_i) - \sigma_i(S_j S_i) - S_i S_j
= - S_j \sigma_i(S_j)
\end{equation}
multiplied by ${\sigma_i\sigma_j(\theta_i) \sigma_i(\theta_j)\theta_i}$,
since $\sigma_i(S_j) S_i + \rho_i(S_j) - S_i S_j = 0 = S_j \sigma_i(S_j) + \sigma_i\rho_j(S_i) - \sigma_i(S_j S_i)$, by Lemma \ref{lem:sigSS} (a) and (b). 
Moreover, the coefficient for $H_{\textup{id}}$ in $(H_iH_jH_i)^\#$ can be simplified as
\begin{equation}
\rho_i\rho_j(S_i) - \rho_i(S_j S_i) - S_i \rho_j(S_i) + S_i S_j S_i + \sigma_i(S_j) R_i 
= - S_j \rho_i(S_j) + S_i S_j S_i - S_i \rho_j(S_i)  + \sigma_i(S_j) R_i
\end{equation}
multiplied by ${\sigma_i\sigma_j(\theta_i) \sigma_i(\theta_j)\theta_i}$,
since $S_j \rho_i(S_j) + \rho_i\rho_j(S_i) - \rho_i(S_j S_i) =0$, by Lemma \ref{lem:sigSS} (c). Thus,
\begin{equation}
\begin{split}
(H_iH_jH_i)^\# 
& 
= \big(S_i S_j S_i - S_j \rho_i(S_j)  - S_i \rho_j(S_i) + \sigma_i(S_j) R_i
- S_j \sigma_i(S_j) H_i - S_i \sigma_j(S_i) H_j
\\
&
  + S_j H_i H_j + S_i H_j H_i - H_i H_j H_i \big){\sigma_i\sigma_j(\theta_i) \sigma_i(\theta_j)\theta_i}.
\end{split}
\end{equation}
Therefore, $(H_i H_j H_i - H_j H_i H_j)^\# = 0$ if and only if 
{
$\sigma_i\sigma_j(\theta_i) \sigma_i(\theta_j)\theta_i
=\sigma_j\sigma_i(\theta_j) \sigma_j(\theta_i)\theta_j
$ and 
\[
S_i S_j S_i + \sigma_i(S_j) R_i - S_j S_i S_j - \sigma_j(S_i) R_j = 0.
\]
}
%\begin{equation}
%\begin{split}
%&(H_i H_j H_i - H_j H_i H_j)^\# 
%=  (H_i H_j H_i)^\# - (H_j H_i H_j)^\# 
%\\
%&\quad =  (H_j H_i H_j - H_i H_j H_i) + (S_i S_j S_i + \sigma_i(S_j) R_i) - %S_j S_i S_j - \sigma_j(S_i) R_j.
%\end{split}
%\end{equation}
The statement of the theorem is now verified.
%===================================
\subsection{Proof of Lemma \texorpdfstring{\ref{lem:Msharpground}}{4.3.1}} \label{app:Msharpground}
%===================================
Consider $m \in M^{\otimes d}$ and $b\in B^{\otimes d}$. 
Then, 
the first equation in \eqref{def:M2}\# holds since
\begin{equation}   
\sigma^{M^\#}_i(b \circ m) 
= \sigma^{M}_i(b^\# \cdot m)
= \sigma_i(b^\#)\cdot \sigma^{M}_i( m)
= \sigma_i(b)^\#\cdot \sigma^{M^\#}_i(m)
= \sigma_i(b)\circ \sigma^{M^\#}_i(m),
\end{equation}
where the second equality follows from \eqref{def:M2},
and the third equality follows from \eqref{eq:sharpW}.
Next, the second equation in \eqref{def:M2}\# holds since
\begin{equation}
\begin{split}    
\rho^{M^\#}_i(b \circ m) 
&=(l^M_{S_i{\theta_i}} -{l^M_{\theta_i}}\rho^{M}_i)(b^\# \cdot m)
\\
&= S_i {\theta_i} b^\# \cdot m - {\theta_i}\cdot\rho^M_i(b^\# \cdot m)
\\
&= \sigma_i(b^\#) S_i {\theta_i}\cdot m 
-\rho_i(b^\#){\theta_i}\cdot m
-{\theta_i}\sigma_i(b^\#)\cdot \rho_i^M(m)
%+ \sigma_i(b^\#)\cdot (-\rho^{M}_i(m)) - \rho_i(b^\#)\cdot m
\\
&= (\sigma_i(b^\#) S_i
-\rho_i(b^\#)) {\theta_i}\cdot m
-{\theta_i}\sigma_i(b)^\# \cdot \rho_i^M(m)
\\
&= \rho_i(b)^\# \cdot m 
+ 
\sigma_i(b)^\# \cdot ( S_i {\theta_i}  \cdot m 
- {\theta_i} \cdot \rho_i^M(m))
\\
&= \sigma_i(b)\circ \rho^{M^\#}_i(m) + \rho_i(b)\circ m,
\end{split}
\end{equation}
where the third equality follows from combining 
\eqref{def:M2}, the fact that $\theta_i \in Z(B^{\otimes d})$, 
and the assumption that $Sb = \sigma(b)S$,
and the fourth and fifth equalities follow from \eqref{eq:sharpW}.
Next, after dividing $\chi(\theta)$ (and $\chi(\theta)^2$, resp.) from both sides of the first (and the second, resp.) equality in \eqref{def:M4}\#, these equalities become 
\begin{align}
\label{eq:M4sharp1}&
\chi(S)(\sigma^{M}_i)^2
+(l_{S_i}^M-\rho_i^M)\sigma^{M}_i
+\sigma^{M}_i(l_{S_i}^M-\rho_i^M) 
= \chi(S)(\sigma^{M}_i), 
    \\
\label{eq:M4sharp2}    &
\chi(R) (\sigma^{M}_i)^2+(l_{S_i}^M-\rho_i^M)^2     
=l^{M}_{S_i} (l_{S_i}^M-\rho_i^M) + l^{M}_{R_i}.
\end{align}
Here, \eqref{eq:M4sharp1} rearranges into
\begin{equation}
    2\chi(S)(\sigma^{M}_i)^2
    =\chi(S)(\sigma^{M}_i)^2
+\rho_i^M\sigma^{M}_i
+\sigma^{M}_i\rho_i^M 
-\chi(S)(\sigma^{M}_i),
\end{equation}
whose right-hand side equals zero, thanks to \eqref{def:M4}.
Therefore, \eqref{eq:M4sharp1} holds whenever $\chi(S)=0$ (in particular when $S=0$) or when $\sigma^M= 0$.
Next, \eqref{eq:M4sharp2} rearranges into
\begin{equation}
    \chi(R) (\sigma^{M}_i)^2   
=\rho_i^M (l_{S_i}^M-\rho_i^M) + \chi(R)\id
= \chi(S)\rho_i^M -(\rho_i^M)^2 + \chi(R)\id,
\end{equation}
which is equivalent to the second equality in \eqref{def:M4}.

Regarding \eqref{def:M5}\#, the first equation holds since both sides are obtained from \eqref{def:M5} up to a unit.
The second equation in \eqref{def:M5}\# holds since, for $|i-j|=1$ and $m \in M^{\otimes d}$:
\begin{equation}
\sigma^M_j \sigma^M_i l^M_{S_j} (m)
= \sigma^M_j (\sigma_i(S_j) \cdot \sigma^M_i (m))
= (\sigma_j\sigma_i(S_j)) \cdot \sigma^M_j \sigma^M_i (m)
= l^M_{S_i}\sigma^M_j\sigma^M_i(m) ,        
\end{equation}
where the third equality uses the fact that $S_i = \sigma_j\sigma_i(S_j)$, which follows from \cite[(P8)]{LNX25}.
{For \eqref{def:M6}\#, in the first scenario when $S=0$, it simplifies to
\[
\chi(\theta)^2 \rho^M_i \sigma_j^M \rho^M_i
=
\chi(\theta)^2
(
\chi(S)
\sigma^M_j \rho_i^M \sigma^M_j
+\rho^M_j \rho_i^M \sigma^M_j
+\sigma^M_j \rho_i^M \rho^M_j
),
\]
which follows from \eqref{def:M6}.
In the second scenario when $\sigma^M = 0$, 
\eqref{def:M6}\# holds automatically.
}
%Observe that \eqref{def:M6}\# holds if and only if 
%\begin{equation}\label{eq:Msharp6}
%    (S_i \sigma_j(S_i) - S_jS_i) \cdot \sigma^M_j(m) = (1-\chi(S)) \sigma_j(S_i) (\sigma_j^M)^2(m),
%\end{equation}
%whose left-hand side is always zero by combining $Sb = \sigma(b) S$ and \cite[(P8)]{LNX25}.
%The right-hand side of \eqref{eq:Msharp6} is zero since $(1-\chi(S))\sigma^M(m) = 0$.
Finally, \eqref{def:M7}\# simplifies into
{
\begin{equation}\label{eq:M7sharp1} 
\chi(\theta)^3\chi(S)(
\chi(R) (\sigma_1^M)^2 + (\rho^M_1)^2 - \chi(S) \rho^M_1
)=
\chi(\theta)^3\chi(S)(
\chi(R) (\sigma_2^M)^2 + (\rho^M_2)^2 - \chi(S) \rho^M_2
)
\end{equation}
after first canceling the cubic terms using \eqref{def:M7}.
Applying \eqref{def:M4}, \eqref{eq:M7sharp1} becomes $\chi(\theta)^3\chi(S)(l^M_{R_1} - l^M_{R_2}) = 0$, which holds since $R_1$ and $R_2$ act by the same scalar $\chi(R)$.
}
The statement of the lemma is now verified.
%===================================
\subsection{Proof of Lemma \texorpdfstring{\ref{lem:Mstarground}}{4.5.1}} \label{app:Mstarground}
%===================================
{
Consider $f \in (M^\ast)^{\otimes d}$, $m \in M^{\otimes d}$ and $b\in B^{\otimes d}$. 
Then, the first equation in \eqref{def:M2}* holds since
\begin{equation}
\begin{split}
\sigma^{M^\ast}_i(b \circ f)
= (b\circ f)\sigma^M_i: m
&\mapsto f(b^\ast\cdot\sigma^M_i(m))
\\
&= f\sigma_i(\sigma_i(b)^\ast)\cdot\sigma^M_i(m)
\\
&= f\sigma^M_i(\sigma_i(b)^\ast \cdot m)
\\
&= (\sigma_i(b)\circ \sigma^{M^\ast}_i(f))(m),    
\end{split}
\end{equation}
where 
the second equality follows from \eqref{eq:starW}, and
the third equality follows from \eqref{def:M2}.
The second equation in \eqref{def:M2}* holds since
\begin{equation}
\begin{split}
\rho^{M^\ast}_i(b \circ f)
= (b\circ f)\rho^M_i: m
&\mapsto f(b^\ast\cdot\rho^M_i(m))
\\
&= f(\sigma_i(\sigma_i(b)^\ast))\cdot\rho^M_i(m)
\\
&= f\rho^M_i(\sigma_i(b)^\ast \cdot m) - f(\rho_i\sigma_i(b)^\ast\cdot m)
\\
&= \rho_i^{M^\ast}(f)(\sigma_i(b)^\ast\cdot m) + f(\rho_i(b)^\ast\cdot m)
\\
&=
(\sigma_i(b)\circ \rho_i^{M^\ast}(f) + \rho_i(b)\circ f)(m),
\end{split}
\end{equation}
where the second and fourth equalities follow from \eqref{eq:starW},
and 
the third equality follows from \eqref{def:M2}.
Next, for \eqref{def:M4}*, observe first that 
\[
\rho^{M^\ast} \sigma^{M^\ast} f (m) 
= (\sigma^{M^\ast} f)(\rho^M(m))
=  f(\sigma^{M}\rho^M(m)),
\]
i.e., the order of maps in $\Hom_K(M^{\otimes d}, K)$ is swapped upon passing to the dual. Hence,
\begin{equation}
\begin{split}
\big(
\chi(S) (\sigma^{M^\ast}_i)^2 + \rho^{M^\ast}_i\sigma^{M^\ast}_i + \sigma^{M^\ast}_i\rho^{M^\ast}_i
\big)(f) : m &\mapsto f(
\chi(S) (\sigma^{M}_i)^2(m)  + \sigma^{M}_i\rho^{M}_i(m) + \rho^{M}_i\sigma^{M}_i(m)
)
\\
&= f(l_{S_i}^M \sigma^M_i(m))
\\
&= \big(l_{S_i}^{M^\ast}\sigma^{M^\ast}_i(f)\big)(m),
\end{split}
\end{equation}
\begin{equation}
\begin{split}
\big(
\chi(R) (\sigma^{M^\ast}_i)^2+(\rho^{M^\ast}_i)^2
\big)(f): m &\mapsto f(
\chi(R) (\sigma^{M}_i)^2(m)+(\rho^{M}_i)^2(m)
)
\\
&= f(l^{M}_{S_i} \rho^{M}_i(m) + l^{M}_{R_i}(m))
\\
&= \big((l^{M^\ast}_{S_i} \rho^{M^\ast}_i + l^{M^\ast}_{R_i})(f)\big)(m),
\end{split}
\end{equation}
where we used the $K$-linearity of $f$ and \eqref{def:M4}.
The verification of \eqref{def:M5}*--\eqref{def:M7}* is similar to the proof of \eqref{def:M4}* but simpler, and hence is omitted.
The statement of the lemma is now verified.
}
%===================================

\section*{AI Disclosure}

During the preparation of this manuscript, 
the authors used \textsf{Claude} solely for the purposes of proofreading and grammatical correction of the human-written English text. 
The tool was neither employed to conceptualize mathematical ideas nor to develop proof strategies. 
All definitions, lemmas, propositions, theorems, and proofs presented in this work were exclusively conceived, developed, and verified by the authors. 
The authors critically reviewed and approved all textual modifications and assume full responsibility for the integrity, mathematical correctness, and final contents of the paper.

\end{document}